\documentclass[a4paper]{article}

\usepackage{geometry}
\usepackage[british]{babel}

\makeatletter
\RequirePackage{ifthen}

	\AtBeginDocument{\@ifpackageloaded{showkeys}
	{}
	{}
}
\makeatother

\newcommand\todo[2][a]{}

\usepackage{amsmath}
\usepackage{amssymb}
\usepackage{amsthm}
\usepackage[british]{babel}
\usepackage{bm}
\usepackage{dsfont}
\usepackage{enumerate}
\usepackage{epic}
\usepackage{eepic}
\usepackage[OT2,OT1]{fontenc}
\usepackage{framed}
\usepackage{ifthen}
\usepackage{pifont}
\usepackage{rotating}
\usepackage{stmaryrd}
\usepackage{textcomp}
\usepackage{theoremref}
\usepackage{upgreek}
\usepackage{url}
\usepackage{wasysym}
\input xy
\xyoption{all}

\newcommand{\mc}[1]{{\mathcal{#1}}}
\newcommand{\mf}[1]{{\mathfrak{#1}}}
\newcommand{\bb}[1]{{\mathbb{#1}}}

\newcommand\cyr{%
\renewcommand\rmdefault{wncyr}%
\renewcommand\sfdefault{wncyss}%
\renewcommand\encodingdefault{OT2}%
\normalfont
\selectfont}
\DeclareTextFontCommand{\textcyr}{\cyr}

\DeclareMathOperator{\IM}{Im}

\renewcommand{\Im}{\IM}

\newcommand{\qu}{\overline}
\newcommand{\defeq}{\mathrel{\mathop:}=}
\newcommand{\eqdef}{=\mathrel{\mathop:}}
\newcommand{\defequ}{\mathrel{\mathop:}\hspace*{-0.72ex}&=}
\newcommand{\defequiv}{\mathrel{\mathop:}\Longleftrightarrow}
\newcommand{\cls}{{\rm c.l.s.}}
\DeclareMathOperator{\dom}{dom}
\DeclareMathOperator{\id}{id}
\DeclareMathOperator{\ind}{ind}
\newcommand{\nontang}{\stackrel{\tiny\varangle}{\longrightarrow}}
\renewcommand{\oe}{{\ddot o}}
\DeclareMathOperator{\ran}{ran}

\DeclareMathOperator{\Rbv}{rbv}
\newcommand{\rbv}{\Rbv_z}
\newcommand{\rbvr}{\Rbv_{z,1}}
\newcommand{\rbvs}{\Rbv_{z,2}}
\newcommand{\rbvpr}[1]{\hspace*{2ex}\widehat{\rule{0ex}{1.8ex}}\hspace*{-2ex}\Rbv_{#1}}

\newcommand{\rbvspr}[1]{\hspace*{2ex}\widehat{\rule{0ex}{1.8ex}}\hspace*{-2ex}\Rbv_{#1,2}}
\DeclareMathOperator{\RbvSL}{rbv}
\newcommand{\rbvSL}[1]{\RbvSL_{#1}^{\rm SL}}
\newcommand{\rbvSLr}[1]{\RbvSL_{#1,1}^{\rm SL}}
\newcommand{\rbvSLs}[1]{\RbvSL_{#1,2}^{\rm SL}}

\newcommand{\rbvSLplr}[1]{\RbvSL_{#1,1}^{\rm SL,+}}
\newcommand{\rbvSLpls}[1]{\RbvSL_{#1,2}^{\rm SL,+}}
\newcommand{\rbvSLpr}[1]{\hspace*{2ex}\widehat{\rule{0ex}{1.8ex}}\hspace*{-2ex}\RbvSL_{#1}^{\rm SL}}
\newcommand{\rbvSLrpr}[1]{\hspace*{2ex}\widehat{\rule{0ex}{1.8ex}}\hspace*{-2ex}\RbvSL_{#1,1}^{\rm SL}}
\newcommand{\rbvSLspr}[1]{\hspace*{2ex}\widehat{\rule{0ex}{1.8ex}}\hspace*{-2ex}\RbvSL_{#1,2}^{\rm SL}}

\newcommand{\rbvSchr}[1]{\Rbv_{#1}^{\rm Schr}}
\newcommand{\rbvSchrr}[1]{\Rbv_{#1,1}^{\rm Schr}}
\newcommand{\rbvSchrs}[1]{\Rbv_{#1,2}^{\rm Schr}}
\newcommand{\rbvSchrpr}[1]{\hspace*{2ex}\widehat{\rule{0ex}{1.8ex}}\hspace*{-2ex}\RbvSL_{#1}^{\rm Schr}}
\newcommand{\rbvSchrrpr}[1]{\hspace*{2ex}\widehat{\rule{0ex}{1.8ex}}\hspace*{-2ex}\RbvSL_{#1,1}^{\rm Schr}}
\newcommand{\rbvSchrspr}[1]{\hspace*{2ex}\widehat{\rule{0ex}{1.8ex}}\hspace*{-2ex}\RbvSL_{#1,2}^{\rm Schr}}

\newcommand{\rbvDirac}[1]{\Rbv_{#1}^{\rm Dir}}
\newcommand{\NSL}[1]{\mathfrak{N}^{\rm SL}_{#1}}
\newcommand{\NSchr}[1]{\mathfrak{N}^{\rm Schr}_{#1}}
\newcommand{\pxxh}{p_{x_0,\hat x_0}}
\newcommand{\rd}{\mathrm{d}}
\newcommand{\reg}{{\rm reg}}
\DeclareMathOperator{\Res}{Res}
\DeclareMathOperator{\spn}{span}
\DeclareMathOperator{\supp}{supp}
\DeclareMathOperator{\tr}{tr}

\newlength{\maxlabwidth}

\numberwithin{equation}{section}
\swapnumbers
\theoremstyle{plain}
	\newtheorem{lemma}{Lemma}[section]
	\newtheorem{proposition}[lemma]{Proposition}
	\newtheorem{theorem}[lemma]{Theorem}
	\newtheorem{corollary}[lemma]{Corollary}
	\newtheorem{ntheoreM}[lemma]{}
\theoremstyle{definition}
	\newtheorem{definitioN}[lemma]{Definition}
	\newtheorem{assumptioN}[lemma]{Assumption}
\theoremstyle{remark}
	\newtheorem{remarK}[lemma]{Remark}
	\newtheorem{examplE}[lemma]{Example}
	\newtheorem{nremarK}[lemma]{}
\newcommand{\thlab}[1]{\thlabel{#1}\label{#1.}}
\renewcommand{\qedsymbol}{\raisebox{-2pt}{\large\ding{113}}}

\newenvironment{definition}{\begin{definitioN}}{
	\renewcommand{\qedsymbol}{\defendsymbol}
	\popQED{\qed}\renewcommand{\qedsymbol}{\qedsymbolsave}\end{definitioN}}

\newenvironment{remark}{\begin{remarK}}{
	\renewcommand{\qedsymbol}{\defendsymbol}
	\popQED{\qed}\renewcommand{\qedsymbol}{\qedsymbolsave}\end{remarK}}
\newenvironment{example}{\begin{examplE}}{
	\renewcommand{\qedsymbol}{\defendsymbol}
	\popQED{\qed}\renewcommand{\qedsymbol}{\qedsymbolsave}\end{examplE}}
\newenvironment{nremark}[1]{\begin{nremarK}\textit{\!{#1.}\, }}{
	\renewcommand{\qedsymbol}{\defendsymbol}
	\popQED{\qed}\renewcommand{\qedsymbol}{\qedsymbolsave}\end{nremarK}}
\newenvironment{proofof}[1]{\begin{proof}[\textit{Proof of #1}]}{\end{proof}}

\newcommand{\smmatrix}[4]{\Bigl(\begin{smallmatrix}
\hspace*{-0.2ex} #1 \hspace*{0.2ex} & \hspace*{0.2ex} #2 \hspace*{-0.2ex} \\[0.5ex]
\hspace*{-0.2ex} #3 \hspace*{0.2ex} & \hspace*{0.2ex} #4 \hspace*{-0.2ex}
\end{smallmatrix}\Bigr)}
\newcommand{\eps}{\varepsilon}

\newcommand{\AV}{A_{V}}
\newcommand{\rmo}{\mathrm{o}}
\newcommand{\rmO}{\mathrm{O}}
\newcommand{\wt}{\widetilde}
\newcommand{\wh}{\widehat}
\newcommand{\KSL}{\mathbb{K}_{\rm SL}}
\newcommand{\KSLpl}{\mathbb{K}_{\rm SL}^+}
\newcommand{\KSchr}{\mathbb{K}_{\rm Schr}}
\newcommand{\DeltaSL}{\Delta_{\rm SL}}
\newcommand{\DeltaSLpl}{\Delta_{\rm SL}^+}
\newcommand{\DeltaSchr}{\Delta_{\rm Schr}}

\newcommand{\void}[2]{\ifthenelse{#1>0}{#2}{}}

\newcounter{mytoc}

\begin{document}

\begin{flushleft}
	{\Large\textbf{%
	Direct and Inverse Spectral Theorems for a Class of \\[0.5ex]
	Canonical Systems with Two Singular Endpoints}
	}
	\\[5mm]
	\textsc{Matthias Langer\,\ $\ast$\,\ Harald Woracek}
	\\[6mm]
\end{flushleft}
	{\small
	\textbf{Abstract:}
	Part I of this paper deals with two-dimensional canonical systems $y'(x)=zJH(x)y(x)$,
	$x\in(a,b)$, whose Hamiltonian $H$ is non-negative and locally integrable,
	and where Weyl's limit point case takes place at both endpoints $a$ and $b$.
	We investigate a class of such systems defined by growth restrictions on $H$ towards $a$.
	We develop a direct and inverse spectral theory parallel to the theory of
	Weyl and de~Branges for systems in the limit circle case at $a$.
	Our approach proceeds via --- and is bound to --- Pontryagin space theory.
	It relies on spectral theory and operator models in such spaces,
and on the theory of de~Branges Pontryagin spaces.
	\\[1mm]
	The main results concerning the direct problem are:
	(1) showing existence of regularized boundary values at $a$;
	(2) construction of a singular Weyl coefficient and a scalar spectral measure;
	(3) construction of a Fourier transform and computation of its action and the
	action of its inverse as integral transforms.
	The main results for the inverse problem are: (4) characterization of the
	class of measures that are obtained via the above construction (positive Borel
	measures with power growth at $\pm\infty$); (5) a global uniqueness theorem
	(if Weyl functions or spectral measures coincide, Hamiltonians essentially coincide);
	(6) a local uniqueness theorem (if Weyl functions coincide up to an
	exponentially small error, Hamiltonians essentially coincide up to a certain point).
	\\[1mm]
	In Part II of the paper the results of Part I are applied to Sturm--Liouville equations
	with singular coefficients.
	We investigate classes of equations without potential (in particular,
	equations in impedance form) and Schr\"odinger equations, where
	coefficients are assumed to be singular but subject to growth restrictions.
	In the latter case potentials include Bessel-type potentials but also
	highly oscillatory potentials.
	We obtain corresponding direct and inverse spectral theorems.
	}
	\\
\begin{flushleft}
	{\small
	\textbf{AMS MSC 2010:}
	34B05, 34L40, 34B20, 34A55, 47B50, 47B32
	\\
	\textbf{Keywords:}
	canonical system, Sturm--Liouville equation, singular potential,
	direct and inverse spectral theorems, Pontryagin space, de~Branges space
	}
\end{flushleft}


%
%
\section{Introduction}
\label{sec-introduction}
%
%

By a \emph{Hamiltonian} we understand a function $H$ defined on a (possibly unbounded) interval $(a,b)$, which
takes real and non-negative $2\!\times\!2$-matrices as values, is locally integrable, and does not
vanish on any set of positive measure.
Throughout this paper we assume that Weyl's limit
point case prevails at the endpoint $b$; this means that for one
(and hence for all) $x_0\in(a,b)$, we have $\intop_{x_0}^b\tr H(x)\,\rd x=\infty$.

The \emph{canonical system} associated with $H$ is the differential equation
\begin{equation}\label{A30}
	y'(x) = zJH(x)y(x),\qquad x\in(a,b),
\end{equation}
where $z$ is a complex parameter (the eigenvalue parameter), $J$ is the signature matrix
$J \defeq \smmatrix{0}{-1}{1}{0}$ and $y$ is a $2$-vector-valued function.
Canonical systems appear frequently in natural sciences,
for example in Hamiltonian mechanics or as generalizations of Sturm--Liouville problems,
e.g.\ in the study of a vibrating string with non-homogeneous mass distribution.
They provide a unifying approach to Schr\"odinger operators, Jacobi operators
and Krein strings.  Some selected references are
\cite{albeverio.gesztesy.hoeghkrohn.holden:2005, arnold:1989,
flanders:1989, lifschitz:1989} for relevance in physics, and
\cite{atkinson:1964, kac:1999, kac.krein:1968, remling:2002}
for the relation to scalar second order differential or difference equations.

The theory of canonical systems was developed in works of Stieltjes, Weyl, Markov,
Krein, Kac and de~Branges. There is a vast literature, especially on spectral theory,
ranging from classical papers to very recent work.  As examples
we mention \cite{kac:1950, gohberg.krein:1967, debranges:1968,
winkler:1995, remling:2018, romanov:1408.6022v1, orcutt:1969,kac:1984, sakhnovich:1999,
hassi.snoo.winkler:2000, hassi.remling.snoo:2000,
snoo.winkler:2005, snoo.winkler:2005a,
kaltenbaeck.woracek:hskansys, arov.dym:2008, krein.langer:2014}.
Our standard reference is \cite{hassi.snoo.winkler:2000},
where the spectral theory of canonical systems is developed in a modern
operator-theoretic language.

With a Hamiltonian $H$ one can associate a Hilbert space $L^2(H)$ and a
(minimal) differential operator $S(H)$; see Section~\ref{sec:canonsys}.
The spectral theory of $S(H)$ changes drastically depending on the
left endpoint $a$ being in Weyl's limit circle case (LC) or Weyl's limit point case (LP),
i.e.\ whether for one (and hence for all) $x_0\in(a,b)$
\[
	\text{(LC)}:\quad \int_a^{x_0}\tr H(x)\,\rd x<\infty\qquad\text{ or }\qquad
	\text{(LP)}:\quad \int_a^{x_0}\tr H(x)\,\rd x=\infty.
\]
Note that because of the non-negativity of $H$ the Hamiltonian $H$ is in the
limit circle case if and only if all entries of $H$ are integrable at~$a$.

\bigskip
%
%
\noindent
\textbf{Limit circle case.}
\\[1ex]
Assume that $H$ is in the limit circle case at its left endpoint
(and, as always in this paper, in the limit point case at its right endpoint).
Then the operator $S(H)$ is symmetric with deficiency index $(1,1)$.
A complex-valued function $q_H$, the \emph{Weyl coefficient} of $H$, can be constructed
as follows.  Let $\uptheta(\cdot\,;z)$ and $\upvarphi(\cdot\,;z)$
be the solutions of \eqref{A30} that satisfy the initial conditions
$\uptheta(a;z)=\binom{1}{0}$ and $\upvarphi(a;z)=\binom{0}{1}$, respectively;
note that $H$ is integrable at $a$.
The Weyl coefficient $q_H$ is defined by
\begin{equation}\label{A150}
	q_H(z)\defeq\lim_{x\nearrow b}
	\frac{\uptheta_1(x;z)\tau+\uptheta_2(x;z)}{\upvarphi_1(x;z)\tau+\upvarphi_2(x;z)}\,,
	\qquad z\in\bb C\setminus\bb R,
\end{equation}
with $\tau\in\bb R\cup\{\infty\}$; the limit is independent of $\tau$
since $H$ is in the limit point case at $b$.
The function $q_H$ belongs to the \emph{Nevanlinna class} $\mc N_0$,
i.e.\ it is analytic in
$\bb C\setminus\bb R$, symmetric with respect to the real line in the sense that
$q_H(\qu z)=\qu{q_H(z)}$, $z\in\bb C\setminus\bb R$,
and maps the open upper half-plane $\bb C^+$ into $\bb C^+\cup\bb R$.

The Weyl coefficient $q_H$ can be used to construct a spectral measure and a Fourier transform.
Let $\mu_H$ be the measure in the Herglotz integral representation of $q_H$
(see \eqref{A129} below) appropriately including a possible point mass at $\infty$,
and define an integral transformation $\Theta_H$ by
\[
	(\Theta_H f)(t)\defeq \int_a^b\upvarphi(x;t)^TH(x)f(x)\,\rd x,\quad
	f\in L^2(H),\;\; \sup(\supp f)<b.
\]
Then a direct spectral theorem holds; more precisely, the following is true.
\begin{enumerate}[(1)]
\item The map $\Theta_H$ extends to an isometric isomorphism from $L^2(H)$ onto
	$L^2(\mu_H)$, where we tacitly understand that the space $L^2(\mu_H)$
	appropriately includes a possible point mass at $\infty$.
\item This extension of $\Theta_H$ establishes a unitary equivalence between the
	self-adjoint extension of $S(H)$ that is determined by the boundary condition $y_1(a)=0$
	and the operator $M_{\mu_H}$ of multiplication by the independent variable in
	the space $L^2(\mu_H)$.
\end{enumerate}
This direct theorem shows, in particular, that the mentioned self-adjoint extension
of $S(H)$ has simple spectrum.
\\[2ex]
An inverse spectral theorem was proved by L.~de Branges
in \cite{debranges:1960}--\cite{debranges:1962a}, in particular
\cite[Theorem~XII]{debranges:1961} and \cite[Theorem~VII]{debranges:1962a};
see also \cite{winkler:1995} for an explicit treatment.
These results include the following statements.
\begin{enumerate}[(1)]
\item Let a function $q$ in the Nevanlinna class $\mc N_0$ be given. Then there exists a
	Hamiltonian $H$ that is in the limit circle case at its left endpoint
	(and in the limit point case at its right endpoint) such that $q_H=q$.
\item Let a positive scalar measure $\mu$ with
	$\int_{\bb R}(1+t^2)^{-1}\rd\mu(t)<\infty$ be given (plus a possible point mass at
	$\infty$). Then there exists a Hamiltonian $H$ that is in the limit circle case at
	its left endpoint (and in the limit point case at its right endpoint) such that $\mu=\mu_H$
	(and possible point masses at $\infty$ coincide).
\item Let two Hamiltonians $H_1$ and $H_2$ be given, both being in the limit circle case at their left
	endpoints (and in the limit point case at their right endpoints).
	Then we have $q_{H_1}=q_{H_2}$ if and only if
	$H_1$ and $H_2$ are reparameterizations of each other; the latter means
	that $H_2(x)=H_1(\gamma(x))\gamma'(x)$ with some increasing bijection $\gamma$
	such that $\gamma$ and $\gamma^{-1}$ are absolutely continuous.
\item Let two Hamiltonians $H_1$ and $H_2$ be given, both being in the limit circle case at their left
	endpoints (and in the limit point case at their right endpoints).
	Then we have $\mu_{H_1}=\mu_{H_2}$ (and possible point masses at $\infty$ coincide)
	if and only if there exists a real constant $\alpha$ such that the Hamiltonians
	\[
		H_1,\qquad
		\smmatrix{1}{\alpha}{0}{1}H_2\smmatrix{1}{0}{\alpha}{1}
	\]
	are reparameterizations of each other.
\end{enumerate}

\bigskip
%
%
\noindent
\textbf{Limit point case.}
\\[1ex]
\noindent
If the limit point case prevails (also) at the left endpoint, much less can be said in general.
The operator $S(H)$ is self-adjoint, and its spectral multiplicity cannot exceed $2$.
A $2\!\times\!2$-matrix-valued Weyl coefficient can be defined. Via the
Titchmarsh--Kodaira formula, this leads to a Fourier transform onto an $L^2$-space with
respect to a $2\!\times\!2$-matrix-valued measure; see, e.g.\
\cite{hassi.snoo.winkler:2000}, and \cite{kodaira:1949} or \cite[\S2]{gilbert:1998}
for Schr\"odinger equations.

For Hamiltonians being in the limit point case non-simple spectrum can appear;
and this is not an exceptional case.
The class of all Hamiltonians that have simple spectrum --- despite being in the
limit point case at both endpoints --- can be characterized based on a theorem
of I.\,S.~Kac from the 1960s;
see \cite[Fundamental Theorem]{kac:1962a}\footnote{A proof is given in \cite{kac:1963a} (in Russian).}.
However, given the Hamiltonian $H$,
the condition given in Kac's Theorem is hardly accessible to computation.
To the best of our knowledge an
explicit characterization of simplicity of the spectrum is not known.
An easy-to-check sufficient condition for $S(H)$ having simple spectrum
follows from a result of L.~de~Branges; see \cite[Theorems~40 and 41]{debranges:1968}.

In the study of limit point Hamiltonians with simple spectrum there remain
some major drawbacks compared with the limit circle situation.
Even in the situation of de~Branges' Theorem there is neither a canonical way
to choose a scalar-valued spectral measure $\mu$ nor further information
on properties of $\mu$ can be obtained.
In view of this fact, naturally, there are no inverse statements
asserting existence or uniqueness of a Hamiltonian which would lead to a given measure.

\bigskip
%
%
\noindent
\textbf{The main results of the present paper.}
\\[1ex]
We specify a class $\bb H$ of Hamiltonians, which are
in the limit point case at both endpoints and for which a Weyl theory analogous to
the limit circle case can be developed.
This class $\bb H$ is a proper subclass of the one familiar from de~Branges' theorem
mentioned above, but it is still sufficiently large to cover many cases of interest.

For each Hamiltonian $H\in\bb H$ we prove the following direct spectral results.
\begin{enumerate}[(1)]
\item
	Every solution of equation \eqref{A30} attains regularized boundary values
	at $a$ in the sense that (at most) finitely many
	divergent terms are discarded in a well-defined way (Theorem~\ref{A32});
	the regularization depends on one free parameter $x_0\in(a,b)$.
	One can then define solutions $\uptheta$ and $\upvarphi$
	by prescribing the regularized boundary values.
	Hence an analogue of the Weyl coefficient, which we call
	\emph{singular Weyl coefficient}, can be defined with the help
	of $\uptheta$ and $\upvarphi$ as in \eqref{A150};
	this singular Weyl coefficient depends on the parameter $x_0$ but
	the dependence shows only in an additive real polynomial (Theorem~\ref{A94}).
\item
	A Fourier transform onto an $L^2$-space generated by a scalar measure exists.
	One measure with this property can be constructed in a canonical way
	via the singular Weyl coefficient, and this measure
	is independent of the parameter $x_0$ (Theorem~\ref{A107}).
	The corresponding Fourier transform and its inverse can be written
	as integral transforms (Theorem~\ref{A4}).
\end{enumerate}
Concerning the, now meaningfully posed, inverse spectral problem, we
\begin{enumerate}[(1)]
\setcounter{enumi}{2}
\item
	characterize the class of measures occurring via the
	mentioned construction (Theorem~\ref{A5});
\item
	establish global and local uniqueness results (Theorems~\ref{A6} and \ref{A31});
\item
	establish a one-to-one correspondence between the growth of the Hamiltonian $H$
	at $a$ and the growth of the spectral measure $\mu_H$ at infinity,
	measured by a positive integer $\Delta$ (Theorem~\ref{A107}).
\end{enumerate}

\bigskip
%
%
\noindent
\textbf{Sturm--Liouville equations.}
\\[1ex]
Recently, Sturm--Liouville equations, and in particular Schr\"odinger equations,
with two singular endpoints attracted a lot of attention;
for example, let us mention \cite{gesztesy.zinchenko:2006, fulton:2008,
fulton.langer:2010, fulton.langer.luger:2012,
kostenko.sakhnovich.teschl:2010, kostenko.teschl:2011,
kostenko.sakhnovich.teschl:2012, kostenko.sakhnovich.teschl:2012a,
eckhardt.teschl:2013,kostenko.teschl:2013,eckhardt:2014}.
Sturm--Liouville equations for which the corresponding operator
is bounded from below can be transformed into canonical systems
of the form \eqref{A30}; see \thref{A306}.
We consider two classes of Sturm--Liouville equations in detail:
first, equations without potential, i.e.\
\begin{equation}\label{A33}
	-\bigl(py'\bigr)' = \lambda wy
\end{equation}
with $p(x),w(x)>0$ a.e., $1/p,w$ locally integrable and either $1/p$ or $w$
integrable at $a$.  Such equations,
which are treated in Sections~\ref{sec-SL} and \ref{sec-SLw},
have many applications (see, e.g.\ \cite{bube.burridge:1983}
and \cite{mclaughlin:1986})
and include equations in impedance form, i.e.\ where $p=w$;
see, e.g.\ \cite{albeverio.hryniv.mykytyuk:2005}.
Second, we consider one-dimensional Schr\"odinger equations, i.e.\
\begin{equation}\label{A149}
	-y''+qy = \lambda y
\end{equation}
with $q$ locally integrable.  The class of equations we can treat includes
radial equations for Schr\"odinger equations with spherically symmetric potentials;
the corresponding operators are also called perturbed Bessel operators.
We apply our results on canonical systems to the Sturm--Liouville equations
\eqref{A33} and \eqref{A149}; in particular, we construct
singular Titchmarsh--Weyl coefficients, spectral measures and Fourier transforms,
and we prove inverse spectral theorems.

\bigskip
%
%
\pagebreak[3]
\noindent
\textbf{Methods employed.}
\\[1ex]
In order to establish our present results, we utilize the theory of
indefinite inner product spaces.
Our approach proceeds via Pontryagin space theory,
i.e.\ the theory of indefinite inner product spaces with a finite-dimensional
negative part.
In some sense our approach reaches as far as Pontryagin space models possibly can.
One key idea is to extend the Hamiltonian $H$ to the left by a so-called
indivisible interval so that the original left endpoint $a$ becomes an
interior point where $H$ is singular.
We can then apply the theory of generalized Hamiltonians,
developed in \cite{kaltenbaeck.woracek:p4db}--\cite{kaltenbaeck.woracek:p6db}
and also \cite{langer.woracek:gpinf},
for which corresponding operator models act in Pontryagin spaces
(in general, a generalized Hamiltonian can have a finite number of
interior singularities).

We use operator-theoretic tools like the spectral theory of self-adjoint relations,
models for generalized Nevanlinna functions and for generalized Hamiltonians,
as well as the theory of de~Branges Pontryagin spaces of entire functions.
In particular, proofs rely heavily on the theory developed
in \cite{langer.woracek:gpinf} and \cite{langer.woracek:ninfrep} and in
\cite{kaltenbaeck.woracek:p4db}--\cite{kaltenbaeck.woracek:p6db}.
We would like to mention that the underlying relation in the Pontryagin space is of the most
intriguing (but also most difficult to handle) kind: it is a proper relation having infinity
as a singular critical point with a neutral algebraic eigenspace.

\bigskip
%
%
\noindent
\textbf{Organization of the manuscript.}
\\[1ex]
The paper is divided into sections according to the following table.
\vspace*{-0mm}
\begin{center}
	\rule{90mm}{0.5pt}\\[2mm]
	{\textbf{Table of contents}}
\end{center}
\vspace*{-1ex}
{\footnotesize
PART I: General Theory
\begin{enumerate}
\setcounter{enumi}{1}
\item The two basic classes \hfill p.~\pageref{sec-H+M}
\item Preliminaries from indefinite theory \hfill p.~\pageref{sec-preliminaries}
\item Construction of the spectral measure \hfill p.~\pageref{sec-measure}
\item The Fourier transform \hfill p.~\pageref{sec-fourier}
\item Inverse theorems \hfill p.~\pageref{sec-inverse}
\setcounter{mytoc}{\value{enumi}}
\end{enumerate}
\vspace*{1ex}
PART II: Applications to Sturm--Liouville Equations
\begin{enumerate}
\setcounter{enumi}{\value{mytoc}}
\item Sturm--Liouville equations without potential: singular $1/p$ \hfill p.~\pageref{sec-SL}
\item Sturm--Liouville equations without potential: singular $w$ \hfill p.~\pageref{sec-SLw}
\item Schr\"odinger equations \hfill p.~\pageref{sec-Schroedinger}
\end{enumerate}
}
\vspace*{-2ex}
\begin{center}
	\rule{90mm}{0.5pt}
\end{center}

\noindent
In Section~\ref{sec-H+M} we introduce the class $\bb H$ of Hamiltonians
that is treated in our paper.  The definition involves a certain growth condition of
the Hamiltonian at the left endpoint $a$.
We associate a positive integer, $\Delta(H)$, with each $H\in\bb H$, which measures
the growth of $H$ at $a$.
Further, we define a class $\bb M$ of Borel measures on $\bb R$ that satisfy
a certain growth condition at infinity; this class will turn out
to be the set of spectral measures of Hamiltonians from $\bb H$.
In Section~\ref{sec-preliminaries} we recall the definition and certain properties
of generalized Nevanlinna functions and the operator that is connected
with equation \eqref{A30}.  Moreover, we recall the notion of
generalized Hamiltonians, a certain subclass of generalized Hamiltonians
that have only one interior singularity, and corresponding operator models.
In Section~\ref{sec-measure} we show that solutions of \eqref{A30} attain
regularized boundary values at $a$ (\thref{A32}), we construct singular
Weyl coefficients (\thref{A94}) and construct a spectral measure
via a Stieltjes-type inversion formula (\thref{A107}).
The Fourier transform is constructed in Section~\ref{sec-fourier} (\thref{A4});
this shows, in particular, that the spectrum is simple.
Inverse spectral theorems (existence and global and local uniqueness theorem)
are proved in Section~\ref{sec-inverse} (Theorems~\ref{A5}, \ref{A6} and \ref{A31}).

In the second part of the paper we consider Sturm--Liouville equations.
First we consider equations of the form \eqref{A33}.
The case when $1/p$ is not integrable at $a$ is considered in Section~\ref{sec-SL};
the case when $w$ is not integrable at $a$ is studied in Section~\ref{sec-SLw}.
Finally, Schr\"odinger equations of the form \eqref{A149} are investigated
in Section~\ref{sec-Schroedinger}.

\bigskip
%
%
\noindent
\textbf{Acknowledgements.}
\\
The first author gratefully acknowledges the support of the Nuffield Foundation,
grant no.\ NAL/01159/G,
and the Engineering and Physical Sciences Research Council (EPSRC),
grant no.\ EP/E037844/1.
The second author was supported by a joint project of the Austrian Science Fund (FWF, I\,1536--N25)
and the Russian Foundation for Basic Research (RFBR, 13-01-91002-ANF).

\bigskip

%
%
\begin{center}
{\Large\textbf{PART I: \\[1ex] General Theory}}
\end{center}
%
%

\noindent
In the first part, which comprises Sections~\ref{sec-H+M}--\ref{sec-inverse},
the direct and inverse spectral theory of canonical systems with two singular
endpoints is developed.

%
%
\section{The two basic classes}
\label{sec-H+M}
%
%

We start with the definition and a brief discussion of the two major objects of our investigation.
These are a class $\bb H$ of Hamiltonians and a class $\bb M$ of measures, which will
turn out to correspond to each other.

\subsection{The class {\boldmath$\bb H$} of Hamiltonians}

Let us state the definition of Hamiltonians again explicitly:
by a \emph{Hamiltonian} $H=(h_{ij})_{i,j=1}^2$ we understand a function defined on some
(non-empty and possibly unbounded)
interval $(a,b)$ whose values are real, non-negative $2\!\times\!2$-matrices, which
is locally integrable and which does not vanish on any set of positive measure.
In the rest of the paper we shall also write $\dom(H)\defeq(a,b)$
if $H$ is defined on $(a,b)$.

We say that two Hamiltonians $H_1$ and $H_2$ defined on intervals $(a_1,b_1)$
and $(a_2,b_2)$, respectively, are \emph{reparameterizations} of each other
if there exists an increasing bijection $\gamma:(a_2,b_2)\to(a_1,b_1)$
such that $\gamma$ and $\gamma^{-1}$ are both absolutely continuous and
\begin{equation}\label{A72}
	H_2(x)=H_1\bigl(\gamma(x)\bigr)\cdot\gamma'(x), \qquad x\in(a_2,b_2) \;\;\text{a.e.}
\end{equation}
Note that in this situation $y$ is a solution of \eqref{A30} with $H=H_1$
if and only if $\tilde y$, where $\tilde y(x)=y(\gamma(x))$,
is a solution of \eqref{A30} with $H=H_2$.

\begin{remark}\thlab{A101}
	As a rule of thumb, Hamiltonians which are reparameterizations of each other share all
	their essential properties. For a detailed and explicit exposition of reparameterizations
	in an up-to-date language, see \cite{winkler.woracek:nnham} (in particular, Theorem~3.8
	therein).
\end{remark}

We also recall the notion of indivisible intervals.
An interval $(\alpha,\beta)\subseteq(a,b)$ is called
\emph{$H$-indivisible} (or just \emph{indivisible}) \emph{of type $\phi$} if
\begin{equation}\label{A136}
	H(x) = h(x)\xi_\phi\xi_\phi^T, \qquad x\in(\alpha,\beta),
\end{equation}
where $\xi_\phi=(\cos\phi,\sin\phi)^T$ and $h$ is a locally integrable
function that is positive almost everywhere;
see, e.g.\ \cite{kac:1984}.
An indivisible interval $(\alpha,\beta)$ is called \emph{maximal} if it is not
contained in any larger indivisible interval.

\begin{definition}\thlab{A1}
	Let $H=(h_{ij})_{i,j=1}^2$ be a Hamiltonian defined on $(a,b)$.
	We say that $H$ belongs to the class $\bb H$ if $H$ is in the limit point case
	at both endpoints,
	the interval $(a,b)$ is neither one indivisible interval
	nor the union of two indivisible intervals,
	and $H$ satisfies the following conditions {\rm(I)}, {\rm(HS)} and {\rm($\Delta$)}.
	\begin{itemize}
	\item[{\rm(I)}] For one (and hence for all) $x_0\in(a,b)$,
		\[
			\intop_a^{x_0} h_{22}(x)\,\rd x<\infty.
		\]
	\item[{\rm(HS)}] For one (and hence for all) $x_0\in(a,b)$,
		\[
			\intop_a^{x_0} \intop_a^x h_{22}(t)\,\rd t\, h_{11}(x)\,\rd x<\infty.
		\]
	\item[{\rm($\Delta$)}] Let $x_0\in(a,b)$ and
		define functions $X_k:(a,x_0]\to\bb C^2$
		recursively by
		\begin{alignat*}{2}
			X_0(x)&\defeq\binom{1}{0}, && x\in(a,x_0], \\[1ex]
			X_k(x)&\defeq\intop_{x_0}^x JH(t)X_{k-1}(t)\,\rd t,\qquad && x\in(a,x_0],\;k\in\bb N.
		\end{alignat*}
		There exists a number $N\in\bb N_0$ such that
		\begin{equation}\label{A2}
			L^2\big(H|_{(a,x_0)}\big)\cap\spn\big\{X_k:\,k\leq N\big\}\neq\{0\}.
		\end{equation}
	\end{itemize}
	If $H\in\bb H$, we denote by $\Delta(H)$ the smallest non-negative integer $N$ such that
	\eqref{A2} holds.
\end{definition}

\noindent
It is proved in \cite[Lemma~3.12]{kaltenbaeck.woracek:p4db} that this
definition is justified, namely that the validity of {\rm($\Delta$)} and
the number $\Delta(H)$ do not depend on the choice of~$x_0$
(for {\rm(I)} and {\rm(HS)} this is trivial to check).

Notice that, for $H\in\bb H$, we always have $\Delta(H)>0$. This follows since
we assume limit point case at $a$. Namely, for each $x_0\in(a,b)$,
the constant function $(0,\,1)^T$ belongs to $L^2(H|_{(a,x_0)})$ by {\rm(I)}, and hence
the constant $(1,\,0)^T$ cannot be in this space.

\begin{remark}\thlab{A133}
	We assume that $(a,b)$ is neither one indivisible interval nor the union of
	two indivisible intervals since, otherwise, the corresponding space $L^2(H)$
	(defined in Section~\ref{sec:canonsys})
	and hence also the Fourier transform would be trivial.
\end{remark}

\begin{remark}\thlab{A91}
	The conditions {\rm(I)} and {\rm(HS)} are, up to a normalization and exchanging
	upper and lower rows,
	precisely the conditions of de~Branges' Theorem \cite[Theorem~41]{debranges:1968}.
	Note that under the conditions {\rm(I)} and {\rm(HS)} any self-adjoint realization
	corresponding to $H|_{(a,x_0)}$ has a Hilbert--Schmidt resolvent.
	The additional condition {\rm($\Delta$)} arose only recently in the context of
	indefinite canonical systems; we recall more details in \S3.2 below.
\end{remark}

\noindent
In general it is difficult to decide whether a given Hamiltonian satisfies {\rm($\Delta$)}.
Contrasting {\rm(I)} and {\rm(HS)} the condition {\rm($\Delta$)} is of recursive nature and not accessible
by simple computation. An easier-to-handle (though still recursive) criterion for the validity
of {\rm($\Delta$)} is available for Hamiltonians of diagonal form,
cf.\ \cite[Theorem~3.7]{winkler.woracek:del} and Section~\ref{sec-SL}.
Using this criterion, various examples can be constructed.
The following two examples are taken from \cite[Corollary~3.14 and Example~3.15]{winkler.woracek:del}.

\begin{example}\thlab{A92}
	Let $\alpha\in\bb R$ and set
	\[
		H_\alpha(x) \defeq \begin{pmatrix} x^{-\alpha} & 0 \\[1ex] 0 & 1 \end{pmatrix},
		\qquad x\in(0,\infty).
	\]
	Then $H_\alpha$ is in the limit point case at $\infty$ and satisfies {\rm(I)} at $0$.
	Depending on the value of $\alpha$, the following conditions hold:
	\[
		\begin{array}{l||c|l}
			\text{\rm value of $\alpha$} & \text{\rm (LP)/(LC) at $0$} &
			\text{\rm(HS) and ($\Delta$)}
			\\[2pt] \hline
			\raisebox{-1mm}{\rule{0pt}{5mm}}
			\alpha<1 & \text{\rm (LC)} & \text{\rm both hold (trivially)}
			\\[2pt]
			1\leq\alpha<2 & \text{\rm (LP)} & \text{\rm both hold}
			\\[2pt]
			\alpha\geq 2 & \text{\rm (LP)} & \text{\rm none holds}
		\end{array}
	\]
	Hence, we have $H_\alpha\in\bb H$ for each $\alpha\in[1,2)$ but not for other values of
	$\alpha$.

	The number $\Delta(H_\alpha)$ can be computed, namely,
	\[
		\Delta(H_\alpha)=n,\qquad \text{when }\alpha\in\Big(2-\frac 1n,2-\frac 1{n+1}\Big)
		\;\text{ with }\; n\in\bb N.
	\]
	This shows that, for a Hamiltonian of class $\bb H$, there are no a priori
	restrictions on the value of the number $\Delta(H)$.
	Computing $\Delta(H_\alpha)$ for $\alpha=2-\frac 1n$ with $n\in\bb N$ is equally well possible,
	but requires more elaborate computations.  These have not been carried out in
	\cite{winkler.woracek:del} but will be made available elsewhere.
\end{example}

\begin{example}\thlab{A93}
	Consider the Hamiltonian
	\[
		H(x) \defeq \begin{pmatrix} (x\ln x)^{-2} & 0 \\[1ex] 0 & 1 \end{pmatrix},
		\qquad x\in(0,1).
	\]
	This Hamiltonian is in the limit point case at $0$ and at $1$, satisfies {\rm(I)} and
	{\rm(HS)} at $0$, but does not satisfy {\rm($\Delta$)}.

	This example shows that the presently considered class $\bb H$ is a
	proper subclass of the one treated in \cite[Theorem~41]{debranges:1968}.
\end{example}

\subsection{The class {\boldmath$\bb M$} of measures}

By a \emph{positive Borel measure on $\bb R$} we understand
a (not necessarily finite) positive measure defined on the $\sigma$-algebra of
all Borel subsets of $\bb R$ which takes finite values on compact sets.

\begin{definition}\thlab{A3}
	Let $\mu$ be a positive Borel measure on $\bb R$.
	We say that $\mu$ belongs to the class $\bb M$ if there exists a number
	$N\in\bb N_0$ such that
	\begin{equation}\label{A22}
		\intop_{\bb R}\frac{\rd\mu(t)}{(1+t^2)^{N+1}}<\infty.
	\end{equation}
	If $\mu\in\bb M$, we denote by $\Delta(\mu)$ the smallest non-negative integer $N$
	such that \eqref{A22} holds.
\end{definition}

\noindent
This class of measures is known from Pontryagin space theory.
A measure $\mu$ belongs to $\bb M$ if and only if it is the measure in the
distributional representation of some generalized Nevanlinna function
with a certain spectral behaviour, cf.\ \cite{krein.langer:1977} and
\cite[Theorems~2.8 and 3.9]{langer.woracek:ninfrep}.
In the classical (positive definite) setting, this corresponds to the fact
that a positive Borel measure $\mu$ satisfies $\int_{\bb R}(1+t^2)^{-1}\,\rd\mu(t)$
if and only if it is the measure in the Herglotz integral
representation of some Nevanlinna function.  We recall details in \S3.1 below.

%
%
\section{Preliminaries from indefinite theory}
\label{sec-preliminaries}
%
%

Our approach to direct and inverse spectral theory for Hamiltonians of class $\bb H$
is based on the theory of indefinite canonical systems
and their Pontryagin space operator models as developed in
\cite{kaltenbaeck.woracek:p4db}--\cite{kaltenbaeck.woracek:p6db} and further in
\cite{langer.woracek:esmod, langer.woracek:gpinf}.
In this preliminary section we recall the relevant notions and theorems.
For the theory of Pontryagin spaces we refer the reader, e.g.\ to \cite{bognar:1974}.

\subsection{Generalized Nevanlinna functions and the class $\mc N_{<\infty}^{(\infty)}$}

As we already mentioned in the introduction, a function $q$ is said to be a
\emph{Nevanlinna function} if it is analytic in $\bb C\setminus\bb R$, satisfies $q(\qu z)=\qu{q(z)}$
for $z\in\bb C\setminus\bb R$ and $\Im q(z)\geq 0$ for $z\in\bb C^+$. We denote the set of all Nevanlinna
functions by $\mc N_0$.

In Pontryagin space theory an indefinite analogue of this notion appears
and plays a significant role; see, e.g.\ \cite{krein.langer:1973, krein.langer:1977}.

\begin{definition}\thlab{A115}
	A function $q$ is called a \emph{generalized Nevanlinna function} if it is meromorphic
	in $\bb C\setminus\bb R$ and has the following properties (i) and (ii):
	\begin{enumerate}[{\rm(i)}]
		\item $q(\qu z)=\qu{q(z)}$ for $z\in\rho(q)$, where $\rho(q)$ denotes the
			domain of analyticity of $q$ in $\bb C\setminus\bb R$;
		\item the reproducing kernel
			\[
				K_q(w,z) \defeq \frac{q(z)-\qu{q(w)}}{z-\qu w}\,,
				\qquad z,w\in\rho(q),
			\]
			has a finite number of negative squares;
			the latter means that there exists a $\kappa\in\bb N_0$ so that
			for every choice of $n\in\bb N$ and $z_1,\ldots,z_n\in\rho(q)$
			the matrices $(K_q(z_i,z_j))_{i,j=1}^n$ have at most $\kappa$ negative eigenvalues.
	\end{enumerate}
	We denote the set of all generalized Nevanlinna functions by $\mc N_{<\infty}$.
	Moreover, if $q\in N_{<\infty}$, we denote the actual number of negative squares
	of the kernel $K_q$ (i.e.\ the minimal $\kappa$ in the above definition)
	by $\ind_-q$.  Further, we set $\mc N_\kappa\defeq\{q\in\mc N_{<\infty}:\ind_-q=\kappa\}$
	for $\kappa\in\bb N_0$.
\end{definition}

\noindent
That this definition is indeed an extension of the definition of $\mc N_0$, i.e.\
that the class $\mc N_0$ in Definition~\ref{A115} coincides with the class defined
before Definition~\ref{A115} is a classical result,
which can be traced back to as far as \cite{herglotz:1911} or \cite{pick:1915}.

Let $q\in\mc N_0$. Using the representation of the positive harmonic function
$\Im q$ as a Poisson integral, one easily
obtains a representation of $q$ with a Cauchy-type integral.

\begin{nremark}{Herglotz integral representation of $\mc N_0$-functions}\thlab{A117}
	A function $q$ belongs to the class $\mc N_0$ if and only if it can be represented in the form
	\begin{equation}\label{A129}
		q(z)=a+bz+\int_{\bb R}\Big(\frac 1{t-z}-\frac t{1+t^2}\Big)\,\rd\mu(t),
		\qquad z\in\bb C\setminus\bb R,
	\end{equation}
	with $a\in\bb R$, $b\geq 0$ and a positive Borel measure $\mu$ satisfying
	$\int_{\bb R}(1+t^2)^{-1}\,\rd\mu(t)<\infty$.
\end{nremark}

\noindent
The analogue of this integral representation in the indefinite setting is a
distributional representation of a generalized Nevanlinna function. In essence this is shown in
\cite{krein.langer:1977}, where an integral representation of $q\in\mc N_{<\infty}$ was given
without using the language of distributions. The distributional viewpoint was first mentioned
in \cite[Introduction, p.~253]{jonas.langer.textorius:1992}, established more thoroughly in
\cite[Corollary~3.5]{kaltenbaeck.woracek:p2db} and refined in
\cite[Proposition~5.4]{kaltenbaeck.winkler.woracek:nksym}.
The formulation given below is taken from our paper \cite{langer.woracek:ninfrep}.
This paper contains several results which are crucial for the present discussion
and is our standard reference in the context of distributional representations.

Before we can provide the actual statement, we need to introduce some notation.
First, we denote by $\bb R(z)$ the set of all rational functions with
real coefficients.
Second, we denote by $\qu{\bb R}$ the one-point compactification of the real line
considered as a $C^\infty$-manifold in the usual way.
Moreover, for each $z\in\bb C\setminus\bb R$, let $\beta_z:\qu{\bb R}\to\bb C$ be defined by
\[
	\beta_z(t) \defeq
	\begin{cases}
		\dfrac{1+xz}{x-z}\,, &\quad x\in\bb R, \\[2ex]
		z, &\quad x=\infty.
	\end{cases}
\]
Third, for a function $f$, set $f^\#(z)\defeq\qu{f(\qu{z})}$ whenever $\qu z\in\bb C$ is
in the domain of $f$.  Further, we denote by $\widetilde{\mc D}'(\qu{\bb R})$ the set of
all distributional densities on $\qu{\bb R}$; see, e.g.\ \cite{hoermander:1990}
or \cite{langer.woracek:ninfrep}.
With each $\upphi\in\widetilde{\mc D}'(\qu{\bb R})$
one can associate a linear functional on $C^\infty(\qu{\bb R})$,
which is again denoted by $\upphi$; see \cite[(2.2)]{langer.woracek:ninfrep}.
Next, we denote by $\mc F(\qu{\bb R})$ the set of all $\upphi\in\widetilde{\mc D}'(\qu{\bb R})$
for which there exists a finite subset $F$ of $\qu{\bb R}$ such that $\upphi$ acts as a
positive measure on $\qu{\bb R}\setminus F$;
for details see \cite[Definitions~2.1 and 2.3]{langer.woracek:ninfrep}.
Finally, $\mc F_{\{\infty\}}$ is the set of $\upphi\in\mc F(\qu{\bb R})$
that act as a positive measure on $\bb R$.  For $\upphi\in\mc F_{\{\infty\}}$
we denote by $\mu_\upphi$ the unique positive Borel measure on $\bb R$ such that
\begin{equation}\label{A272}
	\upphi(f) = \intop_{\bb R} f(x)\frac{\rd\mu_\upphi(x)}{1+x^2}\,,
	\qquad f\in C^\infty(\qu{\bb R}),\,\supp f\subseteq\bb R;
\end{equation}
see \cite[Definition~2.4]{langer.woracek:ninfrep}.

\begin{nremark}{Distributional representation of $\mc N_{<\infty}$-functions
	{\rm\cite[Proposition~5.4]{kaltenbaeck.winkler.woracek:nksym}}}\thlab{A118}
	Let $\upphi\in\mc F(\qu{\bb R})$ and $r\in\bb R(z)$. Then the function
	\begin{equation}\label{A34}
		q(z) \defeq r(z)+\upphi(\beta_z)
	\end{equation}
	belongs to $\mc N_{<\infty}$.

	Conversely, let $q\in\mc N_{<\infty}$ be given. Then there exist unique
	$\upphi\in\mc F(\qu{\bb R})$ and $r\in\bb R(z)$ such that
	\begin{enumerate}[{\rm(i)}]
		\item the representation \eqref{A34} holds;
		\item $r$ is analytic on $\bb R$ and remains bounded for $|z|\to\infty$.
	\end{enumerate}
\end{nremark}

\noindent
In the present paper the following subclass of $\mc N_{<\infty}$ plays a central role.

\begin{definition}\thlab{A116}
	We denote by $\mc N_\kappa^{(\infty)}$, $\kappa\in\bb N_0$,
	the set of all functions $q\in\mc N_\kappa$ such that
	\begin{equation}\label{A7}
		\lim_{z\nontang i\infty}\frac{q(z)}{z^{2\kappa-1}}\in(-\infty,0)
		\qquad \text{or}\qquad
		\lim_{z\nontang i\infty}\Bigl\lvert\frac{q(z)}{z^{2\kappa-1}}\Bigr\rvert=\infty,
	\end{equation}
	where $\nontang$ denotes the non-tangential limit, i.e.\ $z\to\infty$ inside some
	Stolz angle $\{z\in\bb C:\varepsilon\leq\arg z\leq\pi-\varepsilon\}$ for
	one (and hence for all) $\varepsilon\in(0,\frac\pi2)$.
	Moreover, set
	\[
		\mc N_{<\infty}^{(\infty)} \defeq \bigcup_{\kappa\in\bb N_0} \mc N_\kappa^{(\infty)}.
	\]
\end{definition}

\noindent
Let us stress that the significance of the condition in this definition is not that
\eqref{A7} holds \emph{for some $\kappa$}, but that it holds \emph{exactly for $\kappa=\ind_-q$}.
Note that $\mc N_0^{(\infty)}=\mc N_0$.

The classes $\mc N_\kappa^{(\infty)}$ appeared often in the recent literature,
where they are also denoted by $\mc N_\kappa^\infty$.
Let us mention, for instance, \cite{dijksma.shondin:2000,
dijksma.langer.shondin.zeinstra:2000, fulton.langer:2010, kurasov.luger:2011},
where Sturm--Liouville equations with singular endpoints or singular perturbations
of self-adjoint operators were studied,
and \cite{dijksma.langer.shondin:2004} in connection with rank one perturbations
at infinite coupling, and \cite{dijksma.luger.shondin:2006, dijksma.luger.shondin:2009,
dijksma.luger.shondin:2010}
where operator models of such functions were investigated.
The class $\mc N_{<\infty}^{(\infty)}$ has an operator-theoretic interpretation,
namely, the self-adjoint relation in the operator/relation representation\footnote{A detailed account
on the operator representation of scalar-valued generalized Nevanlinna functions can be found in
\cite[\S1]{krein.langer:1977}. For a slightly different viewpoint and results for
operator- (or matrix-) valued functions, see \cite[\S3]{krein.langer:1973}, \cite{daho.langer:1985}, or some
of the vast more recent literature on operator models.} of $q$
has $\infty$ has its only spectral point of non-positive type;
for details see also \cite[\S5]{langer.woracek:ninfrep}.

For our present considerations it is essential that the distributional representation
of a generalized Nevanlinna function $q$ takes a simple form if $q\in\mc N_{<\infty}^{(\infty)}$.
The following result is contained in \cite[Theorem~3.9\,(i), (ii)]{langer.woracek:ninfrep}.

\begin{nremark}{Distributional representation of $\mc N_{<\infty}^{(\infty)}$-functions}\thlab{A50}
	A function $q$ belongs to the class $\mc N_{<\infty}^{(\infty)}$ if and only if it
	can be represented as
	\begin{equation}\label{A273}
		q(z) \defeq r+\upphi(\beta_z)
	\end{equation}
	with a real constant $r$ and a distributional density $\upphi\in\mc F_{\{\infty\}}$.
	Denote by $\mu_q$ the measure in \eqref{A272} that is connected with the
	distributional density $\upphi$ in \eqref{A273}, i.e.\ $\mu_q\defeq\mu_{\upphi_q}$
	with notation from \eqref{A272}.
	Then a Stieltjes inversion formula is valid:
	\begin{equation}\label{A142}
		\mu_q\bigl([a,b]\bigr) = \frac{1}{\pi}\lim_{\eps\searrow0}\lim_{\delta\searrow0}
		\int_{a-\eps}^{b+\eps} \Im q(t+i\delta)\rd t, \qquad [a,b]\subseteq\bb R;
	\end{equation}
	see \cite[Theorem~3.9\,(ii)]{langer.woracek:ninfrep}.
\end{nremark}

\begin{nremark}{Operator model}\thlab{A139}
	With a distributional density $\upphi\in\mc F_{\{\infty\}}$ one can associate
	a Pontryagin space $\Pi(\upphi)$, which is the completion of $C^\infty(\qu{\bb R})$
	with respect to the inner product
	\[
		[f,g]_\upphi \defeq \upphi(f\qu{g}), \qquad f,g\in C^\infty(\qu{\bb R}),
	\]
	and a self-adjoint relation $A_\upphi$ in $\Pi(\upphi)$;
	see \cite{jonas.langer.textorius:1992}, \cite{kaltenbaeck.woracek:p2db}
	or \cite[\S5]{langer.woracek:ninfrep}.
	The space $\Pi(\upphi)$ contains the following set
	\begin{equation}\label{A141}
		\Bigl\{f\in L^2\Bigl(\frac{\mu_\upphi(x)}{1+x^2}\Bigr):
		\supp f\text{ is compact}\Bigr\}.
	\end{equation}
	Let $\mc E_{A_\upphi}(\infty)$ be the algebraic eigenspace at infinity
	of $A_\upphi$.  By \cite[Theorem~5.3]{langer.woracek:ninfrep} there exists
	an isometric, continuous, surjective map
	\[
		\psi(\upphi): \mc E_{A_\upphi}(\infty)^{[\perp]}
		\to L^2\Bigl(\frac{\mu_\upphi(x)}{1+x^2}\Bigr),
	\]
	which acts as the identity on functions from the set in \eqref{A141}.
\end{nremark}

\noindent
Recall from \cite{dijksma.langer.luger.shondin:2000}
that every function $q\in\mc N_{\kappa}^{(\infty)}$ can be written as
$q(z) = p(z)q_0(z)$ where $q_0\in\mc N_0$ and $p$ is a monic real polynomial of
degree $2\kappa$.  Let $\mu_q$ be the measure associated with $q$ as in \ref{A50}
and let $\mu_0$ be the measure in the integral representation \eqref{A129}
of $q_0$.  Then the Stieltjes inversion formula \eqref{A142} implies that
\[
	\mu_q\bigl([a,b]\bigr) = \int_{[a,b]} p(t)\rd\mu_0(t)
\]
for every finite interval $[a,b]$ with $\mu_q(\{a\})=\mu_q(\{b\})=0$.
This, together with \cite[Corollary~3.1]{langer.luger.matsaev:2011}
immediately yields the following lemma.

\begin{lemma}\thlab{A87}
	Let $\kappa\in\bb N$, let $q_n\in\mc N_\kappa^{(\infty)}$ for $n\in\bb N$
	and let $q\in\mc N_\kappa^{(\infty)}$ such that $q_n(z) \to q(z)$ locally
	uniformly on $\bb C\setminus\bb R$.  Moreover, let $\mu_{q_n}$ and $\mu_q$
	be the measures that are associated with $q_n$ and $q$, respectively,
	as in \ref{A50}.  Then, for every interval $[a,b]$
	with $\mu_q(\{a\})=\mu_q(\{b\})=0$, we have
	\[
		\lim_{n\to\infty}\mu_{q_n}\bigl([a,b]\bigr) = \mu_q\bigl([a,b]\bigr).
	\]
	\popQED\qed
\end{lemma}

\noindent
In this lemma the assumption that $\ind_-q=\ind_-q_n$ is crucial.

\subsection{The operator associated with a canonical system}
\label{sec:canonsys}

We recall the definition of the space $L^2(H)$ and the corresponding operator.
Note that the notion of indivisible intervals and the vector $\xi_\phi$ were
defined in Section~\ref{sec-H+M}.
The space $L^2(H)$ is the space
of measurable functions $f$ defined on $(a,b)$ with values in $\bb C^2$ which satisfy
$\int_a^b f^*Hf<\infty$ and have the property that $\xi_\phi^T f$ is constant
on every indivisible interval of type $\phi$, factorized with respect
to the equivalence relation $=_H$ where
\[
	f =_H g \qquad\Longleftrightarrow\qquad H(f-g) = 0 \quad \text{a.e.}
\]
In the space $L^2(H)$ the \emph{operator} $T(H)$ is defined via its graph as
\begin{equation}\label{A143}
\begin{aligned}
	\hspace*{-1.5ex} T(H) \defeq &\, \Bigl\{(f;g)\in\bigl(L^2(H)\bigr)^2\!:
	\exists\, \text{ representatives } \hat f,\,\hat g \text{ of } f,\,g \text{ such that } \\
	&\,\hat f \text{ is locally absolutely continuous and }
	\hat f' = JH\hat g \;\;\text{a.e.\ on } (a,b)\Bigr\}. \!\!
\end{aligned}
\end{equation}
Since $H$ is in the limit point case at both endpoints, the operator $T(H)$
is self-adjoint; see, e.g.\ \cite[\S 6]{hassi.snoo.winkler:2000}.
If $H$ is in the limit circle case at the left endpoint, then $T(H)$
is the maximal operator (or relation);
its adjoint, the minimal operator $S(H)$, is a symmetric operator.

\subsection{General Hamiltonians}

The definition of the Pontryagin space analogue of a Hamiltonian and
the description of its associated canonical system and the operator model
is quite long and involved.  Here we give only an intuitive picture; for a
complete and logically sound formulation we refer to
\cite[\S8]{kaltenbaeck.woracek:p4db} or \cite[Definitions~2.16--2.18]{langer.woracek:gpinf}.
The latter paper is our standard reference in the context of general Hamiltonians.

A \emph{general Hamiltonian $\mf h$} is a collection of data
\begin{align*}
	\mf h:\qquad & n\in\bb N_0,\;\;
		-\infty\leq\sigma_0<\sigma_1<\ldots<\sigma_{n+1}\leq\infty;
	\\[1ex]
	& \text{Hamiltonians }
	H_i:(\sigma_i,\sigma_{i+1})\to\bb R^{2\times 2},\;\; i=1,\ldots,n;
	\\[1ex]
	& \oe_i\in\bb N_0,\;\; b_{i,1},\ldots,b_{i,\oe+1}\in\bb R,\;\;
	d_{i,0},\ldots,d_{i,2\Delta_i-1}\in\bb R,\quad i=1,\ldots,n
	;
	\\[1ex]
	& E\subseteq\{\sigma_0,\sigma_{n+1}\}\cup\bigcup_{i=0}^n(\sigma_i,\sigma_{i+1});
\end{align*}
where, among others, the following conditions have to be satisfied:
\begin{itemize}
\item[--] $H_0$ is in the limit circle case at $\sigma_0$; if $n\ge1$, then
	$H_{i-1}$ is in the limit point case at $\sigma_i$ and $H_i$ is in the limit point case
	at $\sigma_i$ for $i=1,\ldots,n$;
\item[--] the growth of the Hamiltonians $H_i$ towards the points
	$\sigma_1,\ldots,\sigma_n$ is restricted; the number $\Delta_i\in\bb N$ is a certain
	measure for this growth;
\item[--] two adjacent Hamiltonians satisfy an interface condition at their common endpoint.
\end{itemize}
The general Hamiltonian is called \emph{singular} if $H_n$ is in the limit point
case at $\sigma_{n+1}$; it is called \emph{regular} if $H_n$ is in the limit circle
case at $\sigma_{n+1}$.
Moreover, let $H$ be the function defined on $\bigcup_{i=0}^n (\sigma_i,\sigma_{i+1})$
such that $H|_{(\sigma_i,\sigma_{i+1})}=H_i$ for $i=0,\dots,n$.

The following visualization of the canonical system associated with a general
Hamiltonian may be helpful:
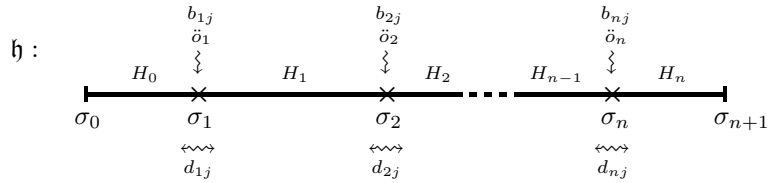
\begin{center}
\setlength\unitlength{1cm}
\begin{picture}(12,3.6)(-6,-1.8)
\put(-5.5,0.5){$\mf h:$}
\put(-4.65,-0.4){$\sigma_0$}
\put(-3.9,0.2){$\scriptstyle H_0$}
\put(-3.15,-0.088){\boldmath$\times$}
\put(-3.1,0.6){\begin{rotate}{270}{$\rightsquigarrow$}\end{rotate}}
\put(-3.15,1.0){$\scriptstyle b_{1j}$}
\put(-3.1,0.7){$\scriptstyle\oe_1$}
\put(-3.15,-0.4){$\sigma_1$}
\put(-3.25,-0.8){$\leftrightsquigarrow$}
\put(-3.2,-1.05){$\scriptstyle d_{1j}$}
\put(-1.9,0.2){$\scriptstyle H_1$}
\put(-0.65,-0.088){\boldmath$\times$}
\put(-0.6,0.6){\begin{rotate}{270}{$\rightsquigarrow$}\end{rotate}}
\put(-0.65,1.0){$\scriptstyle b_{2j}$}
\put(-0.6,0.7){$\scriptstyle\oe_2$}
\put(-0.65,-0.4){$\sigma_2$}
\put(-0.75,-0.8){$\leftrightsquigarrow$}
\put(-0.7,-1.05){$\scriptstyle d_{2j}$}
\put(0,0.2){$\scriptstyle H_2$}
\put(1.4,0.2){$\scriptstyle H_{n-1}$}
\put(2.35,-0.088){\boldmath$\times$}
\put(2.4,0.6){\begin{rotate}{270}{$\rightsquigarrow$}\end{rotate}}
\put(2.35,1.0){$\scriptstyle b_{nj}$}
\put(2.4,0.7){$\scriptstyle\oe_n$}
\put(2.35,-0.4){$\sigma_n$}
\put(2.25,-0.8){$\leftrightsquigarrow$}
\put(2.3,-1.05){$\scriptstyle d_{nj}$}
\put(3.1,0.2){$\scriptstyle H_n$}
\put(3.85,-0.4){$\sigma_{n+1}$}
\Thicklines
\put(-4.5,0){\line(1,0){5}}
\put(-4.5,-0.1){\line(0,1){0.2}}
\put(0.6,0){\line(1,0){0.1}}
\put(0.8,0){\line(1,0){0.1}}
\put(1,0){\line(1,0){0.1}}
\put(1.2,0){\line(1,0){2.8}}
\put(4,-0.1){\line(0,1){0.2}}
\end{picture}
\end{center}
At $\sigma_0$ an initial condition can be prescribed.
The points $\sigma_1,\ldots,\sigma_n$ are inner singularities:
the Hamiltonian function is in the limit point case from both sides.
On the interval $(\sigma_i,\sigma_{i+1})$ a solution of the system behaves according to the
canonical differential equation $y'(x)=zJH_i(x)y(x)$. The data $\oe_i,b_{ij},d_{ij}$ describe what
happens to a solution when passing through the singularity $\sigma_i$.
Thereby, $\oe_i,b_{ij}$ correspond to a point interaction inside the singularity $\sigma_i$,
whereas $d_{ij}$ correspond to a local interaction of the Hamiltonians to the left
and to the right of the singularity $\sigma_i$.
The set $E$ (which we did not indicate in the picture) is used to quantitatively describe the
influence of the interface conditions manifested by the data part $d_{ij}$.
At the points from the set $E$ the interval $(\sigma_0,\sigma_{n+1})$ is split
into smaller pieces that contain at most one singularity.

It can be proved that the spectral theory of singular general Hamiltonians defined
in this way is the full Pontryagin space analogue of the theory of classical Hamiltonians
(being in the limit circle at their left and in the limit point case at their right endpoint).
\begin{enumerate}[$(1)$]
\item With a singular general Hamiltonian $\mf h$ a boundary triple
	$(\mc P(\mf h),T(\mf h),\Gamma(\mf h))$ can be associated and a Weyl coefficient $q_{\mf h}$
	can be constructed; see \cite[Definition~8.5 and Theorem~8.7]{kaltenbaeck.woracek:p4db} and
	\cite[Theorem~5.1 and Definition~5.2]{kaltenbaeck.woracek:p5db}.
\item The Weyl coefficient $q_{\mf h}$ (see \ref{A52} below) belongs to the class $\mc N_{<\infty}$
	and can be interpreted as a $Q$-function of the minimal
	operator $S(\mf h)\defeq T(\mf h)^*$, which is a completely non-self-adjoint symmetry with
	deficiency index $(1,1)$; see \cite[Theorem~8.7]{kaltenbaeck.woracek:p4db} and
	\cite[Proposition~5.19 and Corollary~6.5]{kaltenbaeck.woracek:p5db}.
\item An inverse spectral theorem holds, which states that each generalized Nevanlinna function is the
	Weyl coefficient of some singular general Hamiltonian and that the general Hamiltonian is,
	up to reparameterization, uniquely determined by its Weyl coefficients; see
	\cite[Theorem~1.4]{kaltenbaeck.woracek:p6db} and \cite[Remark~3.38]{kaltenbaeck.woracek:p5db}.

	A local uniqueness theorem holds, which states that beginning sections of Hamiltonians are uniquely
	determined (up to reparameterization) by the asymptotic behaviour of the Weyl function
	towards $i\infty$; see \cite[Theorem~1.2 and Remark~1.3]{langer.woracek:lokinv}.
\end{enumerate}
We should point out that the term `reparameterization', which we used without further notice
in item $(3)$, actually requires some explanation.  Not only the Hamiltonians $H_i$ may be
reparameterized (in the sense of the classical theory, see \eqref{A72})
but also the data $d_{ij}$ and $E$ may be changed according to certain rules;
see \cite[Remark~3.38]{kaltenbaeck.woracek:p5db}.

\begin{nremark}{The fundamental solution}\thlab{A52}
	A crucial concept in the theory of general Hamiltonians is the \emph{fundamental solution}
	associated with a general Hamiltonian $\mf h$. This is the indefinite analogue of the
	fundamental matrix solution of the system \eqref{A30} in the positive definite case.

	It is shown in \cite[\S5a-c]{kaltenbaeck.woracek:p5db} that,
	with a general Hamiltonian $\mf h$,
	a chain $\omega_{\mf h}$ of entire $2\!\times\!2$-matrix functions is associated,
	namely, $\omega_{\mf h}=\omega_{\mf h}(x;z)$,
	$x\in[\sigma_0,\sigma_{n+1})\setminus\{\sigma_1,\ldots,\sigma_n\}$, $z\in\bb C$, so that,
	for fixed $x$, the function $\omega_{\mf h}$ is entire in $z$ and, for fixed $z$, it satisfies
	\begin{equation}\label{A70}
		\frac{\partial}{\partial x}\omega_{\mf h}(x;z)J = z\omega_{\mf h}(x;z)H(x),
		\qquad x\in(\sigma_0,\sigma_{n+1})\setminus\{\sigma_1,\ldots,\sigma_n\}.
	\end{equation}
	Note that the rows of $\omega_{\mf h}$ satisfy the differential equation \eqref{A30}.
	Moreover,
	\begin{equation}\label{A74}
		\omega_{\mf h}(x;0) = I, \quad \det\omega_{\mf h}(x;z)=1, \qquad
		x\in[\sigma_0,\sigma_{n+1})\setminus\{\sigma_1,\ldots,\sigma_n\},\;
		z\in\bb C
	\end{equation}
	and $\omega_{\mf h}(\sigma_0;z)=I$, $z\in\bb C$.
	The chain $\omega_{\mf h}$ is used to construct the Weyl coefficient $q_{\mf h}$
	using a similar limiting procedure as in the positive definite case:
	if $\omega_{\mf h}=(\omega_{\mf h,ij})_{i,j=1}^2$, then
	\[
		q_{\mf h}(z) = \lim_{x\nearrow\sigma_{n+1}}
		\frac{\omega_{\mf h,11}(x;z)\tau+\omega_{\mf h,12}(x;z)}{\omega_{\mf h,21}(x;z)\tau+\omega_{\mf h,22}(x;z)}
	\]
	for $\tau\in\bb R\cup\{\infty\}$ and $z$ in the domain of holomorphy of $q_{\mf h}$
	(which is $\bb C\setminus\bb R$ with at most a finite number of points removed);
	the limit is locally uniform in $z$ and independent of $\tau$;
	see \cite[Lemma~8.2]{kaltenbaeck.woracek:p2db}.

	In the present paper we use some specific properties of fundamental solutions;
	detailed references are provided at the appropriate places.  Here we only would like
	to mention that two general Hamiltonians that are reparameterizations of each other
	give rise to the same fundamental solutions (up to reparameterization in the
	sense of \cite[Definition~3.4]{kaltenbaeck.woracek:p5db}) and hence to the same Weyl coefficients;
	this is shown in \cite[Theorem~1.6]{kaltenbaeck.woracek:p6db}.
\end{nremark}

\begin{nremark}{Splitting of general Hamiltonians}\thlab{A140}
	Let $\mf h$ be a general Hamiltonian, let $s\in\bigcup_{i=0}^n(\sigma_i,\sigma_{i+1})$
	and assume that $s$ is not inner point of an indivisible interval.
	Then `restrictions' of $\mf h$ to the intervals $(\sigma_0,s)$ and $(s,\sigma_{n+1})$
	can be defined, which are denoted by $\mf h_{\Lsh s}$ and $\mf h_{\Rsh s}$,
	respectively; see \cite[Definition~3.47]{kaltenbaeck.woracek:p5db}
	and also \cite[\S2.19]{langer.woracek:gpinf}.
\end{nremark}

\subsection{The class $\mf H_0$}

In the present paper those general Hamiltonians are of interest whose Weyl coefficients
belong to the class $\mc N_{<\infty}^{(\infty)}$.  They can be characterized in a neat way;
see \cite[Theorem~3.1]{langer.woracek:gpinf}.
For the notion of indivisible intervals see Section~\ref{sec-H+M}, \eqref{A136}.

\begin{definition}\thlab{A29}
	We say that a singular general Hamiltonian $\mf h$ belongs to the class $\mf H_0$ if
	\begin{enumerate}[{\rm(i)}]
	\item $\mf h$ has exactly one singularity, i.e.\ $H$ is defined on a set of the form
		$(\sigma_0,\sigma_1)\cup(\sigma_1,\sigma_2)$;
	\item the interval $(\sigma_0,\sigma_1)$ is indivisible of type $0$, i.e.\
		the Hamiltonian function $H_0$ of $\mf h$ on $(\sigma_0,\sigma_1)$ is of the form
		$H_0(x)=h_0(x)\smmatrix{1}{0}{0}{0}$
		with some scalar function $h_0$.
	\end{enumerate}
\vspace*{-3ex}
\end{definition}

\noindent
It follows from the definition of a general Hamiltonian and the form of $H_0$ that
the component $(H_1)_{22}$ is integrable on $(\sigma_1,x_0)$ for some (and hence all)
$x_0\in(\sigma_1,\sigma_2)$, i.e.\ $H_1$ satisfies condition (I) in Definition~\ref{A1}.

\begin{nremark}{The relation $\mf H_0\leftrightsquigarrow\mc N_{<\infty}^{(\infty)}$}\thlab{A51}
	The content of \cite[Theorem~3.1]{langer.woracek:gpinf} is the following:
	a general Hamiltonian $\mf h$ belongs to $\mf H_0$ if and only if its
	Weyl coefficient $q_{\mf h}$ belongs to $\mc N_{<\infty}^{(\infty)}\setminus\mc N_0$.

	Thereby, the negative index of $q_{\mf h}$ can be expressed in terms of the
	general Hamiltonian $\mf h$:
	\begin{equation}\label{A135}
		\ind_-q_{\mf h}=\Delta_1+\left\lfloor\frac{\oe_1}2\right\rfloor+
		\begin{cases}
			1, \quad & \oe_1\text{ odd},\,b_{1,1}>0, \\[0.5ex]
			0, \quad & \text{otherwise}.
		\end{cases}
	\end{equation}
	This is a particular instance of a general formula shown in \cite[Theorem~1.4]{kaltenbaeck.woracek:p6db},
	namely, that $\ind_-q_{\mf h}=\ind_-\mf h$ where $\ind_-\mf h$ is given
	by the formula \cite[(2.13)]{langer.woracek:gpinf}.
\end{nremark}

\noindent
Often it is convenient to use a particular form of a general Hamiltonian
from the class $\mf H_0$, namely
\begin{equation}\label{A113}
\begin{aligned}
	\mf h:\quad & \sigma_0=-1,\; \sigma_1=0,\; \sigma_2=\infty;\qquad
	E=\{-1,x_0,\infty\};
	\\[1ex]
	& H_0(x)=x^{-2}\smmatrix{1}{0}{0}{0},
	\quad x\in(-1,0);\qquad H_1(x),\quad x\in(0,\infty);
	\\[1ex]
	& \oe_1\in\bb N_0,\quad b_{1,1},\ldots,b_{1,\oe_1+1}\in\bb R,
	\quad d_{1,0},\ldots,d_{1,2\Delta_1-1}\in\bb R,
\end{aligned}
\end{equation}
where $x_0\in(0,\infty)$; one can choose $b_{1,\oe_1+1}=0$ if no interval
of the form $(0,\eps)$ with $\eps>0$ is indivisible.
For every given general Hamiltonian $\mf g\in\mf H_0$ and given $x_0\in(0,\infty)$
there exists an $\mf h$ as in \eqref{A113} which is a reparameterization of~$\mf g$.

\subsection{The operator model}

The original definition of the boundary triple associated with a general Hamiltonian given in
\cite{kaltenbaeck.woracek:p4db} is involved and quite abstract (using a completion procedure).
The boundary triple associated with a general Hamiltonian $\mf h$ that has only one
singularity can be described isomorphically in a more concrete way; see
\cite[Definition~2.14 and Theorem~2.15]{langer.woracek:esmod}.
Since we deal with the operator model in some depth, we recall its concrete description
for a general Hamiltonian of the form \eqref{A113}.
In the above mentioned reference the component $h_{11}$ is integrable around $\sigma_1$
instead of the component $h_{22}$.  One has to apply a rotation isomorphism as
defined in \cite[Definition~2.4]{kaltenbaeck.woracek:p5db} and also discussed
in \cite[\S2.g]{langer.woracek:gpinf} to transform the model
from \cite{langer.woracek:esmod} to the current situation.

First we need the following fact, which was shown in
\cite[Lemma~3.10]{kaltenbaeck.woracek:p4db}; also here one has to apply a
rotation isomorphism.

\begin{nremark}{The functions $\mf w_k$}\thlab{A18}
	Let $H\in\bb H$ with $\dom(H)=(a,b)$.  For each $x_0\in(a,b)$,
	there exists a unique sequence $(\mf w_k)_{k\in\bb N_0}$ of absolutely continuous
	real $2$-vector functions on $(a,b)$ such that
	\begin{equation}\label{A114}
	\begin{alignedat}{2}
		&\mf w_0=\binom{1}{0},
			\\[1ex]
		&\mf w_{l+1}'=JH\mf w_l, & & l\geq 0,
			\hspace*{-20ex}
			\\[1ex]
		&\mf w_l|_{(a,x_0)}\in L^2\big(H|_{(a,x_0)}\big),\quad & & l\geq\Delta(H),
			\\[1ex]
		&\mf w_l(x_0)\in\spn\biggl\{\binom{1}{0}\biggr\},\quad & & l\geq 0.
	\end{alignedat}
	\end{equation}
	Let $H_1,H_2\in\bb H$, and assume that $H_1$ and $H_2$ are reparameterizations of each other, say,
	$H_2(x)=H_1(\gamma(x))\gamma'(x)$,
	where $\gamma$ is an increasing bijection such that $\gamma$ and $\gamma^{-1}$ are absolutely continuous.
	Let $x_1\in\dom(H_1)$, set $x_2\defeq\gamma^{-1}(x_1)$
	and let $\mf w_{1;l}$ and $\mf w_{2;l}$
	be the corresponding sequences of functions for $H_1$ and $H_2$ respectively.
	with $x_0$ replaced by $x_1$ and $x_2$ respectively.
	Then, as a simple calculation shows, one has
	\[
		\mf w_{2;l}=\mf w_{1;l}\circ\gamma,\qquad l\geq 0.
	\]
\end{nremark}

\medskip

\noindent
In the following item \ref{A112} we recall the above mentioned isomorphic form of the
operator model.  We restrict ourselves to the case that is needed in the present paper
(this leads to a significant simplification of the formulae).

\begin{nremark}{The boundary triple $(\mc P(\mf h),T(\mf h),\Gamma(\mf h))$}\thlab{A112}
	Let $\mf h\in\mf H_0$ be given by the data as in \eqref{A113},
	and assume, in addition, that $\oe_1=0$ and that no interval $(0,\varepsilon)$
	with $\varepsilon>0$ is indivisible.  Due to the growth
	restriction imposed on the Hamiltonian functions of a general Hamiltonian in its definition (cf.\
	\cite[Definitions~2.16--2.18]{langer.woracek:gpinf}) the Hamiltonian function $H_1$ of
	$\mf h$ satisfies {\rm(I)}, {\rm(HS)} and {\rm($\Delta$)}, i.e.\ it belongs to the class $\bb H$.
	Let $\mf w_l$ be the corresponding functions \eqref{A114}
	and denote by $\mathds{1}_{\Lsh x_0}$ the indicator function of the interval $(0,x_0]$.

	First, we define the base space $\mc P(\mf h)$ of the boundary triple
	$(\mc P(\mf h),T(\mf h),\Gamma(\mf h))$.  Set $\Delta\defeq\Delta(H_1)=\Delta_1$ and
	\[
		L^2_\Delta(H_1)\defeq L^2(H_1)\,\dot+\,
		\spn\big\{\mf w_k\mathds{1}_{\Lsh x_0}\!:\,k=0,\ldots,\Delta-1\big\}.
	\]
	Then $\mc P(\mf h)$ is the linear space
	\[
		\mc P(\mf h)\defeq L^2_\Delta(H_1)\times\bb C^\Delta
	\]
	endowed with an inner product as follows.
	Let $F=(f;\upxi),G=(g;\upeta)\in\mc P(\mf h)$, where
	$\upxi=(\xi_k)_{k=0}^{\Delta-1}$, $\upeta=(\eta_k)_{k=0}^{\Delta-1}$, and denote by
	$\uplambda=(\lambda_k)_{k=0}^{\Delta-1}$ and $\upmu=(\mu_k)_{k=0}^{\Delta-1}$ the
	unique coefficients such that
	\begin{equation}\label{A119}
		\begin{aligned}
			\tilde f &\defeq f-\sum_{l=0}^{\Delta-1}\lambda_l\mf w_l\mathds{1}_{\Lsh x_0}
			\in L^2(H_1), \\[1ex]
			\tilde g &\defeq g-\sum_{l=0}^{\Delta-1}\mu_l\mf w_l\mathds{1}_{\Lsh x_0}
			\in L^2(H_1).
		\end{aligned}
	\end{equation}
	Then
	\[
	    [F,G]\defeq
	    (\tilde f,\tilde g)_{L^2(H_1)}+\sum_{k=0}^{\Delta-1}\lambda_k\overline{\eta_k}
	    +\sum_{k=0}^{\Delta-1}\xi_k\overline{\mu_k}.
	\]
	Second, we define the maximal relation $T(\mf h)$. Set
	\begin{multline*}
		T_{\Delta,\max}(H_1)\defeq \big\{(f;g)\in L^2_\Delta(H_1)\times L^2_\Delta(H_1):\,
		\exists \hat f\text{ absolutely continuous} \\
		\text{representative of $f$ s.t.\ }\hat f'=JH_1g\big\}.
	\end{multline*}
	Then a pair $(F;G)$ of elements $F=(f;\upxi),G=(g;\upeta)\in\mc P(\mf h)$ belongs
	to $T(\mf h)$ if and only if (with $\uplambda$ and $\upmu$ again as in \eqref{A119})
	\begin{enumerate}[(i)]
	\item $(f;g)\in T_{\Delta,\max}(H_1)$;
	\item for each $k\in\{0,\ldots,\Delta-2\}$,\vspace*{-1ex}
		\[
			\xi_k=
			\eta_{k+1}+\frac 12\mu_{\Delta-1}d_{\Delta+k}+
			\frac 12\lambda_0d_k
			-\mf w_{k+1}(x_0)_1f(x_0)_2;
		\]
	\item
		$\displaystyle
			\xi_{\Delta-1}=
			\int\limits_0^{x_0}\mf w_\Delta^*H_1\tilde g+
			\frac 12\sum_{l=0}^{\Delta-1}\lambda_ld_{l+\Delta-1}+
			\mu_{\Delta-1}d_{2\Delta-1}
			-\mf w_\Delta(x_0)_1f(x_0)_2
		$.
	\end{enumerate}
	Here $\mf w_k(x_0)_2$ denotes the lower component of the vector
	$\mf w_k(x_0)$ and
	$f(x_0)=(f(x_0)_1,f(x_0)_2)^T$ denotes the value at $x_0$ of the unique absolutely
	continuous representative $\hat f$ with $\hat f'=JH_1g$ (uniqueness of this representative
	follows since $H_1$ does not end indivisibly towards $0$, cf.\
	\cite[Lemma~3.5]{hassi.snoo.winkler:2000}).

	Finally, we define the boundary relation $\Gamma(\mf h)$: for $(F;G)\in T(\mf h)$, we set
	\[
		\Gamma(\mf h)(F;G)\defeq
		\begin{pmatrix}
			-\lambda_0 \\[2ex]
			\eta_0-f(x_0)_2+\dfrac 12\sum\limits_{l=0}^{\Delta-1}\mu_ld_l
		\end{pmatrix}.
	\]
\end{nremark}

\noindent
The space $\mc P(\mf h)$ and the relation $T(\mf h)$ are related
to the space $L^2(H_1)$ and the maximal relation $T_{\max}(H_1)$ therein as follows.

\begin{nremark}{The map $\psi(\mf h)$}\thlab{A120}
	The original definition of the map $\psi(\mf h)$ that establishes this relation is
	again implicit, cf.\ \cite[Definitions~8.5 and 4.10]{kaltenbaeck.woracek:p4db}.
	However, based on \cite[(4.12)]{kaltenbaeck.woracek:p4db}, it is easy to obtain the
	following description in the concrete model space $\mc P(\mf h)$ introduced above.

	Let $\mf h$ be a general Hamiltonian as in \ref{A112}.
	Then we denote by $\psi(\mf h):\mc P(\mf h)\to L^2_\Delta(H_1)$ the projection onto
	the first component of $\mc P(\mf h)$, i.e.\
	\begin{equation}\label{A71}
		\psi(\mf h)(f;\upxi) \defeq f,\qquad (f;\upxi)\in\mc P(\mf h).
	\end{equation}
	This map satisfies
	\[
		\big(\psi(\mf h)\times\psi(\mf h)\big)T(\mf h)=T_{\Delta,\max}(H_1);
	\]
	see \cite[Remark~2.11]{langer.woracek:esmod}.
	Moreover, it is obvious that $\psi(\mf h)$ maps
	$L^2(H_1)\times\bb C^\Delta$ isometrically and surjectively onto $L^2(H_1)$;
	note that the elements in $\bb C^\Delta$ are neutral.
	Using that
	\[
		L^2(H_1)\times\bb C^\Delta=\big(\{0\}\times\bb C^\Delta\big)^{[\perp]}
	\]
	we can deduce that
	\begin{equation}\label{A56}
	\begin{gathered}
		\psi(\mf h)\big((\{0\}\times\bb C^\Delta)^{[\perp]}\big)=L^2(H_1),
		\\[1ex]
		\big(\psi(\mf h)\times\psi(\mf h)\big)
		\Big(\big(T(\mf h)\cap(\{0\}\times\bb C^\Delta)^{[\perp]}\big)^2\Big)
		=T_{\max}(H_1);
	\end{gathered}
	\end{equation}
	here $[\perp]$ denotes the orthogonal companion with respect to the
	inner product $[\,\cdot,\,\cdot\,]$, i.e.\ $\mc M^{[\perp]}=\{x:[x,y]=0\text{ for all }y\in\mc M\}$.
\end{nremark}

%
%
%
%
%
%
%

%
\subsection{The basic identification}
\label{A73}

The class $\bb H$ of Hamiltonians can be identified with the class
$\mf H_0$ of general Hamiltonians up to the parameters $\oe_1$, $b_{1,j}$, $d_{1,j}$.
This is nearly obvious, but is a crucial observation for our approach.
Hence we point it out in this prominent way.

\begin{nremark}{The relation $\mf H_0\rightsquigarrow\bb H$}\thlab{A109}
	Let $\mf h\in\mf H_0$ be given by the data \eqref{A113}.
	Then $H_1\in\bb H$ and $\Delta(H_1)=\Delta_1$.
\end{nremark}

\begin{nremark}{The relation $\bb H\rightsquigarrow\mf H_0$}\thlab{A108}
	Let $H\in\bb H$ be given, assume that $H$ is defined on $(0,\infty)$ and choose $x_0\in(0,\infty)$.
	Then we associate with $H$ a general Hamiltonian $\mf h\in\mf H_0$ as in \eqref{A113}
	with $H_1=H$ and $\oe_1$, $b_{1,j}$, $d_{1,j}$ arbitrary.
	Again one has $\Delta_1=\Delta(H_1)$, and the negative index of $q_{\mf h}$ is
	given by \eqref{A135}.
\end{nremark}

%
%
\section{Construction of the spectral measure}
\label{sec-measure}
%
%

Let a Hamiltonian $H\in\bb H$ be given. In this section we complete the following tasks:
(1) we show that each solution of the canonical system \eqref{A30}
attains regularized boundary values;
(2) we construct a family of functions from the class $\mc N_{<\infty}^{(\infty)}$
to which we refer as \emph{singular Weyl coefficients} of $H$; and
(3) we construct a positive
Borel measure of class $\bb M$ to which we refer as the \emph{spectral measure} of $H$.
Most of these facts follow relatively easily by using the basic identification \ref{A108}
and previous results from \cite{langer.woracek:gpinf,langer.woracek:ninfrep}.

In order to formulate the theorems, one more notation is needed.

\begin{nremark}{The defect spaces $\mf N_z$}\thlab{A102}
	Let $H\in\bb H$ and $z\in\bb C$.  We denote the set of all locally absolutely continuous
	solutions of the differential equation \eqref{A30} by $\mf N_z$ and speak of the
	\emph{defect space of $H$ at the point $z$}.  Clearly, $\mf N_z$ is a linear space of dimension $2$.
	Note that, for $z=0$, this space is trivial in the sense that
	it consists of all constant functions.

	Let $H_1,H_2\in\bb H$ and assume that $H_1$ and $H_2$ are reparameterizations of each other:
	$H_2(x)=H_1(\gamma(x))\gamma'(x)$
	where $\gamma$ is an increasing bijection such that $\gamma$ and $\gamma^{-1}$ are absolutely continuous.
	Then a simple calculation shows that the mapping $\uppsi\mapsto\uppsi\circ\gamma$ is a bijection between the
	corresponding defect spaces $\mf N_{1,z}$ and $\mf N_{2,z}$.
\end{nremark}

\medskip

\noindent
In the next theorem we show that each solution of \eqref{A30} assumes regularized
boundary values at the left endpoint.
These regularized boundary values will be used later to fix a fundamental system
of solutions.

\begin{theorem}[\textbf{Regularized boundary values}]\thlab{A32}
	Let $H\in\bb H$ with $\dom(H)=(a,b)$.
	Then, for each fixed $x_0\in(a,b)$, the following statements hold {\rm(}the functions
	$\mf w_k$ are as in \ref{A18}{\rm)}.
	\begin{enumerate}[{\rm(i)}]
	\item For each $z\in\bb C$ and each solution $\uppsi=(\uppsi_1,\uppsi_2)^T\in\mf N_z$
		the boundary value
		\[
			\rbvr\uppsi \defeq \lim_{x\searrow a}\uppsi_1(x)
		\]
		and the regularized boundary value
		\begin{equation}
			\rbvs\uppsi
			\defeq
			-\lim_{x\searrow a}\Biggl[
			\sum_{l=0}^{\Delta(H)} z^l\bigl(\mf w_l(x)\bigr)^*J
			\biggl(\uppsi(x)-\lim_{t\searrow a}\uppsi_1(t)
			\hspace*{-1ex}\sum_{k=\Delta(H)+1}^{2\Delta(H)-l}\hspace*{-1ex}z^k\mf w_k(x)\biggr)
			\Biggr]
			\label{A106}
		\end{equation}
		exist.
	\item For $z\in\bb C$ define
		\[
			\rbv:
			\left\{
			\begin{array}{ccl}
				\mf N_z & \to & \bb C^2, \\[1ex]
				\uppsi & \mapsto & (\rbvr\uppsi,\;\rbvs\uppsi)^T.
			\end{array}
			\right.
		\]
		Then $\rbv$ is a bijection from $\mf N_z$ onto $\bb C^2$.
	\item For each $z\in\bb C\setminus\{0\}$ there exists an {\rm(}up to scalar multiples{\rm)} unique
		solution $\uppsi=(\uppsi_1,\uppsi_2)\in\mf N_z\setminus\{0\}$ such that
		$\lim_{x\searrow a}\uppsi_2(x)$ exists.

		This solution is characterized by the property that
		$\uppsi|_{(a,x_0)}\in L^2(H|_{(a,x_0)})$, and also by the property that
		$\rbvr\uppsi=0$ {\rm(}and $\uppsi\neq0${\rm)}.

		If\, $\uppsi$ is such that $\lim_{x\searrow a}\uppsi_2(x)$ exists, then
		\[
			\rbvs\uppsi = \lim_{x\searrow a}\uppsi_2(x).
		\]
	\end{enumerate}
	In contrast to $\rbvr\uppsi$, the regularized boundary value $\rbvs\uppsi$
	depends on the choice of\, $x_0$ since the $\mf w_k$ depend on $x_0$.
	This dependence is controlled as follows.
	\begin{enumerate}[{\rm(i)}]
	\setcounter{enumi}{3}
	\item Let $x_0,\hat x_0\in(a,b)$, and let $\rbv$ and $\rbvpr{z}$
		be the correspondingly defined regularized boundary value mappings.
		Then there exists a polynomial $p(z)$ with real coefficients which has no constant
		term and whose degree does not exceed $2\Delta(H)$ such that
		\[
			\rbvspr{z}\uppsi=\rbvs\uppsi+p(z)\rbvr\uppsi,\qquad
			\uppsi\in\mf N_z,\;\; z\in\bb C.
		\]
	\end{enumerate}
\end{theorem}

\begin{remark}\thlab{A75}
	For $z=0$, solutions $\uppsi$ of \eqref{A30} are constant, and for such $\uppsi$
	the relation
	\[
		\rbv\uppsi = \uppsi(x), \qquad x\in(a,b),
	\]
	holds.
\end{remark}

\begin{proof}[Proof of Theorem~\ref{A32}]
	There is no loss of generality in assuming that $H$ is defined on $(0,\infty)$.
	This follows since the functions $\mf w_k$ transform naturally by composition
	when performing a reparameterization, cf.\ \ref{A18}.

	Let $\mf h$ be the general Hamiltonian given by the data \eqref{A113} with $H_1=H$ as
	in the basic identification \ref{A108} with $\oe_1$, $b_{1,j}$, $d_{1,j}$ all equal to $0$.
	Items (i) and (ii) follow immediately from \cite[Theorem~5.1]{langer.woracek:gpinf};
	we just need to match notation.  Comparing the respective definitions we can deduce that
	\[
		\rbvr\uppsi\,\widehat=\,\Rbv_{\rm r}(z)\uppsi,\quad\rbvs\uppsi\,\widehat=\,-\Rbv_{\rm s}(z)\uppsi,\quad
		\rbv\uppsi\,\widehat=\,\Rbv(z)\uppsi,
	\]
	where the expressions on the right-hand sides are generalized boundary values corresponding
	to the general Hamiltonian $\mf h$ as in \cite{langer.woracek:gpinf}.
	Item (iii) follows directly from \cite[Theorem~5.2]{langer.woracek:gpinf}.
	Only the proof of item (iv) requires an argument.
	
	Let $\wh{\mf h}$ be the general Hamiltonian which is constituted by the same data as $\mf h$
	with the exception that we take $\wh E \defeq \{-1,\hat x_0,\infty\}$ instead of $E$.
	Then $\rbvpr{z}\uppsi=\Rbv(\wh{\mf h},z)\uppsi$, where
	$\Rbv(\wh{\mf h},z)$ denotes the regularized boundary value map defined for $\wh{\mf h}$ as
	in \cite[Theorem~5.1]{langer.woracek:gpinf}.

	By \cite[Proposition~8.11]{kaltenbaeck.woracek:p4db} there exist numbers
	$d_0,\dots,d_{2\Delta(H)-1}\in\bb R$ such that the general Hamiltonian
	$\mf g$ defined as in \eqref{A113} with $H_1=H$, $\oe_1=0$ and $b_{1,1}=0$
	is  a reparameterization of $\wh{\mf h}$.
	Thereby, the increasing bijection between the
	domains of $\wh{\mf h}$ and $\mf g$ is the identity map.  Clearly, the defect spaces $\mf N_z$
	and the functions $\mf w_l$ in \eqref{A114} built with the base point $x_0$ for $\mf h$
	and for $\mf g$, respectively, coincide.

	The fundamental solutions $\omega_{\hat{\mf h}}$ and $\omega_{\mf g}$ coincide; see \ref{A52}.
	Let $\Rbv(\mf g,z)$ be the regularized boundary value map defined for $\mf g$
	and let $\Rbv_{\rm r}(\mf g,z),\Rbv_{\rm s}(\mf g,z)$ be its components
	as in \cite[Theorem~5.1]{langer.woracek:gpinf}.
	By \cite[Remark~5.8]{langer.woracek:gpinf}, equality of
	fundamental solutions implies equality of regularized boundary values, i.e.\
	\[
		\Rbv_{\rm r}(\wh{\mf h},z)=\Rbv_{\rm r}(\mf g,z),\qquad \Rbv_{\rm s}(\wh{\mf h},z)=\Rbv_{\rm s}(\mf g,z).
	\]
	The first of these equalities is of course trivial; both sides, applied to a solution
	$\uppsi=(\uppsi_1,\uppsi_2)^T\in\mf N_z$, are equal to $\lim_{x\searrow a}\uppsi_1(x)$.
	The second equality tells us that
	\[
		\rbvspr{z}\uppsi=-\Rbv_{\rm s}(\mf g,z)\uppsi,\quad \uppsi\in\mf N_z.
	\]
	Comparing the definition of $\Rbv_{\rm s}(\mf g,z)$ in \cite[(5.3)]{langer.woracek:gpinf}
	with the definition of $\rbvs$ in \eqref{A106} we obtain that
	\[
		\Rbv_{\rm s}(\mf g,z)\uppsi=-\rbvs\uppsi+\rbvr\uppsi\cdot
		\sum_{l=1}^{2\Delta(H)}z^ld_{l-1},\quad \uppsi\in\mf N_z,\ z\in\bb C.
	\]
	The assertion in item (iv) thus follows with the polynomial
	\[
		p(z) \defeq -\sum_{l=1}^{2\Delta(H)}z^ld_{l-1}.
		\qedhere
	\]
\end{proof}

\begin{remark}\thlab{A36}
	In the above proof we have defined the general Hamiltonian $\mf h$ via the basic identification using
	$\oe_1$, $b_{1,j}$, $d_{1,j}$ all equal to $0$.  This may seem artificial, and thus requires an explanation.
	To this end, revisit \cite[(5.3)]{langer.woracek:gpinf}.  If we had used other values for
	$\oe_1$, $b_{1,j}$, $d_{1,j}$, then the regularized boundary values of $\mf h$ as defined in
	\cite{langer.woracek:gpinf} would have changed by the summand
	\begin{equation}\label{A54}
		\Bigl(\,\lim_{t\searrow a}\uppsi_1(t)\Bigr)\Biggl(\,\sum_{l=1}^{2\Delta(H)}z^ld_{1,l-1}-
		\sum_{l=0}^{\oe_1}z^{2\Delta(H)+l}b_{1,\oe_1+1-l}\Biggr).
	\end{equation}
	This summand is independent of $x$ and hence contains no information about the asymptotic behaviour
	of $\uppsi$.  We regard the inclusion of a summand \eqref{A54} as a distracting complication from the
	point of our presentation and hence use the choice of vanishing $\oe_1$, $b_{1,j}$, $d_{1,j}$.

	Of course, notions intrinsic for $H$ must not depend on the choice of
	parameters in the basic identification. Thus we shall keep track of
	the influence of $\oe_1$, $b_{1,j}$, $d_{1,j}$.
\end{remark}

\medskip

\noindent
In the next theorem a fundamental system of solutions of \eqref{A30} is constructed.
Since $H$ is not integrable at the left endpoint, this is a non-trivial task.
We fix solutions with the help of the regularized boundary values from Theorem~\ref{A32}.
With this fundamental system of solutions we then construct a singular Weyl coefficient,
which will be used later to obtain a spectral measure.
For the definition of the class $\mc N_{\kappa}^{(\infty)}$ see Definition~\ref{A116}.

\begin{theorem}[\textbf{Singular Weyl coefficients}]\thlab{A94}
	Let $H\in\bb H$ with $\dom(H)=(a,b)$.
	Then, for each fixed $x_0\in(a,b)$, the following statements hold.
	\begin{enumerate}[{\rm(i)}]
	\item For each $z\in\bb C$ denote by $\uptheta(\cdot\,;z)=(\uptheta_1(\cdot\,;z),\uptheta_2(\cdot\,;z))^T$ and
		$\upvarphi(\cdot\,;z)=(\upvarphi_1(\cdot\,;z),\upvarphi_2(\cdot\,;z))^T$ the
		unique elements of\, $\mf N_z$ such that
		\begin{equation}\label{A110}
			\rbv\uptheta(\cdot\,;z)=(1,\,0)^T, \qquad
			\rbv\upvarphi(\cdot\,;z)=(0,\,1)^T.
		\end{equation}
		Then, for each $x\in(a,b)$, the functions $\uptheta(x;\cdot)$ and $\upvarphi(x;\cdot)$
		are entire of finite exponential type\footnote{If the integral in \eqref{A134} is $0$,
		then $\uptheta(x;\cdot)$ and $\upvarphi(x;\cdot)$ are either of minimal exponential type
		or of order less than $1$.}
		\begin{equation}\label{A134}
			\int_a^x \sqrt{\det H(t)}\,\rd t,
		\end{equation}
		and they satisfy $\uptheta_1(x;z)\upvarphi_2(x;z)-\uptheta_2(x;z)\upvarphi_1(x;z)=1$
		for $z\in\bb C$.

		Moreover, let $a_+=\inf\bigl\{x\in(a,b):\int_a^x h_{22}(t)\rd t>0\bigr\}$.
		Then, for each $z\in\bb C$, the following relations hold:
		\begin{equation}\label{A132}
		\begin{alignedat}{2}
			&\lim_{x\searrow a}\uptheta_1(x;z) = 1,
			\qquad &
			&\lim_{x\searrow a}\frac{\uptheta_1(x;z)}{\int_x^{x_0}h_{11}(t)\rd t} = -z,
			\\[1ex]
			&\lim_{x\searrow a_+}\frac{\upvarphi_1(x;z)}{\int_a^x h_{22}(t)\rd t} = -z,
			\hspace*{7ex} &
			&\lim_{x\searrow a}\upvarphi_2(x;z) = 1.
		\end{alignedat}
		\end{equation}
	\item
		For each $\tau\in\bb R\cup\{\infty\}$, the limit
		\begin{equation}\label{A35}
			q_H(z)\defeq\lim_{x\nearrow b}
			\frac{\uptheta_1(x;z)\tau+\uptheta_2(x;z)}{\upvarphi_1(x;z)\tau+\upvarphi_2(x;z)}\,,
			\qquad z\in\bb C\setminus\bb R,
		\end{equation}
		exists locally uniformly on $\bb C\setminus\bb R$, defines an analytic function in $z$
		on $\bb C\setminus\bb R$ and does not depend on $\tau$
		{\rm(}here the fraction on the right-hand side of \eqref{A35} is interpreted
		as $\frac{\uptheta_1(x;z)}{\upvarphi_1(x;z)}$ if $\tau=\infty${\rm)}.
		The function $q_H$ belongs to the class $\mc N_{\Delta(H)}^{(\infty)}$.
	\item We have
		\begin{equation}\label{A23}
			\uptheta(\cdot\,;z)-q_H(z)\upvarphi(\cdot\,;z)\in L^2\bigl(H|_{(x_0,b)}\bigr),
			\qquad z\in\bb C\setminus\bb R,
		\end{equation}
		and this property characterizes the value $q_H(z)$ for each $z\in\bb C\setminus\bb R$.
	\end{enumerate}
	The function $q_H$ depends on the choice of $x_0$, which is controlled as follows.
	\begin{enumerate}[{\rm(i)}]
	\setcounter{enumi}{3}
	\item Let $x_0,\hat x_0\in(a,b)$ and let $q_H$ and $\wh q_H$ be the correspondingly defined
		functions \eqref{A35}.
		Then there exists a polynomial $p$ with real coefficients which has
		no constant term and whose degree does not exceed $2\Delta(H)$ such that
		\[
			\wh q_H(z)=q_H(z)+p(z).
		\]
	\end{enumerate}
\end{theorem}

\begin{proof}
	Again we may assume without loss of generality that $H$ is defined on $(0,\infty)$.
	Let $\mf h$ be the general Hamiltonian \eqref{A113} with $H_1=H$
	and $\oe_1$, $b_{1,j}$, $d_{1,j}$ all equal to $0$,
	so that $\rbv\uppsi=\Rbv(z)\uppsi$, where $\Rbv(z)\uppsi$ is
	as in \cite[Theorem~5.1]{langer.woracek:gpinf}.
	It follows from \cite[Corollary~5.7]{langer.woracek:gpinf} that the
	fundamental solution $\omega_{\mf h}$ from \ref{A52} associated with $\mf h$ is given by
	\begin{equation}\label{A131}
		\omega_{\mf h}(x;z)
		= \begin{pmatrix}
			\uptheta_1(x;z) & \uptheta_2(x;z)
			\\[1ex]
			\upvarphi_1(x;z) & \upvarphi_2(x;z)
		\end{pmatrix}.
	\end{equation}
	The first properties of $\uptheta$ and $\upvarphi$ mentioned in (i) are immediate;
	see \ref{A52}.
	The formula for the exponential type follows from \cite[Theorem~4.1]{langer.woracek:expty}
	if we observe that $\det H_0(x)=0$ for $x\in(-1,0)$ with $H_0$ from \eqref{A113}.
	The limit relations in \eqref{A132} is a consequence of
	\cite[Theorem~4.1 (with $\alpha=0$), Remark~4.2\,(iii) and Lemma~4.14]{langer.woracek:gpinf}.

	Corollary~5.7 in \cite{langer.woracek:gpinf} also implies that the limit in \eqref{A35}
	exists locally uniformly, that $q_H$ is characterized by \eqref{A23} and
	that $q_H$ coincides with the Weyl coefficient $q_{\mf h}$ of the general Hamiltonian $\mf h$.
	In particular, this shows that $q_H$ belongs to the class $\mc N_{\Delta(H)}^{(\infty)}$;
	see \ref{A51} and note that $\Delta(H)=\Delta_1$ and $\oe_1=0$.

	For the proof of item (iv), consider again the general Hamiltonians $\wh{\mf h}$ and
	$\mf g$ as in the proof of Theorem~\ref{A32}\,(iv).  Since they are reparameterizations
	of each other, their Weyl coefficients coincide.
	An application of \cite[Corollary~5.9]{langer.woracek:gpinf} with $\mf h$ and $\mf g$ gives
	\[
		\wh q_H(z)-q_H(z) = q_{\wh{\mf h}}(z)-q_{\mf h}(z)=q_{\mf g}(z)-q_{\mf h}(z)=
		-\sum_{l=1}^{2\Delta(H)}z^ld_{l-1},
	\]
	which shows (iv).
\end{proof}

\noindent
Note that for $z=0$ one has
\begin{equation}\label{A76}
	\uptheta(x;0) = \binom{1}{0}, \quad
	\upvarphi(x;0) = \binom{0}{1}, \qquad
	x\in(a,b),
\end{equation}
which follows from \eqref{A75}.

\begin{remark}\thlab{A96}
	Let us study the influence of the parameters $\oe_1$, $b_{1,j}$, $d_{1,j}$ in \ref{A108} on the
	above proof.  If we chose other parameters than all equal to $0$ and hence
	alter the regularized boundary values by a polynomial summand, then the same would happen to
	the function $q_H$.  In fact, revisiting \cite[Corollary~5.9]{langer.woracek:gpinf} we would
	pass from $q_H$ to $q_H+p$ where $p$ is a polynomial with real coefficients with $p(0)=0$.

	Notice that, conversely, each summand $p\in\bb R[z]$ with $p(0)=0$ can be produced
	by a proper choice of $\oe_1$, $b_{1,j}$, $d_{1,j}$.  Moreover, changing the base point $x_0$
	{\rm(}which is the second arbitrariness in our basic identification{\rm)} also manifests only in
	adding a polynomial summand $p\in\bb R[z]$ with $p(0)=0$.
\end{remark}

\medskip

\noindent
In order to handle the arbitrariness in the basic identification in a structurally clean way, we
introduce an equivalence relation on the set of Weyl coefficients. Namely, we set
\begin{equation}\label{A98}
	q_1\sim q_2\quad \defequiv\quad q_1-q_2\in\bb R[z],\; (q_1-q_2)(0)=0.
\end{equation}
Clearly, this is an equivalence relation on $\mc N_{<\infty}^{(\infty)}$.
In this context, remember that $q\in\mc N_{<\infty}^{(\infty)}$
implies that $q+p\in\mc N_{<\infty}^{(\infty)}$ for all $p\in\bb R[z]$.

\begin{definition}\thlab{A99}
	Let $H\in\bb H$ be given. Then we denote by $[q]_H$ the equivalence class modulo
	the relation \eqref{A98} which contains some (and hence any) function $q_H$ constructed in
	Theorem~\ref{A94}.
	
	We speak of $[q]_H$ as \emph{the singular Weyl coefficient of $H$}.
	Each representative $q_H$ of $[q]_H$ is
	called \emph{a {\rm(!)} singular Weyl coefficient of $H$}.
\end{definition}

\noindent
By this definition we achieve that the singular Weyl coefficient $[q]_H$ of $H\in\bb H$
is nothing but the equivalence class which consists of all Weyl coefficients of general
Hamiltonians associated with $H$ by the basic identification.

In the following theorem a measure is constructed with the help of
the singular Weyl coefficient and a Stieltjes-type inversion formula.

\begin{theorem}[\textbf{The spectral measure}]\thlab{A107}
	Let $H\in\bb H$ be given. Then there exists a unique positive Borel measure $\mu_H$ with
	\begin{equation}\label{A95}
		\mu_H\big([s_1,s_2]\big)=\frac{1}{\pi}\lim_{\varepsilon\searrow 0}
		\lim_{\delta\searrow 0}\intop_{s_1-\varepsilon}^{s_2+\varepsilon}
		\Im q_H(t+i\delta)\,\rd t,\quad -\infty<s_1<s_2<\infty,
	\end{equation}
	where $q_H\in[q]_H$ is any singular Weyl coefficient of\, $H$.
	We have $\mu_H\in\bb M$ and $\Delta(\mu_H)=\Delta(H)$.
\end{theorem}

\begin{proof}
	Since $q_H\in\mc N_{<\infty}^{(\infty)}$, it has a representation $q_H(z)=r+\upphi(\beta_z)$
	with $r\in\bb R$ and a distributional density $\upphi\in\mc F_{\{\infty\}}$,
	i.e.\ $\upphi$ coincides with a measure $\mu_\upphi$ on $\bb R$; see \ref{A50}.
	It follows from \eqref{A142} (see also \cite[Theorem~3.9\,(ii)]{langer.woracek:ninfrep})
	that the measure $\mu_\upphi$ is given by the right-hand side of \eqref{A95}
	on the set of closed intervals; in particular, the double limit exists.
	Moreover, \cite[Theorem~2.8\,(ii)]{langer.woracek:ninfrep} implies that $\mu_H\in\bb M$.

	The fact that $\mu_H$ does not depend on the choice of $q_H\in[q]_H$ is clear since
	a summand that is a real polynomial yields no contribution
	in the Stieltjes inversion formula.

	Finally, we show that $\Delta(\mu_H)=\Delta(H)$.
	It follows from \cite[Theorem~3.9\,(v)]{langer.woracek:ninfrep} that
	\begin{align*}
		\Delta(\mu_H) &= \min\bigl\{\ind_-(q_H+p):p\in\bb R[z]\bigr\} \\[1ex]
		&= \min\bigl\{\ind_-(q_H+p):p\in\bb R[z],\,p(0)=0\bigr\}.
	\end{align*}
	For every $p\in\bb R[z]$ with $p(0)=0$ there exist $\oe_1$, $b_{1,j}$, $d_{1,j}$
	such that the Weyl coefficient of the general Hamiltonian $\mf h$ in \eqref{A113} with $H_1=H$
	is equal to $q_H+p$; and for each choice of $\oe_1$, $b_{1,j}$, $d_{1,j}$ the
	Weyl coefficient is of this form; see \ref{A96}.
	It follows from this and \eqref{A135} that
	\begin{align*}
		\Delta(\mu_H) &= \min\bigl\{\ind_-\mf h: \mf h\text{ as in \eqref{A113} with $H_1=H$}
		\text{ and $\oe_1$, $b_{1,j}$, $d_{1,j}$ arbitrary}\bigr\}
		\\[1ex]
		&= \Delta_1 = \Delta(H).
		\qedhere
	\end{align*}
\end{proof}

\begin{definition}\thlab{A69}
	Let $H\in\bb H$ be given. Then we call the unique positive Borel measure $\mu_H$
	defined by \eqref{A95} the \emph{spectral measure of $H$}.
\end{definition}

\noindent
The choice of this terminology is justified by Theorem~\ref{A4} below where we
construct a Fourier transform into the space $L^2(\mu_H)$.
Before we establish this Fourier transform, let us mention one simple observation.
Namely, it is almost immediate that Hamiltonians which are reparameterizations of each other
give rise to the same singular Weyl coefficients and the same spectral measures.
We provide a slightly more exhaustive variant of this fact.

\begin{proposition}\thlab{A97}
	Let $H\in\bb H$ and $\alpha\in\bb R$. Then the Hamiltonian
	\[
		H_\alpha \defeq
		\begin{pmatrix} 1 & \alpha \\ 0 & 1 \end{pmatrix}H
		\begin{pmatrix} 1 & 0 \\ \alpha & 1 \end{pmatrix}
	\]
	belongs to $\bb H$.
	
	Let, in addition, $\tilde H\in\bb H$ be given and assume that $\tilde H$
	and $H_\alpha$ are reparameterizations of each other.  Then
	\begin{enumerate}[{\rm(i)}]
	\item for each pair of singular Weyl coefficients $q_{\tilde H}$ and $q_H$
		of $\tilde H$ and $H$, respectively, the difference $q_{\tilde H}-q_H$
		is a real polynomial whose constant term equals $\alpha$;
	\item $\mu_{\tilde H}=\mu_H$.
	\end{enumerate}
\end{proposition}

\begin{proof}
	Let $\mf h$ be the general Hamiltonian as in \eqref{A113} with $H_1=H$ and let
	$\omega_{\mf h}$ be its fundamental solution.
	It follows from \cite[Lemma~10.2]{kaltenbaeck.woracek:p2db} that
	\[
		\omega_\alpha \defeq
		\begin{pmatrix} 1 & \alpha \\ 0 & 1 \end{pmatrix}\omega_{\mf h}
		\begin{pmatrix} 1 & -\alpha \\ 0 & 1 \end{pmatrix}
	\]
	is the fundamental solution of some general Hamiltonian $\mf h_\alpha$.
	This factorization immediately yields
	\begin{equation}\label{A53}
		q_{\mf h_\alpha}=q_{\mf h}+\alpha,
	\end{equation}
	which implies that $q_{\mf h_\alpha}\in\mc N_{<\infty}^{(\infty)}$,
	and hence $\mf h_\alpha\in\mf H_0$.
	A short computation, based on \cite[Corollary~5.6]{kaltenbaeck.woracek:p5db},
	shows that the Hamiltonian function of $\mf h_\alpha$ on
	$(0,\infty)$ is equal to $H_\alpha$, and therefore $H_\alpha\in\bb H$.  The functions
	$q_{\mf h_\alpha}$ and $q_{\mf h}$ are singular Weyl coefficients of $H_\alpha$
	and $H$, respectively.  Relation \eqref{A53}, together with Theorem~\ref{A94}\,(iv)
	implies that $q_{H_\alpha}-q_H$ is a real polynomial whose constant term is equal to $\alpha$.

	Assume now that $\tilde H$ is a reparameterizations of $H_\alpha$,
	say \eqref{A72} holds with $H_1=H_\alpha$, $H_2=\tilde H$ and some
	$\gamma:(0,\infty)\to(0,\infty)$.
	We can build general Hamiltonians $\tilde{\mf h}$ and $\mf h_\alpha$ via \ref{A108}
	with $\oe_1=\tilde\oe_1=0$, $b_{1,1}=\tilde b_{1,1}=0$, $d_{1,j}=\tilde d_{1,j}=0$
	and some $x_0,\tilde x_0\in(0,\infty)$ such that $\gamma(\tilde x_0)=x_0$.
	Then $\tilde{\mf h}$ and $\mf h_\alpha$ are reparameterizations of each other;
	see \cite[Remark~3.38]{kaltenbaeck.woracek:p5db}
	and \cite[Proposition~8.13]{kaltenbaeck.woracek:p4db}.
	Hence $\tilde{\mf h}$ and $\mf h_\alpha$ have the same Weyl coefficients;
	see \cite[Theorem~1.4]{kaltenbaeck.woracek:p6db}.
	This shows that $q_{\tilde H}=q_{H_\alpha}$, which in turn implies (i).

	For (ii) observe that an entire summand yields no contribution in the
	Stieltjes inversion formula.
\end{proof}

%
%
\section{The Fourier transform}
\label{sec-fourier}
%
%

For a positive Borel measure $\mu$ on $\bb R$ we denote by $M_{\mu}$
the operator of multiplication by the independent variable in $L^2(\mu)$.
In this section we prove that for each $H\in\bb H$ there exists
a unitary operator $\Theta_H$ from $L^2(H)$ onto $L^2(\mu_H)$, the Fourier transform
connected with $H$, which establishes unitary equivalence of $T(H)$ and $M_{\mu_H}$.
Both $\Theta_H$ and its inverse act as integral transformations.
These results are the most involved ones in the paper.  Their proofs require
to go deeply into the operator model of a general Hamiltonian $\mf h$.
The essential ingredients are the following:
\begin{enumerate}[(1)]
\item
	the spectral theory of the model relation in the Pontryagin space $\mc P(\mf h)$,
	in particular, the spectral decomposition of a self-adjoint relation in a Pontryagin space,
\item
	$Q$-function theory to relate the model relation connected with $\mf h$
	to the model relation of a distributional density,
\item
	the interpretation of a fundamental solution matrix as a generalized $u$-resolvent matrix.
\end{enumerate}

\begin{theorem}[\textbf{The Fourier transform}]\thlab{A4}
	Let $H\in\bb H$ with $\dom(H)=(a,b)$ be given,
	let $T(H)$ be the self-adjoint operator as defined in Section~\ref{sec:canonsys}
	and let $\mu_H$ be the spectral measure associated with $H$ via \eqref{A95}.
	Moreover, let $\upvarphi(\cdot\,;z)=(\upvarphi_1(\cdot\,;z),\upvarphi_2(\cdot\,;z))^T$ be the
	unique element of\, $\mf N_z$ with $\rbv\upvarphi(\cdot\,;z)=(0,\,1)^T$
	as in Theorem~\ref{A94}.
	Then the following statements hold.
	\begin{enumerate}[{\rm(i)}]
	\item
		The map defined by
		\begin{equation}\label{A121}
		\begin{aligned}
			(\Theta_Hf)(t) \defeq \intop_a^b \upvarphi(x;t)^T H(x)f(x)\,\rd x,
			\qquad t\in\bb R, \hspace*{18ex} &
				\\[-1ex]
			f\in L^2(H),\;\; \sup(\supp f)<b, &
		\end{aligned}
		\end{equation}
		extends to an isometric isomorphism from $L^2(H)$ onto $L^2(\mu_H)$,
		which we again denote by $\Theta_H$.
	\item
		The operator $\Theta_H$ establishes a unitary equivalence between $T(H)$ and
		$M_{\mu_H}$, i.e.\ we have
		\[
			\Theta_H\circ T(H) = M_{\mu_H}\circ\Theta_H.
		\]
	\item
		On the subspace of compactly supported functions, also the inverse of $\Theta_H$ acts as an
		integral transformation, namely,
		\begin{equation}\label{A126}
		\begin{aligned}
			(\Theta_H^{-1}g)(x)=\intop_{-\infty}^\infty g(t)\upvarphi(x;t)\,\rd\mu_H(t),
			\qquad x\in(a,b), \hspace*{15ex} &
				\\[-1ex]
			g\in L^2(\mu_H),\;\; \supp g\text{ compact}. &
		\end{aligned}
		\end{equation}
	\end{enumerate}
\end{theorem}

\begin{remark}\thlab{A65}
	Note that the integrals in \eqref{A121} and \eqref{A126} are always well defined.
	For the latter this is obvious since $\upvarphi(x;t)$ is continuous in $t$;
	for the former it follows from Theorem~\ref{A32}\,{\rm(iii)}.
\end{remark}

\medskip

\noindent
As an additional result we prove a connection between the point mass at $0$ of
the spectral measure and the behaviour of $H$ at $b$.

\begin{proposition}\thlab{A77}
	Let $H=(h_{ij})_{i,j=1}^2\in\bb H$, defined on $(a,b)$, and let $\mu_H$ be the
	spectral measure associated with $H$ via \eqref{A95}.
	Then $\mu_H(\{0\})>0$ if and only if
	\begin{equation}\label{A78}
		\intop_a^b h_{22}(x)\rd x < \infty.
	\end{equation}
	If \eqref{A78} is satisfied, then
	\begin{equation}\label{A79}
		\mu_H\big(\{0\}\bigr) = \Biggl[\,\intop_a^b h_{22}(x)\rd x\Biggr]^{-1}.
	\end{equation}
\end{proposition}

\bigskip

\noindent
Note that in any case,
\begin{equation}\label{A130}
	\mu_H\bigl(\{0\}\bigr) = -\lim_{y\searrow0} iyq_H(iy)
\end{equation}
by \cite[(3.8)]{langer.woracek:ninfrep}.

\medskip

The rest of this section is devoted to the proof of Theorem~\ref{A4} and Proposition~\ref{A77}.
We split the proof of Theorem~\ref{A4} into three parts, which are contained in three
separate subsections.
First, in \S\ref{A81}, we construct a Fourier transform $\Theta_H$ from $L^2(H)$ onto
$L^2(\mu_H)$ in an abstract way.  In \S\ref{A82} we show that this map
acts as asserted in \eqref{A121}.  The formula for $\Theta_H^{-1}$ is proved in \S\ref{A83}.
Finally, we prove Proposition~\ref{A77} in \S\ref{A84}.

Since all statements in Theorem~\ref{A4} and Proposition~\ref{A77} are invariant
under reparameterizations,
we can assume without loss of generality that $\dom H=(0,\infty)$.
Throughout these four subsections, keep $H\in\bb H$ fixed and let $\mf h$ be the
general Hamiltonian defined in the basic identification \ref{A108}
with $\oe_1$, $b_{1,j}$, $d_{1,j}$ all equal to $0$.

\subsection{Construction of a Fourier transform}
\label{A81}

Let us first consider the case when $(0,c)$ is a maximal $H$-indivisible interval
of type $0$.  Then the space $L^2(H)$ can be identified with $L^2(H|_{(c,\infty)})$
and $\upvarphi(c;z)=(0,\;1)^T$; see \cite[proof of Theorem~5.1, p.~541]{langer.woracek:gpinf}.
Now items (i) and (ii) in Theorem~\ref{A4} follow from \cite[Theorem~III]{debranges:1961a}.
For the rest of \S\ref{A81} and \S\ref{A82} we assume that
$H$ does not start with an indivisible interval at $0$.
Hence we can use the boundary triple $(\mc P(\mf h),T(\mf h),\Gamma(\mf h))$
associated with the general Hamiltonian $\mf h$ in the form in \ref{A112}.

The Fourier transform $\Theta_H$ is constructed by combining several mappings.
We provide a comprehensive summary in Figure~\ref{A125} on page \pageref{A125}.
The reader may find it helpful to visit this summary already while going through the construction.

It is shown in \cite[Proposition~5.19]{kaltenbaeck.woracek:p5db} that
the Weyl coefficient $q_{\mf h}$ is a $Q$-function of the minimal
relation $S(\mf h) \defeq T(\mf h)^*$.
In fact, denote by $\pi_{l,1}$ the projection from $\bb C^2\times\bb C^2$
onto the upper entry of the first vector
component\footnote{Here `$l,1$' stands for `left vector, first entry'.
This is a generic notation: for example, $\pi_l$ is the projection
from $\bb C^2\times\bb C^2$ onto the first vector component,
$\pi_{r,2}$ onto the lower entry of the second vector component, etc.
The use of `left' and `right' is motivated by the fact that the
vectors correspond to boundary values at the left and right endpoint, respectively.}, set
\begin{equation}\label{A123}
	A\defeq \ker(\pi_{l,1}\circ\Gamma(\mf h)),
\end{equation}
and let $\varepsilon_z$ be the defect elements with
\begin{equation}\label{A25}
	\big(\pi_l\circ\Gamma(\mf h)\big)(\varepsilon_z;z\varepsilon_z)
	= \binom 1{-q_{\mf h}(z)}.
\end{equation}
Then $q_{\mf h}$ is the $Q$-function of $S(\mf h)$ induced by $A$
and $(\varepsilon_z)_{z\in\rho(A)}$.
In particular, since $q_{\mf h}$ is analytic in $\bb C\setminus\bb R$
and $S(\mf h)$ is completely non-self-adjoint, we have $\bb C\setminus\bb R\subseteq\rho(A)$.

Let $\upphi_{\mf h}$ be the distributional density in the representation \eqref{A34}
of $q_{\mf h}$, and let $\Pi(\upphi_{\mf h})$, $A_{\upphi_{\mf h}}$
and $\psi(\upphi_{\mf h})$ be as in \ref{A139}.
By \cite[Proposition~3.4, Proof of Corollary~3.5]{kaltenbaeck.woracek:p2db},
there exists an isometric isomorphism (where we use the base point $z_0=i$)
\begin{equation}\label{A122}
	\Theta_{\mf h}:\mc P(\mf h)\to\Pi(\upphi_{\mf h})\quad\text{with}\quad
	(\Theta_{\mf h}\times\Theta_{\mf h})(A)=A_{\upphi_{\mf h}}.
\end{equation}
This isomorphism is determined by its action on defect elements, namely,
\begin{equation}\label{A49}
	\Theta_{\mf h}(\varepsilon_z) = \hat\varepsilon_z, \qquad z\in\rho(A),
\end{equation}
where $\hat\varepsilon_z\in\Pi(\upphi_{\mf h})$ is defined by
\begin{equation}\label{A145}
	\hat\varepsilon_z(t)
	\defeq \begin{cases}
		\dfrac{t-i}{t-z}\,, &t\in\bb R, \\[2ex]
		1, & t=\infty.
	\end{cases}
\end{equation}
Denote by $\mc E_A(\infty)$ and $\mc E_{\upphi_{\mf h}}(\infty)$ the algebraic eigenspaces
at infinity of the relations $A$ and $A_{\upphi_{\mf h}}$, respectively.
By \cite[Lemma~3.2\,(d)]{langer.woracek:gpinf} we have $\mc E_A(\infty)=\{0\}\times\bb C^\Delta$;
in particular, $\mc E_A(\infty)$ is neutral.
It follows from \eqref{A122} that
\[
	\Theta_{\mf h}(\mc E_A(\infty))=\mc E_{\upphi_{\mf h}}(\infty),
\]
and
\begin{equation}\label{A55}
	\Theta_{\mf h}\bigl(\mc E_A(\infty)^{[\perp]}\bigr)=\mc E_{\upphi_{\mf h}}(\infty)^{[\perp]}.
\end{equation}
These relations imply that also $\mc E_{\upphi_{\mf h}}(\infty)$ is neutral and that
the restriction $\Theta_{\mf h}|_{\mc E_A(\infty)^{[\perp]}}$ is an isometric bijection from
$\mc E_A(\infty)^{[\perp]}$ onto $\mc E_{\upphi_{\mf h}}(\infty)^{[\perp]}$.
Let $\pi_A$ and $\pi_{\upphi_{\mf h}}$ be the following canonical projections:
\begin{align*}
	\pi_A&:\mc E_A(\infty)^{[\perp]}\to
	\raisebox{1pt}{$\mc E_A(\infty)^{[\perp]}$}\big/\raisebox{-2pt}{$\mc E_A(\infty)$}, \\[1ex]
	\pi_{\upphi_{\mf h}}&:\mc E_{\upphi_{\mf h}}(\infty)^{[\perp]}\to
	\raisebox{1pt}{$\mc E_{\upphi_{\mf h}}(\infty)^{[\perp]}$}\big/\raisebox{-2pt}{$\mc E_{\upphi_{\mf h}}(\infty)$}.
\end{align*}
Then there exists an isometric isomorphism
\[
	\Lambda_{\mf h}:\raisebox{0.3ex}{$\mc E_A(\infty)^{[\perp]}$}
	\big/\raisebox{-0.5ex}{$\mc E_A(\infty)$}
	\to
	\raisebox{0.3ex}{$\mc E_{\upphi_{\mf h}}(\infty)^{[\perp]}$}
	\big/\raisebox{-0.5ex}{$\mc E_{\upphi_{\mf h}}(\infty)$}
\]
such that
\[
	\pi_{\upphi_{\mf h}}\circ\Theta_{\mf h}=\Lambda_{\mf h}\circ\pi_A.
\]
By \eqref{A71} we have
\[
	\ker\pi_A=\mc E_A(\infty)=\{0\}\times\bb C^\Delta=\ker\psi(\mf h).
\]
Since $\psi(\mf h)\bigl(\mc E_A(\infty)^{[\perp]}\bigr)=L^2(H)$
by the first relation in \eqref{A56}, we can deduce
that there exists an isometric isomorphism
\[
	\Psi_A:\raisebox{0.3ex}{$\mc E_A(\infty)^{[\perp]}$}\big/
	\raisebox{-0.5ex}{$\mc E_A(\infty)$}\to L^2(H)\quad\text{with}\quad
	\Psi_A\circ\pi_A=\psi(\mf h)|_{\mc E_A(\infty)^{[\perp]}}.
\]
By \cite[Theorem~5.3]{langer.woracek:ninfrep}, $\psi(\upphi_{\mf h})$ maps
$\mc E_{\upphi_{\mf h}}(\infty)^{[\perp]}$ isometrically onto
$L^2\bigl(\frac{\mu_{\upphi_{\mf h}}(x)}{1+x^2}\bigr)$.
Since $L^2\bigl(\frac{\mu_{\upphi_{\mf h}}(x)}{1+x^2}\bigr)$ is non-degenerate, we have
\[
	\ker\psi(\upphi_{\mf h})=\big(\mc E_{\upphi_{\mf h}}(\infty)^{[\perp]}\big)^\circ
	= \mc E_{\upphi_{\mf h}}(\infty)=\ker\pi_{\upphi_{\mf h}},
\]
and we obtain an isometric isomorphism
\[
	\Psi_{\upphi_{\mf h}}:
	\raisebox{0.3ex}{$\mc E_{\upphi_{\mf h}}(\infty)^{[\perp]}$}
	\big/\raisebox{-0.5ex}{$\mc E_{\upphi_{\mf h}}(\infty)$}
	\to L^2\Big(\frac{\mu_{\upphi_{\mf h}}(x) }{1+x^2}\Big)
\]
such that
\[
	\Psi_{\upphi_{\mf h}}\circ\pi_{\upphi_{\mf h}}
	= \psi(\upphi_{\mf h})\big|_{\mc E_{\upphi_{\mf h}}(\infty)^{[\perp]}}.
\]
Finally, let
\[
	U:L^2\Bigl(\frac{\mu_{\upphi_{\mf h}}(x) }{1+x^2}\Bigr)
	\to L^2(\mu_{\upphi_{\mf h}})
\]
be the operator of multiplication by $\frac 1{x-i}$, which is an
isomorphism from $L^2\bigl(\frac{\mu_{\upphi_{\mf h}}(x) }{1+x^2}\bigr)$
onto $L^2(\mu_{\upphi_{\mf h}})$.

Now we define $\Theta_H$ as the composition of the constructed isometric isomorphisms,
namely
\[
	\Theta_H \defeq U\circ\Psi_{\upphi_{\mf h}}\circ\Lambda_{\mf h}\circ\Psi_A^{-1}.
\]
Then $\Theta_H$ is an isometric isomorphism from $L^2(H)$ onto $L^2(\mu_{\upphi_{\mf h}})$.

In order to see how $\Theta_H$ is related to $T(H)$, it is enough to
trace back the defining procedure\footnote{Remember in the following that $T(H)$ is
self-adjoint, and hence $T(H)=T_{\max}(H)=S(H)$ in the notation of several previous papers like
\cite{kaltenbaeck.woracek:p4db}.}.
Using \eqref{A122} and \eqref{A55} we obtain
\[
	(\Lambda_{\mf h}\times\Lambda_{\mf h})
	\Big[(\pi_A\times\pi_A)\big(A\cap(\mc E_A(\infty)^{[\perp]})^2\big)\Big]
	= (\pi_{\upphi_{\mf h}}\times\pi_{\upphi_{\mf h}})
	\big(A_{\upphi_{\mf h}}\cap(\mc E_{\upphi_{\mf h}}(\infty)^{[\perp]})^2\big).
\]
By \cite[Proposition~4.17\,(iii)]{kaltenbaeck.woracek:p4db} we have
\[
	(\psi(\mf h)\times\psi(\mf h))
	\big(T(\mf h)\cap(\mc E_A(\infty)^{[\perp]})^2\big)
	= T(H),
\]
and hence
\[
	(\Psi_A\times\Psi_A)
	\Big[(\pi_A\times\pi_A)\big(A\cap(\mc E_A(\infty)^{[\perp]})^2\big)\Big]
	\subseteq T(H).
\]
It follows from \cite[Theorem~5.3]{langer.woracek:ninfrep} that
\[
	(\psi(\upphi_{\mf h})\times\psi(\upphi_{\mf h}))
	\big(A_{\upphi_{\mf h}}\cap(\mc E_{\upphi_{\mf h}}(\infty)^{[\perp]})^2\big)
	= M_{(1+x^2)^{-1}\rd\mu_{\upphi_{\mf h}}(x)},
\]
and therefore
\[
	(\Psi_{\upphi_{\mf h}}\times\Psi_{\upphi_{\mf h}})
	\Big[(\pi_{\upphi_{\mf h}}\times\pi_{\upphi_{\mf h}})
	\big(A_{\upphi_{\mf h}}\cap(\mc E_{\upphi_{\mf h}}(\infty)^{[\perp]})^2\big)
	\Big]
	= M_{(1+x^2)^{-1}\rd\mu_{\upphi_{\mf h}}(x)}.
\]
Putting these relations together we obtain
\begin{align*}
	\big(\Theta_H^{-1}\times\Theta_H^{-1}\big) M_{\mu_{\upphi_{\mf h}}}
	&= \big(\Psi_A\Lambda_{\mf h}^{-1}\Psi_{\upphi_{\mf h}}^{-1}\times
	\Psi_A\Lambda_{\mf h}^{-1}\Psi_{\upphi_{\mf h}}^{-1}\big)
	M_{(1+x^2)^{-1}\rd\mu_{\upphi_{\mf h}}(x)}
	\\
	&= \big(\Psi_A\Lambda_{\mf h}^{-1}\times\Psi_A\Lambda_{\mf h}^{-1}\big)
	\Big[(\pi_{\upphi_{\mf h}}\times\pi_{\upphi_{\mf h}})
	\big(A_{\upphi_{\mf h}}\cap(\mc E_{\upphi_{\mf h}}(\infty)^{[\perp]})^2\big)\Big]
	\\
	&= (\Psi_A\times\Psi_A)
	\Big[(\pi_A\times\pi_A)\big(A\cap(\mc E_A(\infty)^{[\perp]})^2\big)\Big]
	\\
	&\subseteq T(H).
\end{align*}
Strict inclusion cannot occur since both $M_{\mu_{\upphi_{\mf h}}}$ and $T(H)$
are self-adjoint.  We conclude that
\[
	(\Theta_H\times\Theta_H)T(H) = M_{\mu_{\upphi_{\mf h}}}.
\]
It remains to remember (from the proof of Theorem~\ref{A107}) that $\mu_H=\mu_{\upphi_{\mf h}}$.

%
\begin{figure}
\begin{framed}
\caption{Construction of $\Theta_H$ in \S5.1}
\label{A125}
\hspace*{0pt}\\[-3mm]
\begin{center}
\begin{tabular}{c}
$
	\xymatrix{
	\mc P(\mf h) \ar[r]^{\Theta_{\mf h}}
			\save[]+<0mm,-9mm>*\txt<8pc>{\begin{rotate}{90}$\subseteq$\end{rotate}} \restore
	& \Pi(\upphi_\mf h)
			\save[]+<0mm,-9mm>*\txt<8pc>{\begin{rotate}{90}$\subseteq$\end{rotate}} \restore
	\\
	\mc E_A(\infty)^{[\perp]}
		\ar[d]_{\pi_A} \ar[r]^{\Theta_{\mf h}|_{\mc E_A(\infty)^{[\perp]}}}
		\ar@/_4pc/[dd]_{\psi(\mf h)}
	& \mc E_{\upphi_{\mf h}}(\infty)^{[\perp]}
		\ar[d]^{\pi_{\upphi_{\mf h}}}
		\ar@/^4pc/[dd]^{\psi(\upphi_{\mf h})}
	\\
	{\scriptsize
	\raisebox{1pt}{$\mc E_A(\infty)^{[\perp]}$}\!\big/
		\raisebox{-2pt}{$\mc E_A(\infty)$}
		}
		\ar[r]_{\Lambda_{\mf h}} \ar[d]_{\Psi_A}
	&
	{\scriptsize
	\raisebox{1pt}{$\mc E_{\upphi_{\mf h}}(\infty)^{[\perp]}$}\!\big/
		\raisebox{-2pt}{$\mc E_{\upphi_{\mf h}}(\infty)$}
		}
		\ar[d]^{\Psi_{\upphi_\mf h}}
	\\
	L^2(H)
		\ar@{-->}[rd]_{\Theta_H}
	& L^2\big(\frac{\mu_{\upphi_{\mf h}}(x)}{1+x^2}\big) \ar[d]^{U}
	\\
	& L^2(\mu_{\upphi_{\mf h}})
			\save[]+<22mm,0mm>*\txt<8pc>{$\mu_H=\mu_{\upphi_{\mf h}}$} \restore
	}
$
\end{tabular}
\end{center}
\end{framed}
\end{figure}
%

%
\subsection{Computation of $\bm{\Theta_H}$ as an integral transform}
\label{A82}

Since $\mf h$ starts with an indivisible interval of type $0$,
\cite[Theorem~6.4]{kaltenbaeck.woracek:p5db} is not applicable; a computation of the full
Fourier transform $\Theta_{\mf h}$ in the spirit of this result is not possible. However, we are
only interested in $\Theta_H$, which is essentially a restriction of $\Theta_{\mf h}$.
And it turns out that the action of this restriction can be
computed. The argument is based on a refined variant of
\cite[Proposition~4.6]{kaltenbaeck.woracek:p2db}.
It requires to go into the details of the constructions
made in \cite{kaltenbaeck.woracek:p2db} and \cite{kaltenbaeck.woracek:p5db}.

Let us lay out the operator-theoretic setup (in five parts).  Thereby we use, without
much further notice, terminology and results from \cite{kaltenbaeck.woracek:krall}.
In particular, we ask the reader to recall definitions and usage of
spaces, like $\mc P_-$, dualities $[\,\cdot\,,\cdot\,]_\pm$ and resolvent-like operators $R_z^\pm$,
as introduced and studied in \cite[\S3]{kaltenbaeck.woracek:krall}.
Moreover, we repeatedly employ terminology and results from \cite{kaltenbaeck.woracek:p5db}.
A comprehensive summary of the involved spaces and relations can be found in Figure~\ref{A127}
on page \pageref{A127}.  We advice the reader to visit this summary
already on going through the construction.
\medskip\pagebreak[3]
\begin{center}
	{\footnotesize\ding{80}\, \ding{80}\, \ding{80}}
	\quad\textit{Part 1}\quad
	{\footnotesize\ding{80}\, \ding{80}\, \ding{80}}
\end{center}
\medskip
\noindent
Let $A$ and $\varepsilon_z$ be as in the previous section, cf.\ \eqref{A123}, \eqref{A25},
so that $q_{\mf h}$ is a $Q$-function of $S(\mf h)$ induced by $A$ and
$(\varepsilon_z)_{z\in\rho(A)}$.  Let $u\in\mc P(\mf h)_- \defeq (S(\mf h)^*)'$ be the
element defined by
\[
	\bigl[u,(f;g)\bigr]_\pm
	\defeq \big(\pi_{l,2}\circ\Gamma(\mf h)\big)(f;g),
	\qquad (f;g)\in S(\mf h)^*.
\]
Moreover, set
\[
	R_z \defeq (A-z)^{-1},\qquad
	R_z^+: \left\{
	\begin{array}{rcl}
		\mc P(\mf h) & \to & S(\mf h)^* \\[1ex]
		f & \mapsto & (R_zf;I+zR_zf)
	\end{array}
	\right.
\]
for $z\in\rho(A)$ and let $R_z^-:\mc P(\mf h)_-\to\mc P(\mf h)$ be the dual of $R_{\qu z}^+$, i.e.\
the unique map with
\[
	[R_z^-v,f]=[v,R_{\qu z}^+f]_\pm,\qquad v\in\mc P(\mf h)_-,\, f\in\mc P(\mf h).
\]
Then, by \cite[Theorem~4.24]{kaltenbaeck.woracek:p5db},
we have $\varepsilon_z=R_z^-u$, $z\in\rho(A)$.
\medskip
\begin{center}
	{\footnotesize\ding{80}\, \ding{80}\, \ding{80}}
	\quad\textit{Part 2}\quad
	{\footnotesize\ding{80}\, \ding{80}\, \ding{80}}
\end{center}
\noindent
Since $q_{\mf h}\in\mc N_{<\infty}^{(\infty)}$, the only critical point of $A$ is the point
$\infty$ and $A$ has no finite spectral points of non-positive type,
cf.\ \cite[Lemma~2.5]{langer.woracek:gpinf} and the paragraph preceding it.
This implies that for each bounded Borel set $\Delta$ the spectral projection $E(\Delta)$
of $A$ is well defined and its range is a Hilbert space.
Moreover, $E_\Delta:\Delta'\mapsto E(\Delta\cap\Delta')$ is the spectral measure
of the bounded self-adjoint operator $A_\Delta \defeq A\cap(\ran E(\Delta)\times\ran E(\Delta))$
in the Hilbert space $\ran E(\Delta)$, and $\sigma(A_\Delta)\subseteq\qu{\Delta}$,
cf.\ \cite[Theorem~II.3.1, p.~34]{langer:1982}\footnote{In \cite{langer:1982}
bounded operators are treated.  The extension to the case of relations is
provided in \cite{dijksma.snoo:1987b}.  In \cite{langer:1982} results are
formulated for $\Delta$ being an interval whose endpoints are not critical
points. In our case, the only critical point is $\infty$.
Moreover, with the usual measure-theoretic extension process one can
define $E(\Delta)$ for each bounded Borel set; cf.\ the first paragraph in the proof
of Lemma~\ref{A21} below.  We tacitly use this fact and often formulate results for
bounded Borel sets although in the original references only intervals were considered.}.

For a bounded Borel set $\Delta$ and elements $f,g\in\mc P(\mf h)$ we thus have a complex measure
$E_{\Delta;f,g}$ on $\bb R$ defined by
\[
	E_{\Delta;f,g}(\Delta') \defeq \big[E(\Delta\cap\Delta')f,g\big],\qquad
	\Delta'\text{ a Borel set on $\bb R$}.
\]
For later use let us list some simple properties of these objects.
\begin{enumerate}[$(1)$]
\item  We have
	\begin{equation}\label{A48}
		E_{\Delta;f,g}=E_{\Delta;E(\Delta)f,g}=E_{\Delta;f,E(\Delta)g},\qquad
		f,g\in\mc P(\mf h),
	\end{equation}
	and
	\begin{equation}\label{A64}
		E_{\Delta;f,g}(\Delta')=\qu{E_{\Delta;g,f}(\Delta')},\qquad
		f,g\in\mc P(\mf h),\;\Delta'\;\text{a Borel set on $\qu{\bb R}$}.
	\end{equation}
\item  Denote by $\|\cdot\|_{\mc P(\mf h)}$ some norm which induces the Pontryagin space
	topology of $\mc P(\mf h)$, $\|E(\Delta)\|$ the corresponding operator norm, and let
	$\|E_{\Delta;f,g}\|$ denote the total variation of the measure $E_{\Delta;f,g}$.
	Then
	\[
		\|E_{\Delta;f,g}\|\leq\|E(\Delta)\|\,\|f\|_{\mc P(\mf h)}\|g\|_{\mc P(\mf h)}.
	\]
\item  Let $T>0$, let $F:[-T,T]\to\bb C$ be continuous and let $f,g\in\mc P(\mf h)$.
	Then the map
	\[
		\Delta\mapsto \int_{[-T,T]}F\,\rd E_{\Delta;f,g},\qquad
		\Delta\text{ a Borel subset of }[-T,T],
	\]
	is a complex Borel measure on $[-T,T]$.
\item  For each bounded Borel set $\Delta$ we have
	\[
		E(\Delta)R_z|_{\ran E(\Delta)}
		= (A_\Delta-z)^{-1} = \int_{\bb R}\frac 1{t-z}\,\rd E_\Delta(t),\qquad
		z\in\rho(A).
	\]
\item  If $\Delta,\Delta'$ are bounded Borel sets with $\Delta\subseteq\Delta'$
	and $f,g\in\mc P(\mf h)$, then
	\[
		E_{\Delta;f,g}\ll E_{\Delta';f,g} \quad\text{ with }\quad
		\frac{\rd E_{\Delta;f,g}}{\rd E_{\Delta';f,g}} = \mathds{1}_\Delta.
	\]
\end{enumerate}

\smallskip\pagebreak[3]
\begin{center}
	{\footnotesize\ding{80}\, \ding{80}\, \ding{80}}
	\quad\textit{Part 3}\quad
	{\footnotesize\ding{80}\, \ding{80}\, \ding{80}}
\end{center}
\noindent
Introduce the set
\begin{equation}\label{A89}
	I_\reg \defeq \bigl\{t\in(0,\infty):\text{$t$ is not inner point of
	an indivisible interval}\bigr\}.
\end{equation}
For $t\in I_{\reg}$ the boundary triple
$(\mc P(\mf h),T(\mf h),\Gamma(\mf h))$ is isomorphic to the boundary triple that is obtained by
pasting the boundary triples corresponding to $\mf h_{\Lsh t}$ and $\mf h_{t\Rsh}$;
for the notation of pasting boundary triples and the mentioned result
see \cite[Definition~3.47]{kaltenbaeck.woracek:p5db} and
\cite[Remark~3.51]{kaltenbaeck.woracek:p5db}
(or \cite[2.2, 2.3 and Lemma~2.5]{langer.woracek:esmod});
for the definition of $\mf h_{\Lsh t}$ and $\mf h_{\Rsh t}$ see \ref{A140}.
In particular, the space $\mc P(\mf h)$ can be decomposed as follows:
\[
	\mc P(\mf h)=\mc P(\mf h_{\Lsh t})\,[\dot +]\,\mc P(\mf h_{t\Rsh}),
\]
and we may consider $\mc P(\mf h_{\Lsh t})$ naturally as a subspace of $\mc P(\mf h)$.
Denote by $P_t$, $t\in I_{\reg}$, the orthogonal projection in $\mc P(\mf h)$
onto $\mc P(\mf h_{\Lsh t})$.  Note that this projection acts as
\begin{equation}\label{A128}
	P_t\big((f;\upxi)\big)=(\mathds{1}_{\Lsh t}f;\upxi),\qquad (f;\upxi)\in\mc P(\mf h).
\end{equation}
Let us consider the relation
\[
	S_1(\mf h_{\Lsh t}) \defeq \ker\bigl[(\pi_{l,1}\!\times\!\pi_r)\Gamma(\mf h_{\Lsh t})\bigr]
	\subseteq \mc P(\mf h_{\Lsh t})^2.
\]
This relation is symmetric and has deficiency index $(1,1)$.
Let $\psi_t(z)\in\mc P(\mf h_{\Lsh t})$, $z\in\bb C$, be the defect elements
of $S(\mf h_{\Lsh t})$ that satisfy
\begin{equation}\label{A40}
	\bigl(\pi_l\circ\Gamma(\mf h_{\Lsh t})\bigr)\bigl(\psi_t(z);z\psi_t(z)\bigr)
	= \binom 01,\qquad z\in\bb C;
\end{equation}
then $\psi_t(z)$ are also defect elements of $S_1(\mf h_{\Lsh t})$, i.e.\
\[
	\ran\bigl(S_1(\mf h_{\Lsh t})-\qu z\bigr)^\perp = \spn\bigl\{\psi_t(z)\bigr\},
	\qquad z\in\bb C.
\]
Considering $S_1(\mf h_{\Lsh t})$ as a linear relation in $\mc P(\mf h)$, we have
$S_1(\mf h_{\Lsh t})\subseteq A$ and hence $A\subseteq S_1(\mf h_{\Lsh t})^*$,
where $S_1(\mf h_{\Lsh t})^*$ denotes the adjoint of $S_1(\mf h_{\Lsh t})$ as a
relation in $\mc P(\mf h)$.  The adjoint $S_1(\mf h_{\Lsh t})^{*_t}$ of $S_1(\mf h_{\Lsh t})$
in $\mc P(\mf h_{\Lsh t})$ is given
by $S_1(\mf h_{\Lsh t})^{*_t}=\ker(\pi_{l,1}\circ\Gamma(\mf h_{\Lsh t}))$,
and it follows that
\begin{equation}\label{A124}
	S_1(\mf h_{\Lsh t})^{*_t}=(P_t\times P_t)(A).
\end{equation}
Let $u_t\in\mc P(\mf h_{\Lsh t})_- \defeq (S_1(\mf h_{\Lsh t})^{*_t})'$ be the unique element with
\begin{equation}\label{A41}
	[u_t,(f;g)]_{\pm_t}\defeq \big(\pi_{l,2}\circ\Gamma(\mf h_{\Lsh t})\big)(f;g),\qquad
	(f;g)\in S_1(\mf h_{\Lsh t})^{*_t},
\end{equation}
where $[\,\cdot\,,\cdot\,]_{\pm_t}$ denotes the duality between $S_1(\mf h_{\Lsh t})^{*_t}$
and $\mc P(\mf h_{\Lsh t})_-$.
Since $(P_t\times P_t)(S(\mf h)^*)=S(\mf h_{\Lsh t})^{*_t}$ and
\[
	(\pi_l\circ\Gamma(\mf h))(f;g)=(\pi_l\circ\Gamma(\mf h_{\Lsh t}))(P_tf;P_tg),
	\qquad (f;g)\in S(\mf h)^*,
\]
we have
\begin{equation}\label{A26}
	[u,(f;g)]_\pm=[u_t,(P_tf;P_tg)]_{\pm_t},\qquad (f;g)\in A,
\end{equation}
where $[\,\cdot\,,\cdot\,]_{\pm}$ denotes the duality between $S(\mf h)^*$ and $\mc P(\mf h)_-$.

\medskip\pagebreak[3]
\begin{center}
	{\footnotesize\ding{80}\, \ding{80}\, \ding{80}}
	\quad\textit{Part 4}\quad
	{\footnotesize\ding{80}\, \ding{80}\, \ding{80}}
\end{center}
\noindent
Let $\psi_t$ be as in \eqref{A40}.
For $t\in I_{\reg}$ and $w\in\bb C\setminus\bb R$ define elements
$\psi_{t,w}(z)\in\mc P(\mf h)$ by
\[
	\psi_{t,w}(z) \defeq \bigl(I+(z-w)R_z\bigr)\psi_t(w),\qquad z\in\rho(A);
\]
remember here that $\bb C\setminus\bb R\subseteq\rho(A)$.
It follows from the resolvent identity that
\begin{equation}\label{A42}
	\psi_{t,w}(z_1) = \bigl(I+(z_1-z_2)R_{z_1}\bigr)\psi_{t,w}(z_2),\qquad z_1,z_2\in\rho(A).
\end{equation}
Moreover, since $\psi_t(w)\perp\ran(S_1(\mf h_{\Lsh t})-\qu w)$ and $A$ is a self-adjoint extension
of $S_1(\mf h_{\Lsh t})$, we have
\begin{equation}\label{A37}
	\psi_{t,w}(z)\perp\ran\bigl(S_1(\mf h_{\Lsh t})-\qu z\bigr),\qquad z\in\rho(A).
\end{equation}
Hence there exists a scalar function $\lambda_{t,w}$, which is analytic on $\rho(A)$, such that
\begin{equation}\label{A39}
	P_t\psi_{t,w}(z)=\lambda_{t,w}(z)\cdot\psi_t(z),\qquad z\in\rho(A).
\end{equation}
Clearly, $\lambda_{t,w}(w)=1$, and the zeros of $\lambda_{t,w}$ form a discrete subset
of $\rho(A)$.  Set
\[
	S_{t,w} \defeq \big\{(f;g)\in A:\,g-\qu zf\perp\psi_{t,w}(z),\; z\in\rho(A)\big\}
	\subseteq \mc P(\mf h)^2.
\]
Relation \eqref{A37} implies that $S_{t,w}$ is a symmetric extension of
$S_1(\mf h_{\Lsh t})\subseteq\mc P(\mf h_{\Lsh t})^2$, which acts in
the larger Pontryagin space $\mc P(\mf h)$ and has deficiency index $(1,1)$.
It follows that $(P_t\times P_t)((S_{t,w})^*)\subseteq S_1(\mf h_{\Lsh t})^{*_t}$.
In fact, due to \eqref{A124}, equality must hold.

Set $\mc P(\mf h)^{t,w}_- \defeq (S_{t,w}^*)'$ and let $P_{t,w}:\mc P(\mf h_{\Lsh t})_-\to\mc P(\mf h)^{t,w}_-$
be the adjoint of the map $P_t\times P_t:S_{t,w}^*\to S_1(\mf h_{\Lsh t})^{*_t}$, i.e.\
\[
	\bigl[P_{t,w}v,(f;g)\bigr]_{\pm_{t,w}} = \bigl[v,(P_tf;P_tg)\bigr]_{\pm_t},\qquad
	v\in\mc P(\mf h_{\Lsh t})_-,\; (f;g)\in(S_{t,w})^*,
\]
where $[\,\cdot\,,\cdot\,]_{\pm_{t,w}}$ is the duality between $S_{t,w}^*$ and $\mc P(\mf h)^{t,w}_-$.
Set $u_{t,w} \defeq P_{t,w}u_t$, so that
\begin{equation}\label{A38}
	\bigl[u_{t,w},(f;g)\bigr]_{\pm_{t,w}} = \bigl[u_t,(P_tf;P_tg)\bigr]_{\pm_t},
	\qquad (f;g)\in S_{t,w}^*.
\end{equation}
Remembering \eqref{A26} we have, in particular,
\begin{equation}\label{A24}
	[u_{t,w},(f;g)]_{\pm_{t,w}}=[u,(f;g)]_\pm,\qquad (f;g)\in A.
\end{equation}
Note here that $A$ is contained in both $S(\mf h)^*$ and $S_{t,w}^*$.

Denote by $R_{t,w;z}^-$ the dual of $R_{\qu z}^+:\mc P(\mf h)\to A\subseteq S_{t,w}^*$
corresponding to the duality $[\,\cdot\,,\cdot\,]_{\pm_{t,w}}$. Then, by \eqref{A24}, we have
\[
	[R_{t,w;z}^-u_{t,w},f]=[u_{t,w},R_{\qu z}^+f]_{\pm_{t,w}}
	= [u,R_{\qu z}^+f]_\pm=[R_z^-u,f],\qquad f\in\mc P(\mf h),
\]
and hence
\begin{equation}\label{A17}
	R_{t,w;z}^-u_{t,w} = R_z^-u = \varepsilon_z,\qquad z\in\rho(A).
\end{equation}

\smallskip\pagebreak[3]
\begin{center}
	{\footnotesize\ding{80}\, \ding{80}\, \ding{80}}
	\quad\textit{Part 5}\quad
	{\footnotesize\ding{80}\, \ding{80}\, \ding{80}}
\end{center}
\noindent
For $t\in I_{\reg}$ let $\phi_t(z)$, $z\in\bb C$, be the defect elements
of $S(\mf h_{\Lsh t})$ with
\begin{equation}\label{A146}
	\bigl(\pi_l\circ\Gamma(\mf h_{\Lsh t})\bigr)\bigl(\phi_t(z);z\phi_t(z)\bigr) = \binom 10,
	\qquad z\in\bb C.
\end{equation}
By \cite[Theorem~4.19]{kaltenbaeck.woracek:p5db}, the map $\Xi_t$, which is defined by
\begin{equation}\label{A144}
	(\Xi_t f)(z) = \begin{pmatrix} [f,\phi_t(\qu z)] \\[1ex] [f,\psi_t(\qu z)] \end{pmatrix},\qquad
	f\in\mc P(\mf h_{\Lsh t}),
\end{equation}
is an isomorphism from $\mc P(\mf h_{\Lsh t})$ onto the reproducing kernel space
$\mf K(\omega_{\mf h}(t))$
with the kernel (here $\omega_{\mf h}$ is the fundamental solution of $\mf h$)
\[
	H_{\omega_{\mf h}(t)}(w,z) = \frac{\omega_{\mf h}(t;z)J\omega_{\mf h}(t;w)^*-J}{z-\qu w}\,.
\]
By \cite[Proposition~4.4]{kaltenbaeck.woracek:p5db} the kernel $H_{\omega_{\mf h}(t)}$
can be written as
\[
	H_{\omega_{\mf h}(t)}(w,z) =
	\begin{pmatrix}
		[\phi_t(z),\phi_t(w)] & [\phi_t(z),\psi_t(w)] \\[1ex]
		[\psi_t(z),\phi_t(w)] & [\psi_t(z),\psi_t(w)]
	\end{pmatrix}.
\]
On the dense subspace $\spn\big(\{\phi_t(z):\,z\in\bb C\}\cup\{\psi_t(z):\,z\in\bb C\}\big)$ of
$\mc P(\mf h_{\Lsh t})$, the action of $\Xi_t$ is determined by linearity and the formula
\begin{equation}\label{A47}
	\Xi_t\bigl(\lambda\phi_t(z)+\mu\psi_t(z)\bigr)
	= H_{\omega_{\mf h}(t)}(\qu z,\cdot)\binom\lambda\mu,\qquad
	\lambda,\mu\in\bb C,\; z\in\bb C.
\end{equation}
This relation is seen as follows ($\zeta\in\bb C$):
\begin{multline*}
	\Xi_t\bigl(\lambda\phi_t(z)+\mu\psi_t(z)\bigr)(\zeta)
	= \lambda\begin{pmatrix}
		[\phi_t(z),\phi_t(\qu\zeta)] \\[1ex]
		[\phi_t(z),\psi_t(\qu\zeta)]
	\end{pmatrix}
	+ \mu\begin{pmatrix}
		[\psi_t(z),\phi_t(\qu\zeta)] \\[1ex]
		[\psi_t(z),\psi_t(\qu\zeta)]
	\end{pmatrix}
	\\[1ex]
	= \bigl(H_{\omega_{\mf h}(t)}(\qu\zeta,z)\bigr)^T\binom{\lambda}{\mu}
	= \bigl(H_{\omega_{\mf h}(t)}(\zeta,\qu z)\bigr)^*\binom{\lambda}{\mu}
	= H_{\omega_{\mf h}(t)}(\qu z,\zeta)\binom{\lambda}{\mu}.
\end{multline*}

%
\begin{figure}
\begin{framed}
\caption{Computation of $\Theta_H$ in \S5.2}
\label{A127}
\hspace*{0pt}\\[-3mm]
\begin{center}
\begin{tabular}{l@{\qquad\quad}l}
	$\mc P(\mf h)=\mc P(\mf h_{\Lsh t})[\dot+]\mc P(\mf h_{t\Rsh})$
	& $\pi_l\circ\Gamma(\mf h)=\pi_l\circ\Gamma(\mf h_{\Lsh t})$
	\\
	\hspace*{-3pt}$\xymatrix{\mc P(\mf h) \ar[r]^{P_t} & \mc P(\mf h_{\Lsh t})}$
	& {\footnotesize $[u,\cdot]_\pm\defeq\pi_{l,2}\circ\Gamma(\mf h)$,\;\;
	$[u_t,\cdot]_{\pm_t}\defeq\pi_{l,2}\circ\Gamma(\mf h_{\Lsh t})$}
	\\
	& {\footnotesize
	$[u_{t,w},\cdot]_{\pm_{t,w}}\defeq\pi_{l,2}\circ\Gamma(\mf h_{\Lsh t})\circ(P_t\times P_t)$}
	\\[2mm] \hline &\\[-2mm]
	$S(\mf h)\defeq \ker\big[\pi_l\circ\Gamma(\mf h)\big]$
	&
	\parbox[t]{50mm}{
	$S(\mf h_{\Lsh t})\defeq\ker\big[(\pi_l\times\pi_r)\circ\Gamma(\mf h_{\Lsh t})\big]$\\
		{\small(\scriptsize defect elements $\psi_t(z),\phi_t(z)$\small)}
		}
	\\[6mm]
	$A\defeq \ker\big[\pi_{l,1}\circ\Gamma(\mf h)\big]$
	&
	\parbox[t]{52mm}{
		$S_1(\mf h_{\Lsh t})\defeq
		\ker\big[(\pi_{l,1}\times\pi_r)\circ\Gamma(\mf h_{\Lsh t})\big]$\\
		{\small(\scriptsize defect elements $\psi_t(z)$\small)}
		}
	\\
	\parbox{34mm}{$S_{t,w}$\\
		\big(\,\parbox{32mm}{\scriptsize defect elements $\psi_{t,w}(z)$,\\
		\scriptsize $P_t\psi_{t,w}(z)=\lambda_{t,w}(z)\psi_t(z)$}\big)
		}
	&
	\\[5mm] \hline &\\[-4mm]
	\multicolumn{2}{c}{\hspace*{-30mm}\parbox{97mm}{$
	\xymatrix@R=2mm@C=3mm{
		&&&&&&&&&
		\\
		&&&
		&&&&& S(\mf h_{\Lsh t})^{*_t}
		&
		\\
		&
		S(\mf h)^* && S_{t,w}^* \ar@{-->}[rrrrr]^{P_t\times P_t}
		&&&&& S_1(\mf h_{\Lsh t})^{*_t}
		\save[]+<0mm,2.3mm>*\txt<8pc>{\begin{rotate}{90}$\subseteq$\end{rotate}} \restore
		&
		\\
		&& \hspace*{-2pt}A \ar@{-->}@/_1pc/[urrrrrr]^{P_t\times P_t}
		\save[]+<-7mm,3mm>*\txt<8pc>{\begin{rotate}{-40}$\supseteq$\end{rotate}} \restore
		\save[]+<2.5mm,-3.5mm>*\txt<8pc>{\begin{rotate}{-40}$\supseteq$\end{rotate}} \restore
		\save[]+<3.5mm,1.5mm>*\txt<8pc>{\begin{rotate}{40}$\subseteq$\end{rotate}} \restore
		\save[]+<-5.5mm,-5.5mm>*\txt<8pc>{\begin{rotate}{40}$\subseteq$\end{rotate}} \restore
		& &&&&& \raisebox{-3pt}{\begin{rotate}{90}$\subseteq$\end{rotate}}
		&
		\\
		&
		S(\mf h) && S_{t,w}
		&&&&& S_1(\mf h_{\Lsh t})\ar@{-->}[lllll]^{\supseteq}
		&
		\\
		&
		\save[]+<-18mm,7.5mm>*\txt<8pc>{\begin{rotate}{90}{\scriptsize
		co-dimension}\end{rotate}} \restore
		\save[]+<-14mm,22mm>*\txt<8pc>{
			\begin{rotate}{-90}
			\rule{0.6pt}{4pt}%
				\rule{5pt}{0.4pt}\hspace*{4pt}
				\raisebox{-3pt}{\begin{rotate}{90}{\scriptsize 1}\end{rotate}}
				\hspace*{-1pt}\rule{5pt}{0.4pt}%
			\rule{0.6pt}{4pt}%
				\rule{5pt}{0.4pt}\hspace*{4pt}
				\raisebox{-3pt}{\begin{rotate}{90}{\scriptsize 1}\end{rotate}}
				\hspace*{-1pt}\rule{5pt}{0.4pt}%
			\rule{0.6pt}{4pt}%
			\end{rotate}
		} \restore
		&&
		&&&&& S(\mf h_{\Lsh t}) \ar@{-->}@/^1pc/[ulllllll]^{\supseteq}
		\save[]+<0mm,2.3mm>*\txt<8pc>{\begin{rotate}{90}$\subseteq$\end{rotate}} \restore
		\save[]+<20mm,22mm>*\txt<8pc>{\begin{rotate}{-90}{\scriptsize
			co-dimension}\end{rotate}} \restore
		\save[]+<15mm,0mm>*\txt<8pc>{
			\begin{rotate}{90}
			\rule{0.6pt}{4pt}%
				\rule{5pt}{0.4pt}\hspace*{-0.5pt}
				\raisebox{1pt}{\begin{rotate}{-90}{\scriptsize 1}\end{rotate}}
				\hspace*{3.5pt}\rule{5pt}{0.4pt}%
			\rule{0.6pt}{4pt}%
				\rule{12pt}{0.4pt}\hspace*{3pt}
				\raisebox{0pt}{\begin{rotate}{-90}{\scriptsize 2}\end{rotate}}
				\hspace*{7pt}\rule{12pt}{0.4pt}%
			\rule{0.6pt}{4pt}%
				\rule{6pt}{0.4pt}\hspace*{-0.5pt}
				\raisebox{1pt}{\begin{rotate}{-90}{\scriptsize 1}\end{rotate}}
				\hspace*{3.5pt}\rule{6pt}{0.4pt}%
			\rule{0.6pt}{4pt}%
			\end{rotate}
		} \restore
		&
		\\
		\ar@{..}[uuuuuu] &&
		\save[]+<0mm,-4mm>*\txt<8pc>{$\underbrace{\rule{30mm}{0pt}}_{\text{in
			$\mc P(\mf h)^2$}}$} \restore
		&&&&&&
		\save[]+<0mm,-4mm>*\txt<8pc>{$\underbrace{\rule{10mm}{0pt}}_{\text{in
			$\mc P(\mf h_{\Lsh t})^2$}}$} \restore
		& \ar@{..}[uuuuuu]
	}
	$}}
	\\[27mm] \hline &\\[-2mm]
\multicolumn{2}{c}{$
\xymatrix{
	& \bb C &&&
	\\
	S(\mf h)^* \ar[ru]^{[u,\cdot]_\pm}
	& S_{t,w}^* \ar[rrr]_{P_t\times P_t} \ar[u]_{[u_{t,w},\cdot]_{\pm_{t,w}}}
	&&& S_1(\mf h_{\Lsh t})^* \ar[lllu]_{[u_t,\cdot]_{\pm_t}}
	\\
	& A \ar[rrru]_{P_t\times P_t}
		\save[]+<1mm,5mm>*\txt<8pc>{\begin{rotate}{90}$\subseteq$\end{rotate}} \restore
		\save[]+<-11mm,6.5mm>*\txt<8pc>{\begin{rotate}{-45}$\supseteq$\end{rotate}} \restore
	&&&
}
$}
\\[-4mm]
\end{tabular}
\end{center}
\end{framed}
\end{figure}
%

\hspace*{1ex}

\begin{center}
	{\footnotesize\ding{80}\qquad \ding{80}\qquad \ding{80}\qquad \ding{80}\qquad
	\ding{80}\qquad \ding{80}\qquad \ding{80}\qquad \ding{80}\qquad \ding{80}}
\end{center}

\hspace*{1ex}

\noindent
In order to prove \eqref{A121}, we proceed in several steps and compute various
inner products and the action of several maps.
The desired result will then follow by putting these formulae together.
A comprehensive overview of the involved maps is provided in the diagram
on page \pageref{A57}.

For $t\in I_{\reg}$ denote by $\Theta_t$ the map $\Xi_t$ followed by projection onto the second component, i.e.\
\begin{equation}\label{A88}
	(\Theta_tf)(z) \defeq [f,\psi_t(\qu z)],\qquad f\in\mc P(\mf h_{\Lsh t}).
\end{equation}
Recall that the entries of $\omega_{\mf h}(t;z)$ are by definition (see
\cite[Definitions~5.3, 4.3]{kaltenbaeck.woracek:p5db}) just the right-hand boundary values
of the elements $\phi_t(z)$ and $\psi_t(z)$:
\begin{align*}
	\bigl(\pi_r\circ\Gamma(\mf h_{\Lsh t})\bigr)\bigl(\phi_t(z);z\phi_t(z)\bigr)
	&= \omega_{\mf h}(t;z)^T\binom 10,
	\\[1ex]
	\bigl(\pi_r\circ\Gamma(\mf h_{\Lsh t})\bigr)\bigl(\psi_t(z);z\psi_t(z)\bigr)
	&= \omega_{\mf h}(t;z)^T\binom 01.
\end{align*}
With the same method that was used in the proof of \cite[Proposition~4.1]{kaltenbaeck.woracek:p2db}
we obtain the following statement.

\begin{lemma}\thlab{A21}
	Let $\Delta$ be a bounded Borel set, $z_0\in\bb C\setminus\bb R$ and $t\in I_{\reg}$. Then
	\begin{equation}\label{A61}
	\begin{aligned}
		\bigl[E(\Delta)f,\psi_{t,w}(\zeta)\bigr]
			= \int\limits_{-\infty}^\infty(\Theta_tf)(x)\cdot(x-z_0)\,
			\rd E_{\Delta;\varepsilon_{z_0},\psi_{t,w}(\zeta)}(x),
			\hspace*{12ex} &
			\\[-1ex]
		f\in\mc P(\mf h_{\Lsh t}),\;\zeta,w\in\bb C\setminus\bb R. &
	\end{aligned}
	\end{equation}
\end{lemma}

\begin{proof}
	The algebra of all bounded Borel sets is the union $\bigcup_{T>0}\mc B_{[-T,T]}$
	of the $\sigma$-algebras $\mc B_{[-T,T]}$ of all Borel subsets of $[-T,T]$.
	Since $\mc B_{[-T,T]}$, as a $\sigma$-algebra, is generated by the set of all closed
	intervals whose endpoints do not carry a point mass of $E$, it is enough to establish \eqref{A61}
	for such intervals.
	
	Throughout the proof fix $f\in\mc P(\mf h_{\Lsh t})$ and an interval
	$[a_-,a_+]$ with $E(\{a_-\})=E(\{a_+\})=0$. Moreover, choose a bounded open interval $\Delta_0$
	which contains $[a_-,a_+]$.
	
	\bigskip

	\noindent\textit{Step 1: some computations.}
	To start with, we compute (indicating which equations are used)
	\begin{align*}
		\bigl[u_{t,w},\bigl(\psi_{t,w}(\qu z);\qu z\psi_{t,w}(\qu z)\bigr)\bigr]_{\pm_{t,w}}&
			\stackrel{(\ref{A38})}{=}\hspace*{2ex}
			\bigl[u_t,\bigl(P_t\psi_{t,w}(\qu z);\qu zP_t\psi_{t,w}(\qu z)\bigr)\bigr]_{\pm_t}
			\\
		&\stackrel{(\ref{A39})}{=}\hspace*{2ex}
			\qu{\lambda_{t,w}(\qu z)}\bigl[u_t,\bigl(\psi_t(\qu z);\qu z\psi_t(\qu z)\bigr)\bigr]_{\pm_t}
			\\
		&\hspace*{-2ex}\stackrel{(\ref{A40},\,\ref{A41})}{=}
			\qu{\lambda_{t,w}(\qu z)},
	\end{align*}
	which yields
	\begin{equation}\label{A10}
		[f,\psi_{t,w}(\qu z)]
		\stackrel{(\ref{A39})}{=} \qu{\lambda_{t,w}(\qu z)}\,[f,\psi_t(\qu z)]
		= (\Theta_tf)(z)\cdot
		\bigl[u_{t,w},\bigl(\psi_{t,w}(\qu z);\qu z\psi_{t,w}(\qu z)\bigr)\bigr]_{\pm_{t,w}}.
	\end{equation}
	Since
	\[
		(\qu z-\qu z_0)R_{\qu z_0}\psi_{t,w}(\qu z)
		\stackrel{(\ref{A42})}{=}
		\psi_{t,w}(\qu z)-\psi_{t,w}(\qu z_0),
	\]
	we have
	\begin{align*}
		(z&-z_0)[\varepsilon_{z_0},\psi_{t,w}(\qu z)]
		\stackrel{(\ref{A17})}{=}
		(z-z_0)\bigl[R_{t,w;z_0}^-u_{t,w},\psi_{t,w}(\qu z)\bigr]
			\\[1ex]
		&= \bigl[u_{t,w},(\qu z-\qu z_0)R_{\qu{z_0}}^+\psi_{t,w}(\qu z)\bigr]_{\pm_{t,w}}
			\\[1ex]
		&= \Bigl[u_{t,w},\Bigl((\qu z-\qu z_0)R_{\qu z_0}\psi_{t,w}(\qu z);
			(\qu z-\qu z_0)\bigl(I+\qu z_0 R_{\qu z_0}\bigr)\psi_{t,w}(\qu z)\Bigr)\Bigr]_{\pm_{t,w}}
			\displaybreak[0]\\[1ex]
		&= \Bigl[u_{t,w},\Bigl(\psi_{t,w}(\qu z)-\psi_{t,w}(\qu z_0);(\qu z-\qu z_0)\psi_{t,w}(\qu z)
			+\qu z_0\bigl(\psi_{t,w}(\qu z)-\psi_{t,w}(\qu z_0)\bigr)\Bigr)\Bigr]_{\pm_{t,w}}
			\\[1ex]
		&= \bigl[u_{t,w},(\psi_{t,w}(\qu z);\qu z\psi_{t,w}(\qu z))\bigr]_{\pm_{t,w}}-
			\bigl[u_{t,w},(\psi_{t,w}(\qu{z_0});\qu{z_0}\psi_{t,w}(\qu{z_0}))\bigr]_{\pm_{t,w}}.
	\end{align*}
	Bringing the very last summand to the left-hand side and substituting the first term on
	the right-hand side into \eqref{A10} we obtain
	\begin{equation}\label{A11}
	\begin{aligned} {}
		[f,\psi_{t,w}(\qu z)]
			= (\Theta_tf)(z)\cdot\Bigl(&(z-z_0)\bigl[\varepsilon_{z_0},\psi_{t,w}(\qu z)\bigr]
			\\[1ex]
		&+\bigl[u_{t,w},\bigl(\psi_{t,w}(\qu{z_0});\qu{z_0}\psi_{t,w}(\qu{z_0})\bigr)\bigr]_{\pm_{t,w}}\Bigr).
	\end{aligned}
	\end{equation}
	For each $z,\zeta\in\bb C\setminus\bb R$, we have
	\begin{equation}\label{A12}
		\big[(I+(z-\qu\zeta)R_z)f,\psi_{t,w}(\zeta)\big]=\big[f,(I+(\qu z-\zeta)R_{\qu z})\psi_{t,w}(\zeta)\big]
		\stackrel{(\ref{A42})}{=}
		[f,\psi_{t,w}(\qu z)]
	\end{equation}
	and
	\begin{equation}\label{A13}
	\begin{aligned}
		\bigl[\varepsilon_{z_0},\psi_{t,w}(\qu z)\bigr]\hspace*{3pt}
		\hspace*{3pt}&\hspace*{-6.5pt}\stackrel{(\ref{A42})}{=}
			\bigl[\varepsilon_{z_0},(I+(\qu z-\zeta)R_{\qu z})\psi_{t,w}(\zeta)\bigr]
			\\[1ex]
		&=\hspace*{1.5ex} \bigl[(I+(z-\qu \zeta)R_z)\varepsilon_{z_0},\psi_{t,w}(\zeta)\bigr]
			\\[1ex]
		&=\hspace*{1.5ex} \bigl[\varepsilon_{z_0},\psi_{t,w}(\zeta)\bigr]
			+(z-\qu \zeta)\bigl[R_zE(\Delta_0)\varepsilon_{z_0},\psi_{t,w}(\zeta)\bigr]
			\\[1ex]
		&\hspace*{1.5ex}\quad+(z-\qu \zeta)\bigl[R_zE(\Delta_0^c)\varepsilon_{z_0},\psi_{t,w}(\zeta)\bigr].
	\end{aligned}
	\end{equation}
	If $z\neq\qu\zeta$, then
	\begin{align}
		&\frac{1}{z-\qu\zeta}[f,\psi_{t,w}(\zeta)]+[R_zf,\psi_{t,w}(\zeta)]
			\stackrel{(\ref{A12})}{=}
			\frac{1}{z-\qu\zeta}[f,\psi_{t,w}(\qu z)]
			\label{A43}\\[1ex]
		&\stackrel{(\ref{A11})}{=}
			\frac{1}{z-\qu\zeta}(\Theta_t f)(z)\cdot\Bigl((z-z_0)\bigl[\varepsilon_{z_0},\psi_{t,w}(\qu z)\bigr]
			+\bigl[u_{t,w},(\psi_{t,w}\bigl(\qu{z_0});\qu{z_0}\psi_{t,w}(\qu{z_0})\bigr)\bigr]_{\pm_{t,w}}\Bigr)
			\notag\displaybreak[0]\\[1ex]
		&\stackrel{(\ref{A13})}{=}
			(\Theta_t f)(z)\cdot\biggl(\frac{z-z_0}{z-\qu\zeta}[\varepsilon_{z_0},\psi_{t,w}(\zeta)]
			\label{A44}\\[1ex]
		&\qquad +(z-z_0)\bigl[R_zE(\Delta_0)\varepsilon_{z_0},\psi_{t,w}(\zeta)\bigr]
		+(z-z_0)\bigl[R_zE(\Delta_0^c)\varepsilon_{z_0},\psi_{t,w}(\zeta)\bigr]
			\label{A45}\\[1ex]
		&\qquad +\frac{1}{z-\qu\zeta}\bigl[u_{t,w},
			\bigl(\psi_{t,w}(\qu{z_0});\qu{z_0}\psi_{t,w}(\qu{z_0})\bigr)\bigr]_{\pm_{t,w}}\biggr).
			\label{A46}
	\end{align}
	
	\bigskip

	\noindent\textit{Step 2: use of the Stieltjes--Lif\v{s}ic inversion formula.}
	We shall apply the Stieltjes--Lif\v{s}ic inversion formula as stated
	in \cite[p.~24, Corollary~II.2 (second formula)]{langer:1982}
	with two minor modifications.
	First, since we assume that the endpoints $a_+$ and $a_-$ carry no point mass of $E$,
	the limit `{\small$\varepsilon\searrow 0$}' is not needed.
	Second, by analyticity we may apply Cauchy's theorem to replace
	the path of integration used in \cite{langer:1982} by the path $\gamma_\delta$
	composed of the two oriented line segments connecting the points
	$a_--i\delta$, $a_+-i\delta$, and $a_++i\delta$, $a_-+i\delta$, respectively.
	Then, for each function $g$ that is analytic in some open neighbourhood of $[a_-,a_+]$ and
	$u,v\in\ran E(\Delta_0)$,
	\begin{equation}\label{A14}
		\lim_{\delta\searrow0}\frac{-1}{2\pi i}
		\int\limits_{\gamma_\delta} g(z)[R_z v,w]\rd z
		= \int_{[a_-,a_+]} g(x)\,\rd E_{\Delta_0;v,w}(x).
	\end{equation}
	Since $\Theta_t f$ is entire and $\zeta\notin\bb R$, the first term on the left-hand side of \eqref{A43},
	the terms in \eqref{A44} and \eqref{A46} and the second term in \eqref{A45} are holomorphic
	in a neighbourhood of $[a_-,a_+]$.  Hence
	\begin{equation}\label{A62}
		\lim_{\delta\searrow0}
		\int\limits_{\gamma_\delta} [R_z f,\psi_{t,w}(\zeta)]\rd z
		= \lim_{\delta\searrow0}
		\int\limits_{\gamma_\delta} (\Theta_t f)(z)\cdot(z-z_0)
		\bigl[R_zE(\Delta_0)\varepsilon_{z_0},\psi_{t,w}(\zeta)\bigr]\rd z.
	\end{equation}
	Applying \eqref{A14} twice, namely with the entire functions $g(z)=1$ and
	$g(z)=(\Theta_t f)(z)\cdot(z-z_0)$, we obtain
	\begin{align*}
		\bigl[E(&[a_-,a_+])f,\psi_{t,w}(\zeta)\bigr]
			\\[1ex]
		&=\bigl[E([a_-,a_+])E(\Delta_0)f,E(\Delta_0)\psi_{t,w}(\zeta)\bigr]
		=\hspace*{-1ex}\intop_{[a_-,a_+]}\hspace*{-10pt} \rd E_{\Delta_0;E(\Delta_0)f,E(\Delta_0)\psi_{t,w}(\zeta)}(x)
			\\[1ex]
		&\hspace*{-1.3ex}\stackrel{(\ref{A14})}{=}
			\lim_{\delta\searrow 0}\frac{-1}{2\pi i}
			\int\limits_{\gamma_\delta} \bigl[R_zE(\Delta_0)f,E(\Delta_0)\psi_{t,w}(\zeta)\bigr]\,\rd z
			\displaybreak[0]\\[1ex]
		&= \lim_{\delta\searrow 0}\frac{-1}{2\pi i}
			\int\limits_{\gamma_\delta} \bigl[R_zf,\psi_{t,w}(\zeta)\bigr]\,\rd z
			\displaybreak[0]\\[1ex]
		&\hspace*{-1.3ex}\stackrel{(\ref{A62})}{=}
			\lim_{\delta\searrow 0}\frac{-1}{2\pi i}
			\int\limits_{\gamma_\delta} (\Theta_t f)(z)\cdot(z-z_0)
			\bigl[R_zE(\Delta_0)\varepsilon_{z_0},E(\Delta_0)\psi_{t,w}(\zeta)\bigr]\rd z
			\\[1ex]
		&\hspace*{-1.3ex}\stackrel{(\ref{A14})}{=}
			\intop_{[a_-,a_+]}(\Theta_tf)(x)\cdot(x-z_0)\,\rd E_{\Delta_0;E(\Delta_0)\varepsilon_{z_0},E(\Delta_0)\psi_{t,w}(\zeta)}(x)
			\\[1ex]
		&=\intop_{[a_-,a_+]}(\Theta_tf)(x)\cdot(x-z_0)\,\rd E_{\Delta_0;\varepsilon_{z_0},\psi_{t,w}(\zeta)}(x).
	\end{align*}
	It remains to remember that
	$\frac{\rd E_{[a_-,a_+];u,v}}{\rd E_{\Delta_0;u,v}}=\mathds{1}_{[a_-,a_+]}$.
\end{proof}

\noindent
The next statement is the key lemma.

\begin{lemma}\thlab{A15}
	Assume that $H$ does not end with an indivisible interval towards $\infty$. Then, for each $s_0\in(0,\infty)$,
	we have
	\begin{equation}\label{A137}
		\cls\big\{\psi_{t,w}(z):\; t\in I_{\reg},\;t\geq s_0,\;z,w\in\bb C\setminus\bb R\big\}=
		\mc P(\mf h)\,[-]\,\mc E_A(\infty),
	\end{equation}
	where $\cls$ stands for `closed linear span' and $I_{\reg}$ is defined in \eqref{A89}.
\end{lemma}

\begin{proof}
	\hspace*{0pt} \\
	\textit{Step 1.}
	Fix a point $s\in I_{\reg}$ and let $M\subseteq\bb C$ be a set with a finite accumulation point.
	Let $\delta_j\in\mc P(\mf h_{\Lsh s})$, $j=0,\ldots,\Delta(H)-1$, be the elements
	$\delta_j \defeq (0,(\delta_{jk})_{k=0}^{\Delta(H)-1})$ where $\delta_{jk}$ is the
	Kronecker delta.  Then
	$(0;\delta_0),(\delta_0;\delta_1),\dots,(\delta_{\Delta(H)-2};\delta_{\Delta(H)-1})\in T(\mf h_{\Lsh s})$,
	and the boundary values of these elements vanish.
	Repeatedly applying the abstract Green identity \cite[(2.6) and Proposition~5.2]{kaltenbaeck.woracek:p4db}
	we obtain, for $k=0,\dots,\Delta(H)-1$ and $z\in C$,
	\[
		[\delta_k,\psi_s(z)] = [\delta_{k-1},z\psi_s(z)] = \ldots = [\delta_0,z^k\psi_s(z)]
		= [0,z^{k+1}\psi_s(z)] = 0.
	\]
	In particular, we have
	\begin{equation}\label{A16}
		\mc P(\mf h_{\Lsh s})\,[-]\,\cls\big\{\psi_s(z):\,z\in M\big\}
		\supseteq \spn\{\delta_0,\ldots,\delta_{\Delta(H)-1}\}.
	\end{equation}
	Applying the isomorphism $\Xi_s:\mc P(\mf h_{\Lsh s})\to\mf K(\omega_{\mf h}(s))$
	from \eqref{A144} and using \eqref{A47} we can deduce that
	\[
		\Xi_s\Bigl(\mc P(\mf h_{\Lsh s})\,[-]\,\cls\{\psi_s(z):\,z\in M\}\Bigr)
		= \Bigl(\cls\Bigl\{H_{\omega_{\mf h}(s)}(\qu z,\cdot)\binom 01:\,z\in M\Bigr\}\Bigr)^{[\perp]}.
	\]
	By analyticity, the space on the right-hand side equals $\ker\pi_-$, where $\pi_-$ denotes
	the projection onto the second component in $\mf K(\omega_{\mf h}(s))$.  As shown in the
	proof of \cite[Lemma~6.3 (subcase 3b)]{woracek:nass}, we have $\dim\ker\pi_-=\Delta(H)$.
	Thus equality must hold in \eqref{A16}.
	
	\bigskip

	\noindent\textit{Step 2.}
	First note that $P_s\delta_k=\delta_k$ where $P_s$ is as in \eqref{A128}, and hence
	\begin{multline*}
		[\psi_{s,w}(z),\delta_k]\!=\![P_s\psi_{s,w}(z),\delta_k]\stackrel{(\ref{A39})}{=}\lambda_{s,w}(z)[\psi_s(z),\delta_k]=0,
			\\[1ex]
		k=0,\ldots,\Delta(H)-1,\; z,w\in\rho(A).
	\end{multline*}
	Together with the fact that
	$\mc E_A(\infty)=\spn\{\delta_0,\dots,\delta_{\Delta(H)-1}\}$ (see
	\cite[Lemma~3.2\,(d)]{langer.woracek:gpinf}) this gives
	\[
		\mc E_A(\infty) \subseteq
		\{\psi_{t,w}(z):\; t\in I_{\reg},t\geq s_0,\,z,w\in\bb C\setminus\bb R\}^{[\perp]},
	\]
	i.e.\ the inclusion `$\subseteq$' in \eqref{A137}.
	To show the reverse inclusion, let $f\in\mc P(\mf h)$ be given and assume
	that $f\,[\perp]\,\psi_{s,w}(z)$, $s\in I_{\reg}$,
	$s\geq s_0$, $z,w\in\bb C\setminus\bb R$. Then, for each fixed $s\geq s_0$ and $z\in\bb C\setminus\bb R$,
	\[
		[P_sf,\psi_s(z)]=[f,\psi_{s,z}(z)]=0,\qquad z\in\bb C\setminus\bb R.
	\]
	By Step~1 we therefore have $P_sf\in\spn\{\delta_0,\ldots,\delta_{\Delta(H)-1}\}$.
	This tells us that $\psi(\mf h_{\Lsh s})P_sf=0$.  Since
	\[
		\psi(\mf h)f|_{(-1,0)\cup(0,s)}=\psi(\mf h_{\Lsh s})P_sf,
	\]
	and $s$ may be chosen arbitrarily large by our hypothesis that $H$ does not
	end indivisibly, it follows that $\psi(\mf h)f=0$.  Hence
	\[
		f\in\spn\{\delta_0,\ldots,\delta_{\Delta(H)-1}\}=\mc E_A(\infty),
	\]
	which proves the reverse inclusion.
\end{proof}

\noindent
If $s_1,s_2\in I_{\reg}$, $s_1<s_2$, then clearly $P_{s_1}\psi_{s_2}(z)=\psi_{s_1}(z)$,
and hence
\[
	\Theta_{s_1}f=\Theta_{s_2}f,\qquad f\in\mc P(\mf h_{\Lsh s_1}),
\]
where $P_s$ and $\Theta_s$ are defined in \eqref{A128} and \eqref{A88}.
Thus a map $\Theta_\infty$ on $\bigcup_{s\in I_{\reg}}\mc P(\mf h_{\Lsh s})$ is
well defined by setting
\begin{equation}\label{A138}
	\Theta_\infty f \defeq \Theta_{s_f}f \qquad \text{with $s_f$
	sufficiently large so that $f\in\mc P(\mf h_{\Lsh s_f})$}.
\end{equation}
Using Lemma~\ref{A15} we can extend Lemma~\ref{A21}.

\begin{proposition}\thlab{A20}
	Let $\Delta$ be a bounded Borel set and let $z_0\in\bb C\setminus\bb R$.  Then
	\begin{equation}\label{A60}
		\begin{aligned}{}
			[E(\Delta)f,g]=\int\limits_{-\infty}^\infty(\Theta_\infty f)(x)\cdot(x-z_0)
				\,\rd E_{\Delta;\varepsilon_{z_0},g}(x), &
				\\
			&\hspace*{-10ex} f\in\!\bigcup_{s\in I_{\reg}}\!\!\!\mc P(\mf h_{\Lsh s}),\;
			g\in\mc P(\mf h).
		\end{aligned}
	\end{equation}
\end{proposition}

\begin{proof}
	Fix $f\in\!\bigcup_{s\in I_{\reg}}\!\!\mc P(\mf h_{\Lsh s})$ and let $s_f$
	be so large that $f\in\mc P(\mf h_{\Lsh s_f})$.
	By Lemma~\ref{A21} the asserted relation holds for all
	$g\in\spn\{\psi_{s,w}(z):\,s\in I_{\reg},\,s\geq s_f,\,z,w\in\bb C\setminus\bb R\}$.

	Both sides of \eqref{A60} depend continuously on $g$. For the left-hand side
	this is obvious, for the right-hand side remember that $E_\Delta$ is
	compactly supported, and hence, for each continuous function $F$ on $\bb R$, the integral
	$\int_{-\infty}^\infty F\,\rd E_\Delta$ exists in the strong operator topology.
	We obtain from Lemma~\ref{A15} that \eqref{A60} holds for all $g\in\mc E_A(\infty)^{[\perp]}$.
	
	Finally, let an arbitrary element $g\in\mc P(\mf h)$ be given.
	Since $E(\Delta)g\in\mc E_A(\infty)^{[\perp]}$,
	we may apply what we have already shown and obtain
	\begin{align*}
		[E(\Delta)f,g] = [E(\Delta)f,E(\Delta)g]
			&= \intop_{-\infty}^\infty(\Theta_\infty f)(x)\cdot(x-z_0)
			\,\rd E_{\Delta;\varepsilon_{z_0},E(\Delta)g}(x)
			\\
		&\hspace*{-1.8ex}\stackrel{(\ref{A48})}
			= \intop_{-\infty}^\infty(\Theta_\infty f)(x)
			\cdot(x-z_0)\,\rd E_{\Delta;\varepsilon_{z_0},g}.
		\qedhere
	\end{align*}
\end{proof}

\begin{lemma}\thlab{A27}
	Let $\Delta$ be a bounded Borel set and $g\in\mc P(\mf h)$. Then
	\begin{equation}\label{A63}
		(x-z_0)\,\rd E_{\Delta;\varepsilon_{z_0},g}
		= (x-z_1)\,\rd E_{\Delta;\varepsilon_{z_1},g},\qquad
		z_0,z_1\in\bb C\setminus\bb R.
	\end{equation}
\end{lemma}

\begin{proof}
	To see this, note that
	\[
		I+(z_0-z_1)(A_\Delta-z_0)^{-1}=\intop_{\bb R}\frac{x-z_1}{x-z_0}\,\rd E_\Delta,\qquad
		z_0,z_1\in\bb C\setminus\bb R.
	\]
	From the identity \cite[(3.2)]{kaltenbaeck.woracek:krall} we obtain
	\[
		\big(I+(z_0-z_1)(A-z_0)^{-1}\big)E(\Delta)R_{z_1}^-=R_{z_0}^-.
	\]
	Hence, for each Borel set $\Delta'$ and $g\in\mc P(\mf h)$,
	\begin{align*}
		E_{\Delta;\varepsilon_{z_0},g}(\Delta')\hspace*{1.5ex} &\hspace*{-1.5ex}\stackrel{(\ref{A17})}{=}
			[E(\Delta\cap\Delta')R_{z_0}^-u,g]=[E(\Delta')E(\Delta)R_{z_0}^-u,E(\Delta)g]
			\\[1ex]
		&= \bigl[E(\Delta')E(\Delta)\bigl(I+(z_0-z_1)(A-z_0)^{-1}\bigr)R_{z_1}^-u,E(\Delta)g\bigr]
			\displaybreak[0]\\[1ex]
		&= \bigl[E(\Delta')\bigl(I+(z_0-z_1)(A_\Delta-z_0)^{-1}\bigr)E(\Delta)R_{z_1}^-u,E(\Delta)g\bigr]
			\displaybreak[0]\\[1ex]
		&\hspace*{-1.5ex}\stackrel{(\ref{A17})}{=}
			\biggl[\,\intop_{\Delta'}\frac{x-z_1}{x-z_0}\,\rd E_\Delta(x)\cdot E(\Delta)\varepsilon_{z_1},E(\Delta)g\biggr]
			\\[1ex]
		&= \intop_{\Delta'}\frac{x-z_1}{x-z_0}\,\rd E_{\Delta;E(\Delta)\varepsilon_{z_1},E(\Delta)g}(x)
			\stackrel{(\ref{A48})}{=} \intop_{\Delta'}\frac{x-z_1}{x-z_0}\,\rd E_{\Delta;\varepsilon_{z_1},g}(x).
	\end{align*}
	It follows that $E_{\Delta;\varepsilon_{z_0},g}\ll E_{\Delta;\varepsilon_{z_1},g}$ and
	$(x-z_0)\,\rd E_{\Delta;\varepsilon_{z_0},g} = (x-z_1)\,\rd E_{\Delta;\varepsilon_{z_1},g}$.
\end{proof}

\noindent
Now we are in position to relate the maps $\Theta_\infty$ and $\Theta_{\mf h}$.

\begin{lemma}\thlab{A19}
	Let $\upphi_{\mf h}$ be the distributional density in the representation \eqref{A50} of
	the Weyl coefficient $q_{\mf h}$ of $\mf h$, let $\psi(\upphi_{\mf h})$ be as
	in \ref{A139}, let $\Theta_{\mf h}$ be the isomorphism acting as in \eqref{A49}
	and let $\Theta_\infty$ be as in \eqref{A138}.
	Then
	\begin{align*}
		(\Theta_\infty f)(x) = \frac 1{x-i}\big(\big[\psi(\upphi_{\mf h})\circ\Theta_{\mf h}\big]f\big)(x)
		\qquad\mu_H\text{--a.e.}, \hspace*{20ex} & \\[1ex]
		f\in\!\bigcup_{s\in I_{\reg}}\mc P(\mf h_{\Lsh s})\cap\mc E_A(\infty)^{[\perp]}. &
	\end{align*}
	In particular, $\Theta_\infty$ maps $\bigcup_{s\in I_{\reg}}\mc P(\mf h_{\Lsh s})$ into $L^2(\mu_H)$.
\end{lemma}

\begin{proof}
	Let $f\in\bigcup_{s\in I_{\reg}}\mc P(\mf h_{\Lsh s})\cap\mc E_A(\infty)^{[\perp]}$
	be given.  For a bounded open interval $\Delta$ and
	$w\in\bb C\setminus\bb R$, we compute the inner product $[E(\Delta)f,\varepsilon_w]$ in two ways.
	
	On one hand, we have
	\begin{equation}\label{A68}
	\begin{aligned}{}
		[E(\Delta)f,\varepsilon_w]\hspace*{6.5pt}
		&\hspace*{-1ex}\stackrel{(\ref{A60})}=
			\int\limits_{\bb R}(\Theta_\infty f)(x)(x-i)\,\rd E_{\Delta;\varepsilon_i,\varepsilon_w}(x)
			\\
		&\hspace*{-3ex}\stackrel{(\ref{A64},\ref{A63})}=
			\int\limits_{\bb R}(\Theta_\infty f)(x)(x-i)\frac{x+i}{x-\qu w}\,\rd E_{\Delta;\varepsilon_i,\varepsilon_i}(x).
	\end{aligned}
	\end{equation}
	Using a standard argument we now relate $E_{\Delta;\varepsilon_i,\varepsilon_i}$ to $\mu_H$.
	Since $q_H$ is the $Q$-function induced by $A$ and the family
	$(\varepsilon_z)_{z\in\rho(A)}$, we have the representation
	\begin{align}
		q_{\mf h}(z) &= \qu{q_{\mf h}(i)}+(z+i)[\varepsilon_i,\varepsilon_i]+(z^2+1)[R_z\varepsilon_i,\varepsilon_i]
			\nonumber
			\\[1ex]
		&= \qu{q_{\mf h}(i)}+(z+i)[\varepsilon_i,\varepsilon_i]
			+(z^2+1)\bigl[R_zE(\Delta^c)\varepsilon_i,E(\Delta^c)\varepsilon_i\bigr]
			\nonumber
			\\[1ex]
		&\quad +(z^2+1)\bigl[(A_\Delta-z)^{-1}E(\Delta)\varepsilon_i,E(\Delta)\varepsilon_i\bigr].
		\label{A67}
	\end{align}
	Let $[a_-,a_+]\subseteq\Delta$ be such that $\mu_H(\{a_-\})=\mu_H(\{a_+\})=0$
	and $E(\{a_-\})=E(\{a_+\})=0$.  Observing that all summands on the
	right-hand side of \eqref{A67} apart from the last one are analytic across $\Delta$
	we compute (where $\gamma_\delta$ is as in the proof of Lemma~\ref{A21})
	\begin{align*}
		\mu_H\big([a_-,a_+]) &
			\stackrel{(\ref{A95})}= \frac 1\pi\lim_{\delta\searrow0}\int\limits_{[a_-,a_+]}
			\Im q_H(x+i\delta)\,\rd x
			\\
		& \hspace*{-2.5ex}\stackrel{{\scriptscriptstyle q_{\mf h}(\qu z)=\qu{q_{\mf h}(z)}}}=
			\frac 1\pi\lim_{\delta\searrow0}\frac{-1}{2i}\int_{\gamma_\delta}q_H(z)\,\rd z
			\displaybreak[0]\\
		& \hspace*{-0.4ex}\stackrel{(\ref{A67})}= \lim_{\delta\searrow0}\frac{-1}{2\pi i}\int_{\gamma_\delta}(z^2+1)
			\bigl[(A_\Delta-z)^{-1}E(\Delta)\varepsilon_i,E(\Delta)\varepsilon_i\bigr]\,\rd z
			\\
		& \hspace*{-0.4ex}\stackrel{(\ref{A14})}=
			\hspace*{-5pt}\int\limits_{[a_-,a_+]}\hspace*{-5pt}(x^2+1)\,\rd E_{\Delta;E(\Delta)\varepsilon_i,E(\Delta)\varepsilon_i}
			\stackrel{(\ref{A48})}= \hspace*{-1ex}
			\int\limits_{[a_-,a_+]}\hspace*{-1ex}(x^2+1)\,\rd E_{\Delta;\varepsilon_i,\varepsilon_i}.
	\end{align*}
	From this it follows that $\mathds{1}_\Delta(x)\,\rd\mu_H(x)\ll \rd E_{\Delta;\varepsilon_i,\varepsilon_i}$ and
	\[
		\frac{\mathds{1}_\Delta(x)\,\rd\mu_H(x)}{\rd E_{\Delta;\varepsilon_i,\varepsilon_i}}=x^2+1.
	\]
	Hence, we can further rewrite the last integral in \eqref{A68}, and obtain
	\begin{equation}\label{A147}
		[E(\Delta)f,\varepsilon_w]=
		\int\limits_{\bb R}(\Theta_\infty f)(x)(x-i)\frac{x+i}{x-\qu w}\,\frac{\mathds{1}_\Delta(x)\rd\mu_H(x)}{1+x^2}\,.
	\end{equation}

	On the other hand, let $\mc E_{\upphi_{\mf h}}(\infty)^{[\perp]}$ be the algebraic
	eigenspace at infinity of $A_{\upphi_{\mf h}}$, let $E_{\upphi_{\mf h}}$ be the spectral
	measure of $A_{\upphi_{\mf h}}$ and let $\hat\varepsilon_w$ be as in \eqref{A145}.
	Since $\Theta_{\mf h}$ is isometric,
	$\ran E_{\upphi_{\mf h}}(\Delta)\subseteq \mc E_{\upphi_{\mf h}}(\infty)^{[\perp]}$
	and $\psi(\upphi_{\mf h})$ is isometric on $\mc E_{\upphi_{\mf h}}(\infty)^{[\perp]}$
	(by \cite[Theorem~5.3]{langer.woracek:ninfrep}), we have
	\begin{align}
		\hspace*{-1ex} [E(\Delta)f,\varepsilon_w]
		&= [E(\Delta)f,\,E(\Delta)\varepsilon_w]
			\nonumber
			\\[1ex]
		&= \bigl[\Theta_{\mf h}E(\Delta)f,\,\Theta_{\mf h}E(\Delta)\varepsilon_w\bigr]_{\Pi(\upphi_{\mf h})}
			\nonumber
			\displaybreak[0]\\[1ex]
		&\hspace*{-1ex}\stackrel{(\ref{A122})}=
			\bigl[E_{\upphi_{\mf h}}(\Delta)\Theta_{\mf h}f,\,E_{\upphi_{\mf h}}(\Delta)
			\Theta_{\mf h}\varepsilon_w\bigr]_{\Pi(\upphi_{\mf h})}
			\nonumber
			\\[1ex]
		&\hspace*{-1ex}\stackrel{(\ref{A49})}=
			\int\limits_{\bb R}\Bigl(\psi(\upphi_{\mf h})E_{\upphi_{\mf h}}(\Delta)\Theta_{\mf h}f\Bigr)(x)
			\qu{\Bigl(\psi(\upphi_{\mf h})E_{\upphi_{\mf h}}(\Delta)\hat\eps_w\Bigr)(x)}
			\,\frac{\rd\mu_H(x)}{1+x^2}.
			\label{A66}
	\end{align}
	The mappings $\psi(\upphi_{\mf h})\circ E_{\upphi_{\mf h}}(\Delta)|_{\mc E_{\upphi_{\mf h}}(\infty)^{[\perp]}}$ and
	$(\,\cdot\,\mathds{1}_\Delta)\circ\psi(\upphi_{\mf h})$ are continuous
	on $\mc E_{\upphi_{\mf h}}(\infty)^{[\perp]}$.
	By the definition of $\psi(\upphi_{\mf h})$
	and \cite[Proposition~3.1]{kaltenbaeck.woracek:p2db}, they coincide
	on all compactly supported functions of $\mc B_2(\upphi_{\mf h})$
	(for this notation see \cite[\S5]{langer.woracek:ninfrep}).
	Their continuity implies that they coincide on $\mc E_{\upphi_{\mf h}}(\infty)^{[\perp]}$.
	The assumption $f\in\mc E_A(\infty)^{[\perp]}$ and the relation \eqref{A55}
	imply that $\Theta_{\mf h}f\in\mc E_{\upphi_{\mf h}}(\infty)^{[\perp]}$.
	Hence
	\[
		\Bigl(\psi(\upphi_{\mf h})E_{\upphi_{\mf h}}(\Delta)\Theta_{\mf h}f\Bigr)(x)
		= \mathds{1}_\Delta(x)\Bigl(\psi(\upphi_{\mf h})\Theta_{\mf h}f\Bigr)(x)
		\qquad x\in\bb R\quad \text{$\mu_H$-a.e}.
	\]
	Using \cite[Proposition~3.1]{kaltenbaeck.woracek:p2db} and the fact that
	$\psi(\upphi_{\mf h})$ acts as the identity on compactly supported functions we obtain
	\[
		\Bigl(\psi(\upphi_{\mf h})E_{\upphi_{\mf h}}(\Delta)\hat\eps_w\Bigr)(x)
		= \Bigl(\psi(\upphi_{\mf h})\bigl(\mathds{1}_\Delta\hat\eps_w\bigr)\Bigr)(x)
		= \mathds{1}_\Delta(x)\cdot\frac{x-i}{x-w}\,.
	\]
	The integral on the right-hand side of \eqref{A66} thus equals
	\begin{equation}\label{A148}
		\int\limits_{\bb R} \Big(\psi(\upphi_{\mf h})\Theta_{\mf h}f\Bigr)(x)
		\cdot\frac{x+i}{x-\qu w} \mathds{1}_\Delta(x)
		\,\frac{\rd\mu_H(x)}{1+x^2}\,.
	\end{equation}
	Since the linear span of the functions $x\mapsto\frac{x+i}{x-\qu w}$, $w\in\rho(A)$,
	is dense in $L^2\big(\frac{\mathds{1}_\Delta(x)}{1+x^2}\mu_H(x)\big)$,
	we conclude from \eqref{A147}, \eqref{A66} and \eqref{A148} that
	\[
		(\Theta_\infty f)(x)(x-i) = \Bigl(\psi(\upphi_{\mf h})\Theta_{\mf h}f\Bigr)(x),\qquad
		x\in\Delta\quad \text{$\mu_H$-a.e}.
	\]
	Since $\Delta$ was an arbitrary bounded open interval, the assertion follows.
\end{proof}

\medskip

\noindent
To finish the proof of \eqref{A121}, let $[a,b]\subseteq(0,\infty)$ and
denote by $\Omega$ the map defined on $L^2(H|_{[a,b]})$ as
\[
	(\Omega f)(z)\defeq \begin{pmatrix}
	\int_{[a,b]}\big[\psi(\mf h_{\Lsh b})\phi_b(\qu z)\big](t)^*H(t)f(t)\,\rd t \\[2ex]
	\int_{[a,b]}\big[\psi(\mf h_{\Lsh b})\psi_b(\qu z)\big](t)^*H(t)f(t)\,\rd t
	\end{pmatrix},
\]
where $\phi_b$ and $\psi_b$ are as in \eqref{A146} and \eqref{A40}, respectively.
The function $[\psi(\mf h_{\Lsh b})\psi_b(\qu z)](t)$ is a solution of \eqref{A30}
with $z$ replaced by $\qu z$, which assumes boundary values at $0$, namely $(0,\ 1)^T$.
The function $\upphi(t;\qu z)$ shares these properties,
and hence we have $[\psi(\mf h_{\Lsh b})\psi_b(\qu z)](t)=\upphi(t;\qu z)$.
Thus the second component of $\Omega$ can be rewritten as
\[
	\int_{[a,b]}\upphi(t,\qu z)^*H(t)f(t)\,\rd t.
\]
The asserted formula \eqref{A121} for the action of $\Theta_H$ is now obtained by putting together
the so far collected knowledge.  Consider the following diagram
(here $\pi_-$ denotes the projection onto the second component,
and references between $\#$ are for the proof of the commutativity of the corresponding part
of the diagram):
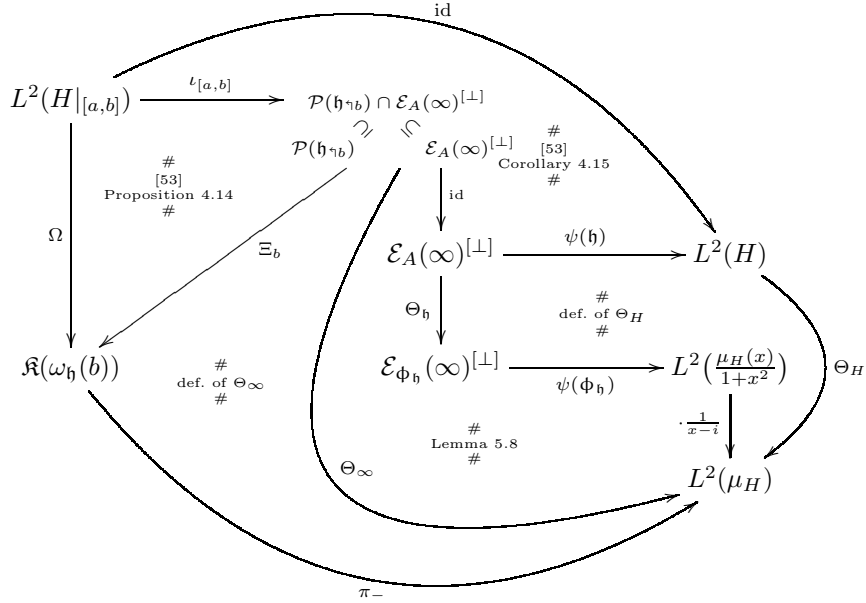
\[
\label{A57}
	\xymatrix{
		L^2(H|_{[a,b]}) \ar[rr]^{\hspace*{-6mm}\iota_{[a,b]}\hspace*{6mm}} \ar[dd]_\Omega
			\ar@/^5pc/[rrrd]^{\id} & &
			\parbox{40mm}{\hspace*{-2mm}\raisebox{-12mm}{{\scriptsize
	\begin{tabular}{cc@{\hspace*{-30pt}}c@{\hspace*{-35pt}}cc}
		&& $\mc P(\mf h_{\Lsh b})\cap\mc E_A(\infty)^{[\perp]}$ && \\
		& \hspace*{-10pt}\raisebox{-3pt}{\begin{rotate}{45}$\supseteq$\end{rotate}}\hspace*{10pt}
			&& \hspace*{3pt}\raisebox{1.5pt}{\begin{rotate}{-45}$\subseteq$\end{rotate}}\hspace*{-3pt} & \\
		$\mc P(\mf h_{\Lsh b})$ && && $\mc E_A(\infty)^{[\perp]}$
	\end{tabular}
	}}\hspace*{2mm}}
		\ar[ddll]^{\Xi_b} \ar@/_10pc/[dddr]^{\Theta_\infty} \ar[d]^{\raisebox{-16pt}{\tiny$\id$}} &
	\\
		& & \mc E_A(\infty)^{[\perp]} \ar[r]^{\psi(\mf h)} \ar[d]_{\Theta_{\mf h}}
			\save[]+<15mm,12mm>*\txt<8pc>{{\tiny Corollary 4.15}} \restore
			\save[]+<15mm,14mm>*\txt<8pc>{{\tiny\cite{kaltenbaeck.woracek:p4db}}} \restore
			\save[]+<15mm,16.5mm>*\txt<8pc>{{\tiny\#}} \restore
			\save[]+<15mm,10mm>*\txt<8pc>{{\tiny\#}} \restore
		& L^2(H) \ar@/^3pc/[dd]^{\Theta_H}
	\\
			\mf K(\omega_{\mf h}(b)) \ar@/_5pc/[drrr]_{\pi_-}
			\save[]+<13mm,23mm>*\txt<8pc>{{\tiny Proposition 4.14}} \restore
			\save[]+<13mm,25mm>*\txt<8pc>{{\tiny\cite{kaltenbaeck.woracek:p4db}}} \restore
			\save[]+<13mm,27.5mm>*\txt<8pc>{{\tiny\#}} \restore
			\save[]+<13mm,21mm>*\txt<8pc>{{\tiny\#}} \restore
			\save[]+<20mm,0.5mm>*\txt<8pc>{{\tiny\#}} \restore
			\save[]+<20mm,-4mm>*\txt<8pc>{{\tiny\#}} \restore
			\save[]+<20mm,-2mm>*\txt<8pc>{{\tiny def.\ of $\Theta_\infty$}} \restore
			&
			& \mc E_{\upphi_{\mf h}}(\infty)^{[\perp]} \ar[r]_{\psi(\upphi_{\mf h})} & L^2\big(\frac{\mu_H(x)}{1+x^2}\big)
				\ar[d]_{\cdot\frac 1{x-i}}
			\save[]+<-17mm,9.5mm>*\txt<8pc>{{\tiny\#}} \restore
			\save[]+<-17mm,5mm>*\txt<8pc>{{\tiny\#}} \restore
			\save[]+<-17mm,7mm>*\txt<8pc>{{\tiny def.\ of $\Theta_H$}} \restore
			\save[]+<-34mm,-8mm>*\txt<8pc>{{\tiny\#}} \restore
			\save[]+<-34mm,-12mm>*\txt<8pc>{{\tiny\#}} \restore
			\save[]+<-34mm,-10mm>*\txt<8pc>{{\tiny Lemma~\ref{A19}}} \restore
	\\
		& & & L^2(\mu_H)
	}
\]
We see that \eqref{A121} holds for each $f\in L^2(H)$ supported on $[a,b]$.
For each fixed $T>0$ both sides of \eqref{A121} depend continuously on $f$ when $f$
varies in the set of all elements of $L^2(H)$ whose support is bounded above by $T$. Remember here
that $\upphi(\,\cdot\,;z)\in L^2(0,T)$. Hence, \eqref{A121} holds for all $f\in L^2(H)$ whose support is
bounded above by $T$.
\renewcommand{\qedsymbol}{{\large\smiley}}
\popQED{\qed}
\renewcommand{\qedsymbol}{\raisebox{-2pt}{\large\ding{113}}}

\subsection{Computation of $\bm{\Theta_H^{-1}}$ as an integral transform}
\label{A83}

The final task in the proof of Theorem~\ref{A4} is to establish the formula \eqref{A126}
for the action of $\Theta_H^{-1}$.

Denote by $E$ the spectral family associated with $T(H)$.
Let $f\in L^2(H)$ with $\sup(\supp f)<b$, let $h\in L^2(H)$ be bounded and with
$\sup(\supp h)<b$, and let $\Delta\subseteq\bb R$ be a finite interval.
Using (i), (ii) of Theorem~\ref{A4} and Fubini's theorem we obtain
\begin{align*}
	\bigl(E(\Delta)f,h\bigr)_{L^2(H)}
	&= \int_{\bb R} \mathds 1_\Delta(t)\bigl(\Theta_H f\bigr)(t)\qu{\bigl(\Theta_H h\bigr)(t)}\,\rd\mu_H(t)
		\\[1ex]
	&= \int_\Delta \bigl(\Theta_H f\bigr)(t) \int_a^b h(x)^*H(x)\upvarphi(x;t)\,\rd x\,\rd\mu_H(t)
		\\[1ex]
	&= \int_a^b h(x)^*H(x) \int_\Delta \bigl(\Theta_H f\bigr)(t)\upvarphi(x;t)\,\rd\mu_H(t)\,\rd x.
\end{align*}
Since the set $\{h\in L^2(H):\sup(\supp h)<b,\,h\text{ bounded}\}$ is dense in $L^2(H)$,
it follows that
\begin{equation}\label{A8}
	\bigl(E(\Delta)f\bigr)(x) = \int_\Delta \bigl(\Theta_H f\bigr)(t)\upvarphi(x;t)\,\rd\mu_H(t)
	\qquad H\text{-a.e.}
\end{equation}
Both sides of this equality depend continuously on $f$ and therefore this relation is
valid for arbitrary $f\in L^2(H)$.

To complete the proof, let $g\in L^2(\mu_H)$ with compact support be given.  Choose a finite
interval $\Delta$ which contains $\supp g$. Since $\Theta_H$ intertwines $T(H)$
with the multiplication operator in $L^2(\mu_H)$, we have
\begin{equation}\label{A9}
	E(\Delta)\circ\Theta_H^{-1}=\Theta_H^{-1}\circ \big(\cdot\mathds{1}_\Delta\big).
\end{equation}
Thus
\begin{align*}
	\big(\Theta_H^{-1}g\big)(x) &= \big(\Theta_H^{-1}(\mathds{1}_\Delta g)\big)(x)
	\stackrel{(\ref{A9})}= \big(E(\Delta)(\Theta_H^{-1}g)\big)(x)
	\\[1ex]
	&\hspace*{-1.5ex}\stackrel{(\ref{A8})}= \int_\Delta g(t)\upvarphi(x;t)\,\rd\mu_H(t)
	= \int_{\bb R} g(t)\upvarphi(x;t)\,\rd\mu_H(t)
	\qquad H\text{-a.e.}
\end{align*}
This finishes the proof of Theorem~\ref{A4}.

\subsection{The connection between the point mass at 0 and the behaviour of $\bm H$}
\label{A84}

Before we prove Proposition~\ref{A77}, we need a lemma.

\begin{lemma}\thlab{A80}
	Let $H\in\bb H$ with $\dom H=(0,\infty)$ and let $\mf h$ be
	the general Hamiltonian as in \ref{A109}.  Moreover, let $\omega_{\mf h}$ be
	the chain of matrices as in \ref{A52}.  Then, for each $x\in(0,\infty)$,
	\begin{equation}\label{A85}
		\biggl[\frac{\partial}{\partial z}\omega_{\mf h,21}(x;z)\biggr]\Bigg|_{z=0}
		= - \int_0^x h_{22}(t)\rd t.
	\end{equation}
\end{lemma}

\begin{proof}
	Let $x_1\in(0,x)$.  Integrating \eqref{A70} we obtain
	\[
		\omega_{\mf h}(x;z) - \omega_{\mf h}(x_1;z)
		= z\int_{x_1}^x \omega_{\mf h}(t;z)H(t)J\rd t
	\]
	for $z\in\bb C$.  If we differentiate both sides with respect to $z$,
	set $z=0$ and use \eqref{A74}, we arrive at
	\begin{equation}\label{A86}
	\begin{aligned}
		\biggl[\frac{\partial}{\partial z}\omega_{\mf h}(x;z)\biggr]\Bigg|_{z=0}
		- \biggl[\frac{\partial}{\partial z}\omega_{\mf h}(x_1;z)\biggr]\Bigg|_{z=0}
		&= \int_{x_1}^x \omega_{\mf h}(t;0)H(t)J\rd t
		\\[1ex]
		&= \int_{x_1}^x H(t)J\rd t.
	\end{aligned}
	\end{equation}
	It follows from \cite[Theorem~4.1]{langer.woracek:gpinf} that
	\[
		\lim_{x_1\searrow0}\biggl[\frac{\partial}{\partial z}\omega_{\mf h,21}(x_1;z)\biggr]\Bigg|_{z=0}
		= 0,
	\]
	which, together with \eqref{A86}, implies \eqref{A85}.
\end{proof}

\begin{proof}[Proof of Proposition~\ref{A77}]
	Assume first that $\mu_H(\{0\})>0$.  Define $g\in L^2(\mu_H)$ by
	\[
		g(t) = \begin{cases} 1, & t=0, \\[0.5ex] 0, & t\ne0. \end{cases}
	\]
	Then
	\[
		(\Theta_H^{-1}g)(x)
		= \int_{-\infty}^\infty g(t)\upvarphi(x;t)\rd\mu_H(t)
		= \mu_H\bigl(\{0\}\bigr)\upvarphi(x;0) = \mu_H\bigl(\{0\}\bigr)\binom01.
	\]
	Since $\Theta_H$ is an isometric isomorphism, we obtain
	\begin{align*}
		\mu_H\bigl(\{0\}\bigr) &= \|g\|_{L^2(\mu_H)}^2
		= \bigl\|\Theta_H^{-1}g\bigr\|_{L^2(H)}^2
		= \bigl[\mu_H\bigl(\{0\}\bigr)\bigr]^2\int_0^\infty \binom01^* H(x)\binom01\rd x
		\\[1ex]
		&= \bigl[\mu_H\bigl(\{0\}\bigr)\bigr]^2\int_0^\infty h_{22}(x)\rd x,
	\end{align*}
	which implies \eqref{A78} and \eqref{A79}.

	Now assume that \eqref{A78} is satisfied.
	Let $c>0$ be large enough such that $(0,c)$ is not an indivisible interval
	and introduce the Hamiltonian function
	\[
		H_c(x) \defeq \begin{cases}
			H(x), & x\in(0,c], \\[1ex]
			\smmatrix{1}{0}{0}{0}, & x\in(c,\infty),
		\end{cases}
	\]
	which belongs to $\bb H$, satisfies $\Delta(H_c)=\Delta(H)$ and
	is in the limit point case at infinity.  Moreover, let $\mf h_c$
	be the corresponding singular general Hamiltonian as in \ref{A109}
	and let $\omega_{\mf h_c}$ be the chain of matrices as in \ref{A52}.
	The singular Weyl coefficient $q_{H_c}$ is given by
	\[
		q_{H_c}(z) = q_{\mf h_c}(z) = \omega_{\mf h_c}(c;z)\star\infty
		= \frac{\omega_{\mf h_c,11}(c;z)}{\omega_{\mf h_c,21}(c;z)}\,,
		\qquad z\in\bb C\setminus\bb R;
	\]
	it belongs to $\mc N_{\Delta(H)}^{(\infty)}$ and is meromorphic in $\bb C$.
	Let $\mu_{H_c}$ be the spectral measure associated with $H_c$ via \eqref{A95}.
	Then Lemma~\ref{A80} and \eqref{A74} imply that, for every $\eps>0$,
	\begin{align*}
		\mu_{H_c}\bigl((-\eps,\eps)\bigr)
		&\ge \mu_{H_c}\bigl(\{0\}\bigr) = \Res(q_{H_c};0)
		= -\frac{\omega_{\mf h_c,11}(c;0)}{\bigl[\frac{\partial}{\partial z}\omega_{\mf h_c,21}(c;z)\bigr]\big|_{z=0}}
		\\[1ex]
		&= \biggl[\int_0^c h_{22}(x)\rd x\biggr]^{-1}
		\ge \biggl[\int_0^\infty h_{22}(x)\rd x\biggr]^{-1}
		\eqdef M.
	\end{align*}
	By Theorem~\ref{A94}\,(ii) we have $q_H(z) = \lim_{c\to\infty}q_{H_c}(z)$
	locally uniformly in $\bb C\setminus\bb R$.
	Since $\ind_- q_{H_c} = \Delta(H_c) = \Delta(H) = \ind_- q_H$,
	we can apply Lemma~\ref{A87},
	which implies that, for every $\eps>0$ such that $\mu_H(\{-\eps\})=\mu_H(\{\eps\})=0$,
	\[
		\mu_H\bigl((-\eps,\eps)\bigr)
		= \lim_{c\to\infty}\mu_{H_c}\bigl((-\eps,\eps)\bigr) \ge M.
	\]
	Hence $\mu_H(\{0\})\ge M > 0$.
\end{proof}

%
%
\section{Inverse theorems}
\label{sec-inverse}
%
%

By the procedure described in Theorem~\ref{A107} a map from $\bb H$ to $\bb M$,
namely $H\mapsto\mu_H$, is well defined.  Hence it is meaningful to pose inverse problems.
Concisely formulated, we face the task to determine the range and kernel of the
mapping $H\mapsto\mu_H$.

In this section we complete this task. In fact, we provide somewhat more detailed results.
They include singular Weyl functions and a local version of the uniqueness theorem.
Proofs are again relatively simple; they are carried out in the same manner as in \S4,
using the basic identifications \ref{A109}, \ref{A108}, and
some results taken from the literature.
Recall the notation from Definition~\ref{A99}.
Theorems~\ref{A5} and \ref{A6} are the analogues for our class $\bb H$ of Hamiltonians
with two singular endpoints to de Branges' celebrated inverse spectral theorem.

\begin{theorem}[\textbf{Existence Theorem}]\thlab{A5} \rule{0ex}{1ex}\\[1ex]
	The following statements hold.
	\begin{enumerate}[{\rm(i)}]
	\item Let $q\in\mc N_{<\infty}^{(\infty)}$ with $\ind_-q>0$.
		Then there exists a Hamiltonian	$H\in\bb H$ with $q\in[q]_H$.
	\item Let $\mu\in\bb M$ with $\Delta(\mu)>0$. Then there exists a Hamiltonian
		$H\in\bb H$ with $\mu_H=\mu$.
	\end{enumerate}
\end{theorem}

\begin{proof}
	To show (i) let $q\in\mc N_{<\infty}^{(\infty)}$ be given.
	According to \cite[Theorem~1.4]{kaltenbaeck.woracek:p6db} there exists a
	general Hamiltonian $\mf h_0$ with $q_{\mf h_0}=q$.  Moreover,
	\cite[Theorem~3.1]{langer.woracek:gpinf} implies that $\mf h_0\in \mf H_0$.
	Applying an appropriate reparameterization we may assume that
	\begin{itemize}
	\item[--] $\mf h_0$ is defined on $(-1,0)\cup(0,\infty)$;
	\item[--] the Hamiltonian function of $\mf h$ equals
		$x^{-2}\smmatrix{1}{0}{0}{0}$
		for $x\in(-1,0)$;
	\item[--] $b_{\oe+1}=0$;
	\item[--] $E=\{-1,1,\infty\}$.
	\end{itemize}
	Let $H$ be the Hamiltonian function of $\mf h_0$ on the interval $(0,\infty)$.
	Then, by our basic identification \ref{A109}, we have $H\in\bb H$.
	
	Let $\mf h$ be the general Hamiltonian built from $H$ in the basic identification \ref{A108}
	with $\oe_1$, $b_{1,j}$, $d_{1,j}$ all equal to $0$.
	Then $\mf h$ and $\mf h_0$ differ only in the data part $\oe_1,b_{1,j},d_{1,j}$.
	From \cite[Corollary~5.9]{langer.woracek:gpinf} we obtain that
	\[
		q_H(z)=q_{\mf h}=q_{\mf h_0}-\sum_{l=1}^{2\Delta(H)}z^l d_{1,l-1}
		+\sum_{l=1}^{\oe_1} z^{2\Delta(H)+l}b_{1,\oe_1+1-l},
	\]
	i.e.\ $q_H$ and $q$ differ only by a polynomial with real coefficients and vanishing constant
	term.  Hence $q\in[q]_H$.

	For the proof of (ii) let $\mu\in\bb M$ be given.
	Choose $q\in\mc N_{<\infty}^{(\infty)}$ with $\mu_q=\mu$;
	this is possible by \cite[Theorem~3.9\,(v)]{langer.woracek:ninfrep}\footnote{The
	set on the right-hand side of \cite[Theorem~3.9\,(v)]{langer.woracek:ninfrep} is
	certainly non-empty.}.
	An application of the already proved item (i)
	provides us with a Hamiltonian $H\in\bb H$ such that $q_H-q$ is a real polynomial.
	By the Stieltjes inversion formula \eqref{A142} and the
	definition of $\mu_H$, it follows that $\mu_H=\mu_{q_H}=\mu$.
\end{proof}

\noindent
As we have seen in Proposition~\ref{A97}, Hamiltonians which are --- essentially --- reparameterizations
of each other have --- essentially --- the same singular Weyl coefficients and have the
same spectral measures.
The converse of this fact is an important result.

\begin{theorem}[\textbf{Global Uniqueness Theorem}]\thlab{A6} \rule{0ex}{1ex}\\[1ex]
	Let $H_1,H_2\in\bb H$ be given.
	\begin{enumerate}[{\rm(i)}]
	\item Assume that there exist singular Weyl coefficients $q_{H_1}$ and $q_{H_2}$ of $H_1$ and
		$H_2$, respectively, such that $q_{H_1}-q_{H_2}$ is a real polynomial, and set
		$\alpha \defeq (q_{H_1}-q_{H_2})(0)$. Then the Hamiltonians
		\begin{equation}\label{A111}
			H_1\quad\text{and}\quad
			\begin{pmatrix} 1 & \alpha \\ 0 & 1 \end{pmatrix}
			H_2
			\begin{pmatrix} 1 & 0 \\ \alpha & 1 \end{pmatrix},
		\end{equation}
		are reparameterizations of each other.

		In particular, if\, $[q]_{H_1}=[q]_{H_2}$, then $H_1$ and $H_2$ are
		reparameterizations of each other.
	\item If\, $\mu_{H_1}=\mu_{H_2}$, then there exists a real constant $\alpha$
	such that the Hamiltonians in \eqref{A111} are reparameterizations of each other.
	\end{enumerate}
\end{theorem}

\begin{proof}
	Let $H_1,H_2\in\bb H$ be given. Assume that both Hamiltonians are defined on $(0,\infty)$. This is no
	loss in generality since it can always be achieved by a reparameterization, and reparameterizations change
	neither singular Weyl coefficients nor spectral measures; see Proposition~\ref{A97}.

	First we consider the case when $[q]_{H_1}=[q]_{H_2}$.
	Let $\mf h_1$ and $\mf h_2$ be the general Hamiltonians defined for $H_1$ and $H_2$ by the
	basic identification \ref{A108} with $\oe_1$, $b_{1,j}$, $d_{1,j}$ all equal to $0$
	with some base points $x_1$ and $x_2$.  Then
	\[
		q_{\mf h_1}\in[q]_{H_1}\quad\text{and}\quad q_{\mf h_2}\in[q]_{H_2}.
	\]
	Applying a reparameterization to $\mf h_2$, we can achieve that $x_1=x_2$.
	By \cite[Corollary~5.9]{langer.woracek:gpinf} there exist numbers
	$d_{1,0},\ldots,d_{1,2\Delta(H_2)-1}\in\bb R$ and $\oe_1\in\bb N_0$,
	$b_{1,1},\ldots,b_{1,\oe_1}\in\bb R$ such that the Weyl coefficient of
	the general Hamiltonian $\tilde{\mf h}_2$ constituted by the same data as $\mf h_2$ with exception of
	$d_{1,j},\oe_1,b_{1,j}$ is equal to $q_{\mf h_1}$. By the uniqueness part in
	\cite[Theorem~1.4]{kaltenbaeck.woracek:p6db}, thus, $\mf h_1$ and $\tilde{\mf h}_2$ are
	reparameterizations of each other. In particular, their Hamiltonian functions on $(0,\infty)$ are
	reparameterizations of each other. However, the Hamiltonian function of $\mf h_1$ on this
	interval is $H_1$ and the one of $\tilde{\mf h}_2$ is $H_2$.

	Now assume that some singular Weyl coefficients $q_{H_1}$ and
	$q_{H_2}$ differ by a real constant, say $\alpha$. Consider the Hamiltonian
	$H_0 \defeq \smmatrix{1}{\alpha}{0}{1}H_2\smmatrix{1}{0}{\alpha}{1}$.
	We know from Proposition~\ref{A97} (and its proof) that $H_0\in\bb H$ and that $q_{H_0}-q_{H_2}=\alpha$ when we choose the same base point
	in the construction of $q_{H_1}$ and $q_{H_2}$, respectively. We thus have $q_{H_1}=q_{H_0}$.
	By what we proved in the previous paragraph, this implies that $H_1$ and
	$H_0$ are reparameterizations of each other; hence (i) is shown.
	
	For the proof of (ii) assume that $\mu_{H_1}=\mu_{H_2}$.  Choose $x_0\in(0,\infty)$,
	and let $q_{H_1}$ and $q_{H_2}$ be singular Weyl coefficients of
	$H_1$ and $H_2$, respectively, built with the base point $x_0$.
	By \cite[Theorem~3.9\,(iv)]{langer.woracek:ninfrep} the difference
	$q_{H_1}-q_{H_2}$ is a real polynomial.  Now the already proved item (i) yields \eqref{A111}.
\end{proof}

\noindent
Our last result in this section is a refined version of the above uniqueness theorem.
It asserts that certain beginning sections of two Hamiltonians $H_1,H_2\in\bb H$ coincide
if (and only if) some of their singular Weyl coefficients are exponentially close.
Local uniqueness theorems for one-dimensional Schr\"odinger operators were first
proved by B.~Simon in \cite[Theorem~1.2]{simon:1999}.
For canonical systems with a regular left endpoint a local uniqueness theorem
was proved in \cite[Theorem~1.2]{langer.woracek:lokinv};
see also \cite[Section~4]{langer:2016} for a formulation in terms of transfer functions.

\begin{theorem}[\textbf{Local Uniqueness Theorem}]\thlab{A31} \rule{0ex}{1ex} \\
	Let $H_1,H_2\in\bb H$ with $\dom(H_i)=(a_i,b_i)$, $i=1,2$, be given and set
	\begin{equation}\label{A58}
		s_i(\tau)\defeq\sup\biggl\{x\in(a_i,b_i):\,
		\intop_{a_i}^x\sqrt{\det H_i(\xi)}\;\rd\xi<\tau\biggr\},\qquad \tau>0, \;i=1,2.
	\end{equation}
	Then, for each $\tau>0$, the following statements are equivalent.
	\begin{enumerate}[{\rm(i)}]
	\item The Hamiltonian $H_1|_{(a_1,s_1(\tau))}$ is a reparameterization of $H_2|_{(a_2,s_2(\tau))}$.
	\item There exist singular Weyl coefficients $q_{H_1}$ and $q_{H_2}$ of $H_1$ and $H_2$, respectively,
		and there exists a $\beta\in(0,\pi)$ such that for each $\varepsilon>0$,
		\[
			q_{H_1}(re^{i\beta})-q_{H_2}(re^{i\beta})=
			{\rm O}\bigl(e^{(-2\tau+\varepsilon)r\sin\beta}\bigr),
			\quad r\to\infty.
		\]
	\item There exist singular Weyl coefficients $q_{H_1}$ and $q_{H_2}$ of $H_1$ and $H_2$, respectively,
		and there exists a $k\geq 0$ such that for each $\delta\in(0,\frac\pi 2)$,
		\begin{equation}\label{A104}
			q_{H_1}(z)-q_{H_2}(z)={\rm O}\bigl((\Im z)^ke^{-2\tau\Im z}\bigr),
			\quad |z|\to\infty,\; z\in\Gamma_\delta,
		\end{equation}
		where $\Gamma_\delta$ is the Stolz angle
		$\Gamma_\delta\defeq\{z\in\bb C:\,\delta\leq\arg z\leq\pi-\delta\}$.
	\end{enumerate}
\end{theorem}

\medskip

\noindent
Note that the integral in \eqref{A58} is always finite.
This is a consequence of \cite[Theorem~4.1]{langer.woracek:expty},
which also implies that this integral is equal to the exponential type of each entry
of $\uptheta(x;\cdot)$ and $\upvarphi(x;\cdot)$.

\begin{proofof}{Theorem~\ref{A31}}
	This theorem is a consequence of the indefinite version of
	\cite[Theorem~1.2]{langer.woracek:lokinv} indicated in \cite[Remark~1.3]{langer.woracek:lokinv}.
	Let $H_1,H_2\in\bb H$ be given, and assume w.l.o.g.\ that both are defined on $(0,\infty)$.
	The proof again proceeds via considering general Hamiltonians $\mf h_1$ and $\mf h_2$ built in our basic
	identification \ref{A108} from $H_1$ and $H_2$, respectively.

	Assuming (i) we choose $\oe_1,b_{1,j},d_{1,j}$ all equal to $0$ and the same base point $x_0$
	in the definition of $\mf h_1$ and $\mf h_2$.  Then, by \cite[Theorem~1.2 (indefinite variant)]{langer.woracek:lokinv},
	it follows that $q_{\mf h_1}$ and $q_{\mf h_2}$ are exponentially close in the sense of (ii) and (iii).

	Conversely, assume (ii) or (iii), and choose $\oe^{(1)}_1$, $b^{(1)}_{1,j}$, $d^{(1)}_{1,j}$
	and $x_1$ and $\oe^{(2)}_1$, $b^{(2)}_{1,j}$, $d^{(2)}_{1,j}$ and $x_2$ in the definition
	of $\mf h_1$ and $\mf h_2$ so that
	$q_{H_i}=q_{\mf h_i}$.  This is possible; cf.\ Remark~\ref{A96}. Then, again
	by \cite[Theorem~1.2 (indefinite variant)]{langer.woracek:lokinv}, $\mf h_{1,\Lsh s_1(\tau)}$ and
	$\mf h_{2,\Lsh s_2(\tau)}$ are reparameterizations of each other. In particular, (i) holds.
\end{proofof}

\begin{remark}\thlab{A103}
	If one (and hence all) of the equivalent conditions of Theorem~\ref{A31} hold, then
	\eqref{A104} holds with
	\begin{equation}\label{A90}
		k\defeq8\max\{\Delta(H_1),\Delta(H_2)\}+3.
	\end{equation}
	This can be seen by tracing the proof of
	\cite[Theorem~1.2 (indefinite version)]{langer.woracek:lokinv} as indicated in the
	footnotes in this paper.
	
	The value \eqref{A90} of the constant $k$ in \eqref{A104} is probably not the best possible.
	However, it is noteworthy that \eqref{A90} depends only on $\Delta(H_1)$ and $\Delta(H_2)$.
\end{remark}

\begin{remark}\thlab{A100}
	The dependence of these results on the choices of $q_{H_1}$ and $q_{H_2}$ is not essential.
	If {\rm(ii/iii)} holds with some pair $(q_{H_1},q_{H_2})\in[q]_{H_1}\times[q]_{H_2}$, then
	for each $q_1\in[q]_{H_1}$ there exists a unique element $q_2\in[q]_{H_2}$ such that
	{\rm(ii/iii)} holds for $(q_1,q_2)$. This corresponding function $q_2$ can be determined
	by starting with some $q\in[q]_{H_2}$, computing the polynomial asymptotics of $q-q_{H_1}$
	at $i\infty$, and subtracting this polynomial from $q$.
\end{remark}

Viewing the above remark from a slightly different perspective leads to the following more
effective test for {\rm(ii/iii)} to hold, which removes the dependence on the choice
of $q_{H_1}$ and $q_{H_2}$.

\begin{corollary}\thlab{A105}
	Assume that we are in the situation of Theorem~\ref{A31}.
	Pick some singular Weyl coefficients $q_1\in[q]_{H_1}$ and $q_2\in[q]_{H_2}$ and let
	\[
		q_2(iy)-q_1(iy)=\alpha_ny^n+\alpha_{n-1}y^{n-1}+\ldots+\alpha_1y+{\rm o}(y),\quad y\to\infty.
	\]
	Then {\rm(ii/iii)} of Theorem~\ref{A31} hold if and only if
	the conditions stated in {\rm(ii/iii)} hold with $q_{H_1}=q_1$ and
	\[
		q_{H_2}(z) = q_2(z)-\sum_{l=1}^n\alpha_l(-i)^lz^l.
	\]
	\popQED\qed
\end{corollary}

\noindent
The following corollary of Theorem~\ref{A31} is also worth mentioning.
It says that under the a priori hypothesis of finite exponential type, the global
uniqueness result Theorem~\ref{A6} can be strengthened; and it may be
much easier to establish exponential closeness of singular Weyl coefficients than
their actual equality.

\begin{corollary}\thlab{A59}
	Let $H_1,H_2\in\bb H$ with $\dom(H_i)=(a_i,b_i)$, $i=1,2$, be given and assume that
	\[
		\intop_{a_i}^{b_i}\sqrt{\det H_i(y)}\;\rd y<\infty,\qquad i=1,2.
	\]
	If there exist singular Weyl coefficients $q_{H_1}$ and $q_{H_2}$ of $H_1$ and $H_2$,
	respectively, and there exist $\beta\in(0,\pi)$ and
	$\tau>\max\bigl\{\int_{a_i}^{b_i}\sqrt{\det H_i(y)}\,\rd y:\ i=1,2\bigr\}$
	such that
	\[
		q_{H_1}(re^{i\beta})-q_{H_2}(re^{i\beta})=
		{\rm O}\bigl(e^{-2\tau r\sin\beta}\bigr),
		\qquad r\to\infty,
	\]
	then $H_1$ and $H_2$ are reparameterizations of each other.
	\popQED\qed
\end{corollary}

\bigskip

%
%
\begin{center}
{\Large\textbf{PART II: \\[1ex] Applications to Sturm--Liouville Equations}}
\end{center}
%
%

\noindent
In the remaining sections \ref{sec-SL}--\ref{sec-Schroedinger} we study
scalar second-order differential equations.
Under certain assumptions such equations can be transformed to canonical systems
of the form \eqref{A30} so that the results
from Sections~\ref{sec-measure}--\ref{sec-inverse} can be applied.
In Sections~\ref{sec-SL} and \ref{sec-SLw} we consider Sturm--Liouville equations
of the form \eqref{A208} below where either $1/p$ or $w$ is not integrable at $a$.
In Section~\ref{sec-Schroedinger} we study one-dimensional Schr\"odinger operators
with a singular potential.

%
%
\section{Sturm--Liouville equations without \\ potential: singular {\boldmath$1/p$}}
\label{sec-SL}
%
%

In this section we consider Sturm--Liouville equations of the form
\begin{equation}\label{A208}
	-\bigl(py'\bigr)' = \lambda wy
\end{equation}
on an interval $(a,b)$ with $-\infty\le a<b\le\infty$ where $\lambda\in\bb C$
and the functions $p$ and $w$ satisfy the conditions
\begin{equation}\label{A210}
	p(x)>0,\, w(x)>0 \;\;\text{a.e.}, \qquad \frac{1}{p}\,,w \in L^1_{\rm loc}(a,b).
\end{equation}
In the following we write $\dom(p;w)\defeq(a,b)$.
Moreover, let $L^2(w)$ be the weighted $L^2$-space with
inner product $(f,g)=\int_a^b f\qu{g}w$.

We consider the following class of coefficients.

\begin{definition}\thlab{A213}
	We say that $(p;w)\in\KSL$ if $p$ and $w$ are defined on some interval $(a,b)$
	and they satisfy \eqref{A210} and the following conditions.
	\begin{enumerate}[{\rm(i)}]
	\item
		For one (and hence for all) $x_0\in(a,b)$,
		\begin{equation}\label{A202}
			\intop_a^{x_0} \frac{1}{p(x)}\rd x = \infty
			\qquad\text{and}\qquad
			\intop_a^{x_0} w(x)\rd x < \infty.
		\end{equation}
	\item
		For one (and hence for all) $x_0\in(a,b)$,
		\begin{equation}\label{A203}
			\intop_a^{x_0} \intop_a^x w(t)\rd t\,\frac{1}{p(x)}\rd x < \infty.
		\end{equation}
	\item
		Let $x_0\in(a,b)$ and define functions $w_l$, $l=0,1,\dots$, recursively by
		\begin{equation}\label{A206}
		\begin{aligned}
			w_0(x) &= 1,
				\\[1ex]
			w_{l}(x) &= \begin{cases}
				\displaystyle \,\intop_x^{x_0} \frac{1}{p(t)}w_{l-1}(t)\rd t \quad
				& \text{if $l$ is odd},
				\\[4ex]
				\displaystyle \,\intop_a^x w(t)w_{l-1}(t)\rd t & \text{if $l$ is even}.
			\end{cases}
		\end{aligned}
		\end{equation}
		There exists an $n\in\bb N_0$ such that
		\begin{equation}\label{A209}
			w_n\big|_{(a,x_0)} \in \begin{cases}
				L^2\bigl(\frac{1}{p}\big|_{(a,x_0)}\bigr) & \text{if $n$ is even},
				\\[2ex]
				L^2\bigl(w\big|_{(a,x_0)}\bigr) & \text{if $n$ is odd}.
			\end{cases}
		\end{equation}
	\item
		Equation \eqref{A208} is in the limit point case at $b$, i.e.\
		for $\lambda\in\bb C\setminus\bb R$, equation \eqref{A208} has
		(up to a scalar multiple) only one solution in $L^2(w|_{(x_0,b)})$
		for $x_0\in(a,b)$.
	\end{enumerate}
	If $(p;w)\in\KSL$, we denote by $\DeltaSL(p,w)$ the minimal $n\in\bb N_0$ such that
	\eqref{A209} holds.
\end{definition}

\begin{remark}\thlab{A207}
	\rule{0ex}{1ex}
	\begin{enumerate}[{\rm(i)}]
	\item
		Under the assumption of \eqref{A202}, condition \eqref{A203} is equivalent to
		\[
			\intop_a^{x_0} \intop_x^{x_0} \frac{1}{p(t)}\,\rd t\,w(x)\rd x < \infty;
		\]
		see, e.g.\ \cite[Lemma~4.3]{langer.woracek:gpinf}.
	\item
		Assume that \eqref{A202} holds.  Then \eqref{A203} and \eqref{A209} with $n=1$
		are satisfied
		if and only if equation \eqref{A208} is in the limit circle case at $a$;
		this is true because the solutions of \eqref{A208}
		with $\lambda=0$ are $y(x)=c_1 w_1(x)+c_2$ with $c_1,c_2\in\bb C$
		and the limit circle case prevails at $a$ if and only if all these
		solutions are in $L^2(w|_{(a,x_0)})$.
	\item
		The functions $w_0$ and $w_1$ are solutions of \eqref{A208} with $\lambda=0$.
		Since $w_1(x)\to\infty$ as $x\searrow a$, the function $w_0$ is a principal
		solution and $w_1$ is a non-principal solution, i.e.\ $w_0(x)=\mathrm{o}(w_1(x))$ as $x\searrow0$;
		for the notions of principal and non-principal solutions see, e.g.\ \cite{niessen.zettl:1992}.
		Moreover, one can easily verify that
		\[
			-\frac{1}{w}\bigl(pw_{l+2}'\bigr)' = w_l
			\qquad \text{when $l\in\bb N$ is odd}.
		\]
	\end{enumerate}
\end{remark}

\medskip

\noindent
For given $p$ and $w$ satisfying \eqref{A210} define the Hamiltonian
\begin{equation}\label{A234}
	H(x) \defeq \begin{pmatrix} \dfrac{1}{p(x)} & 0 \\[3ex] 0 & w(x) \end{pmatrix},
	\qquad x\in(a,b).
\end{equation}
If $\uppsi=(\uppsi_1,\uppsi_2)^T$ is a solution of equation \eqref{A30}
with $H$ as in \eqref{A234}, then
\[
	\uppsi_1' = -zw\uppsi_2, \qquad
	\uppsi_2' = z\frac{1}{p}\uppsi_1,
\]
and hence $\uppsi_2$ is a solution of \eqref{A208} with $\lambda=z^2$.
Conversely, if $\psi$ is a solution of \eqref{A208} and $z\in\bb C$ is such
that $z^2=\lambda$, then
\begin{equation}\label{A211}
	\uppsi(x) = \begin{pmatrix} p(x)\psi'(x) \\[1ex] z\psi(x) \end{pmatrix}
\end{equation}
satisfies \eqref{A30} with $H$ as in \eqref{A234}.

In the following assume that $(p;w)\in\KSL$.
The first relation in \eqref{A202} implies that $H$ is in the limit point case at $a$.
Since \eqref{A208} is in the limit point case at $b$,
the Hamiltonian $H$ is also in the limit point case at $b$
because $\uppsi\in L^2(H|_{(x_0,b)})$ with $\uppsi$ as in \eqref{A211}, $z\ne0$
and $x_0\in(a,b)$ implies that $\psi\in L^2(w|_{(x_0,b)})$.
Therefore the operator $T(H)$, which is defined in \eqref{A143}
and acts in the space $L^2(H)=L^2(\frac{1}{p})\oplus L^2(w)$, is self-adjoint.
Since $H(x)$ is invertible for a.e.\ $x\in(a,b)$,
the operator $T(H)$ can be written as
\[
	T(H)f = H^{-1}J^{-1}f' = \begin{pmatrix} pf_2' \\[1ex] -\dfrac{1}{w}f_1' \end{pmatrix}
\]
with maximal domain
\[
	\dom(T(H)) = \biggl\{\binom{f_1}{f_2}: f_1,f_2 \;\;\text{abs.\ cont.},\;
	pf_2'\in L^2\Bigl(\frac{1}{p}\Bigr),\; \frac{1}{w}f_1'\in L^2(w)\biggr\}
\]
Hence $(T(H))^2$ acts as follows
\begin{equation}\label{A250}
	\bigl(T(H)\bigr)^2\binom{f_1}{f_2}
	= \begin{pmatrix}
		-p\Bigl(\dfrac{1}{w}f_1'\Bigr)'\, \\[2ex]
		-\dfrac{1}{w}\bigl(pf_2'\bigr)'
	\end{pmatrix}.
\end{equation}
With the mappings
\begin{equation}\label{A251}
	\iota_2 : \begin{cases}
		L^2(w) \to L^2(H) \\[1ex]
		\displaystyle g \mapsto \binom{0}{g},
	\end{cases}
	\qquad
	P_2 : \begin{cases}
		L^2(H) \to L^2(w) \\[1ex]
		\displaystyle \binom{f}{g} \mapsto g,
	\end{cases}
\end{equation}
we define the self-adjoint operator
\begin{equation}\label{A286}
	A_{p,w} \defeq P_2\bigl(T(H)\bigr)^2\iota_2.
\end{equation}
This operator acts like
\begin{align}
	A_{p,w}y &= -\frac{1}{w}\bigl(py'\bigr)'
	\notag\\[1ex]
	\dom(A_{p,w}) &= \Bigl\{y\in L^2(w): y,py'\text{ locally absolutely continuous},
	\notag\\[1ex]
	&\qquad \int_a^b p(x)\bigl|y'(x)\bigr|^2\rd x<\infty,\;\frac{1}{w}\bigl(py'\bigr)'\in L^2(w)\Bigr\}
	\label{A289}
\end{align}
and is the Friedrichs extension of the minimal operator associated with \eqref{A208}
since all functions in $\dom(A_{p,w})$ are in the form domain;
note that $A_{p,w}$ is non-negative.
If \eqref{A208} is also in the limit point case at $a$ (that is, when $\DeltaSL(p,w)\ge2$),
then $A_{p,w}$ coincides with the maximal operator,
i.e.\ the condition $\int_a^b p|y'|^2<\infty$ is automatically satisfied.
If \eqref{A208} is in the limit circle case at $a$, one can replace the
condition $\int_a^b p|y'|^2<\infty$ in \eqref{A289} by any of the two boundary conditions
\begin{equation}\label{A290}
	\lim_{x\searrow a}\frac{y(x)}{w_1(x)} = 0, \qquad
	\lim_{x\searrow a}p(x)y'(x) = 0;
\end{equation}
see, e.g.\ \cite[Theorem~4.3]{niessen.zettl:1992}.

\begin{remark}\thlab{A274}
	One can also treat the situation when \eqref{A208} is either regular or in the
	limit circle case at $b$.  In the former case one extends $H$ by an indivisible interval
	of infinite length; in the latter case $H$ is in the limit point case.
	In both cases elements in the domain of $A_{p,w}$ defined via \eqref{A286}
	satisfy some boundary condition at $b$.
\end{remark}

Assume that $(p;w)\in\KSL$ and let $H$ be as in \eqref{A234}.
It follows from \cite[Theorem~3.7]{winkler.woracek:del} that $H\in\bb H$,
that the functions $\mf w_l$, $l\in\bb N_0$, defined
in \eqref{A114} are given by
\begin{equation}\label{A296}
	\mf w_l(x)
	= \begin{cases}
		\begin{pmatrix} w_l(x) \\ 0 \end{pmatrix} & \text{if $l$ is even}, \\[3ex]
		\begin{pmatrix} 0 \\ -w_l(x) \end{pmatrix} & \text{if $l$ is odd},
	\end{cases}
\end{equation}
and that
\begin{equation}\label{A247}
	\DeltaSL(p,w)=\Delta(H).
\end{equation}
Therefore we can apply the results from Part~I to the Hamiltonian $H$.
Using the connection between \eqref{A208} and \eqref{A30} we can show that
regularized boundary values of solutions of \eqref{A208} exist at $a$.
Denote by $\NSL{\lambda}$ the set of all solutions of the
Sturm--Liouville equation \eqref{A208}.

\begin{theorem}[\textbf{Regularized boundary values}]\thlab{A220}
	Let $(p;w)\in\KSL$ with $\dom(p;w)=(a,b)$ and set $\Delta\defeq\DeltaSL(p,w)$.
	Then, for $x_0\in(a,b)$, the following statements hold.
	\begin{enumerate}[{\rm(i)}]
	\item
		For each $\lambda\in\bb C$ and each solution $\psi\in\NSL{\lambda}$
		the boundary value
		\begin{equation}\label{A221}
			\rbvSLr{\lambda}\psi \defeq \lim_{x\searrow a}p(x)\psi'(x)
		\end{equation}
		and the regularized boundary value
		\begin{align}
			\rbvSLs{\lambda}\psi &\defeq \lim_{x\searrow a}\Biggl[
				\,\sum_{k=0}^{\lfloor\frac{\Delta-1}{2}\rfloor}
				\lambda^k\Bigl(w_{2k}(x)\psi(x)+w_{2k+1}(x)p(x)\psi'(x)\Bigr)
				\label{A222}
				\\[1ex]
			&\quad +
				\left\{\begin{array}{ll}
					\lambda^{\frac{\Delta}{2}} w_\Delta(x)\psi(x) & \text{if $\Delta$ is even}  \\[1ex]
					0 & \text{if $\Delta$ is odd}
				\end{array}\right\}
				\notag\\[1ex]
			&\quad + \Bigl(\lim_{t\searrow a}p(t)\psi'(t)\Bigr)
				\sum_{k=\lfloor\frac{\Delta+1}{2}\rfloor}^{\Delta-1}
				\sum_{l=0}^{2k-\Delta}
				(-1)^l \lambda^k w_l(x)w_{2k-l+1}(x)\Biggr]
				\hspace*{5ex}
				\notag
		\end{align}
		exist.
	\item
		For each $\lambda\in\bb C$ we define
		\[
			\rbvSL{\lambda}:
			\left\{
			\begin{array}{ccl}
			\NSL{\lambda} & \to & \bb C^2 \\[1ex]
			\psi & \mapsto & \bigr(\rbvSLr{\lambda}\psi,\;\rbvSLs{\lambda}\psi\bigr)^T.
			\end{array}
			\right.
		\]
		Then $\rbvSL{\lambda}$ is a bijection from $\NSL{\lambda}$ onto $\bb C^2$.
	\item
		For each $\lambda\in\bb C$ there exists an {\rm(}up to scalar multiples{\rm)} unique
		solution $\psi\in\NSL{\lambda}\setminus\{0\}$ such that
		$\lim_{x\searrow a}\psi(x)$ exists.

		This solution is characterized by the property that
		$\int_a^{x_0}p|\psi'|^2<\infty$
		and also by the property that
		$\rbvSLr{\lambda}\psi=0$ {\rm(}and $\psi\not\equiv0${\rm)}.

		If $\psi$ is a solution such that $\lim\limits_{x\searrow a}\psi(x)$ exists, then
		$\rbvSLs{\lambda}\psi = \lim\limits_{x\searrow a}\psi(x)$.
	\end{enumerate}
	The regularized boundary value $\rbvSLs{\lambda}$ depends on the choice of\, $x_0$
	in the following way.
	\begin{enumerate}[{\rm(i)}]
	\setcounter{enumi}{3}
	\item Let $x_0,\hat x_0\in(a,b)$, and let $\rbvSL{\lambda}$ and $\rbvSLpr{\lambda}$
		be the correspondingly defined regularized boundary value mappings.
		Then there exists a polynomial $\pxxh$ with real coefficients whose degree
		does not exceed $\Delta-1$ such that
		\[
			\rbvSLspr{\lambda}\psi = \rbvSLs{\lambda}\psi
			+ \pxxh(\lambda)\rbvSLr{\lambda}\psi,\qquad \psi\in\NSL{\lambda},\;\; \lambda\in\bb C.
		\]
		Moreover, clearly, $\rbvSLrpr{\lambda} = \rbvSLr{\lambda}$.
	\end{enumerate}
\end{theorem}

\medskip

\begin{remark}\thlab{A288}\rule{0ex}{1ex}
\begin{enumerate}[{\rm(i)}]
\item
	Let
	\begin{equation}\label{A317}
		W_p(y_1,y_2)(x)\defeq p(x)\Bigl(y_1(x)y_2'(x)-y_1'(x)y_2(x)\Bigr)
	\end{equation}
	be the weighted Wronskian with weight $p$.  Using \eqref{A206} we can rewrite
	the expression that appears within the round brackets in \eqref{A222} as follows:
	\[
		w_{2k}\psi+w_{2k+1}p\psi' = W_p(w_{2k+1},\psi).
	\]
\item
	When $\DeltaSL(p,w)=1$, i.e.\ when \eqref{A208} is in the limit circle case at $a$,
	then $\rbvSLs{\lambda}\psi$ does not depend on $\lambda$
	explicitly, and it takes the form
	\[
		\rbvSLs{\lambda}\psi = \lim_{x\searrow a}\Bigl(\psi(x)+p(x)w_1(x)\psi'(x)\Bigr)
		= \lim_{x\searrow a} W_p(w_1,\psi)(x).
	\]
	Note also that $\rbvSLr{\lambda}\psi = \lim_{x\searrow a} W_p(1,\psi)(x)$.
\item
	Instead of the functions $w_l$ one can use functions $\check{w}_l$
	that are defined by the recurrence relation $\check{w}_0\equiv1$ and
	\[
		\check{w}_{l}(x) = \begin{cases}
			\displaystyle \,\intop_x^{x_0} \frac{1}{p(t)}\check{w}_{l-1}(t)\rd t + c_l \quad
			& \text{if $l$ is odd},
			\\[4ex]
			\displaystyle \,\intop_a^x w(t)\check{w}_{l-1}(t)\rd t & \text{if $l$ is even},
		\end{cases}
	\]
	with arbitrary real numbers $c_l$ for odd $l$.
	To add the extra constants $c_l$ is useful for practical calculations, in particular,
	when $w_l$ has an asymptotic expansion (for $x\searrow a$) in which a constant
	term can be removed by adjusting $c_l$.
	One can show that the corresponding regularized boundary value
	$\check\RbvSL_{\lambda,2}^{\rm SL}$ satisfies
	\[
		\check{\RbvSL}_{\lambda,2}^{\rm SL}\psi
		= \rbvSLs{\lambda}\psi+\check p(\lambda)\rbvSLr{\lambda}\psi
	\]
	with
	\begin{equation}\label{A307}
		\check p(\lambda) = \sum_{k=0}^{\Delta-1} \lambda^k \;
		\sum_{i=0}^k c_{2k+1-2i}\lim_{t\searrow a}v_{2i}(t)
	\end{equation}
	where $v_0\equiv1$ and
	\[
		v_l(x) = \intop_x^{x_0} \intop_a^x w(t)v_{l-2}(t)\rd t\,\frac{1}{p(x)}\rd x,
		\qquad l\;\;\text{even}.
	\]
	The limit $\lim_{t\searrow a}v_l(t)$ exists because of condition \eqref{A203}.
\end{enumerate}
\end{remark}

\medskip

\noindent
Before we prove Theorem~\ref{A220}, we show the following lemma,
where $\mf N_z$ denotes the set of all solutions of \eqref{A30};
see \S\ref{A102}.

\begin{lemma}\thlab{A235}
	Let $z\in\bb C$ and let $\psi$ be a solution of \eqref{A208} with $\lambda\defeq z^2$.
	Moreover, set
	\begin{equation}\label{A236}
		\uppsi(x) \defeq \begin{pmatrix} p(x)\psi'(x) \\[1ex] z\psi(x) \end{pmatrix}, \qquad
		\widehat\uppsi(x) \defeq \begin{pmatrix} p(x)\psi'(x) \\[1ex] -z\psi(x) \end{pmatrix}.
	\end{equation}
	Then
	\begin{equation}\label{A295}
		\uppsi\in\mf N_z,\qquad \widehat\uppsi\in\mf N_{-z},\qquad
		\rbvr\uppsi=\Rbv_{-z,1}\widehat\uppsi
	\end{equation}
	and
	\begin{align*}
		\rbvs\uppsi &= -\Rbv_{-z,2}\widehat\uppsi
			\\[1ex]
		&= z\lim_{x\searrow a}\Biggl[
			\,\sum_{\substack{l:\,0\le l\le\Delta-1 \\[0.5ex] l\text{ \rm even}}}\hspace*{-1ex}
			z^l\Bigl(w_l(x)\psi(x)+w_{l+1}(x)p(x)\psi'(x)\Bigr)
			\\[1ex]
		&\quad + \left\{\begin{array}{ll}
				z^\Delta w_\Delta(x)\psi(x) & \text{if $\Delta$ is even}  \\[1ex]
				0 & \text{if $\Delta$ is odd}
			\end{array}\right\}
			\\[1ex]
		&\quad + \Bigl(\lim_{t \searrow a}p(t)\psi'(t)\Bigr)
			\sum_{k=\lfloor\frac{\Delta+1}{2}\rfloor}^{\Delta-1}
			\sum_{l=0}^{2k-\Delta}
			(-1)^l z^{2k} w_l(x)w_{2k-l+1}(x)\Biggr].
	\end{align*}
\end{lemma}

\begin{proof}
	Let $H$ be as in \eqref{A234} and set $\Delta\defeq\Delta(H)$.
	The relations in \eqref{A295} are clear from the considerations around equation \eqref{A211}
	and the fact that \eqref{A208} does not change when we replace $z$ by $-z$.
	From \eqref{A106} and \eqref{A296} we obtain
	\begin{align*}
		\rbvs\uppsi &= -\lim_{x\searrow a}\Biggl[\,\sum_{l=0}^\Delta
			z^l\bigl(\mf w_l(x)\bigr)^* J\biggl(\uppsi(x)-(\rbvr\uppsi)
			\sum_{k=\Delta+1}^{2\Delta-l}z^k\mf w_k(x)\biggr)\Biggr]
			\\[1ex]
		&= \lim_{x\searrow a}\Biggl[\,\sum_{\substack{l:\,0\le l\le\Delta \\[0.5ex] l\text{ even}}}
			z^l w_l(x)\biggl(\uppsi(x)-(\rbvr\uppsi)
			\sum_{k=\Delta+1}^{2\Delta-l}z^k\mf w_k(x)\biggr)_{\!2}
			\\[1ex]
		&\quad + \sum_{\substack{l:\,1\le l\le\Delta \\[0.5ex] l\text{ odd}}}
			z^l w_l(x)\biggl(\uppsi(x)-(\rbvr\uppsi)
			\sum_{k=\Delta+1}^{2\Delta-l}z^k\mf w_k(x)\biggr)_{\!1}\Biggr]
			\displaybreak[0]\\[1ex]
		&= \lim_{x\searrow a}\Biggl[\,\sum_{\substack{l:\,0\le l\le\Delta \\[0.5ex] l\text{ even}}}
			z^l w_l(x)\biggl(z\psi(x)+\Bigl(\lim_{t\searrow a}p(t)\psi'(t)\Bigr)\hspace*{-1ex}
			\sum_{\substack{k:\,\Delta+1\le k\le 2\Delta-l \\[0.5ex] k\text{ odd}}}\hspace*{-1.5ex}
			z^k w_k(x)\biggr)
			\\[1ex]
		&\quad + \sum_{\substack{l:\,1\le l\le\Delta \\[0.5ex] l\text{ odd}}}
			z^l w_l(x)\biggl(p(x)\psi'(x)-\Bigl(\lim_{t\searrow a}p(t)\psi'(t)\Bigr)\hspace*{-1ex}
			\!\sum_{\substack{k:\,\Delta+1\le k\le 2\Delta-l \\[0.5ex] k\text{ even}}}\hspace*{-1.5ex}
			z^k w_k(x)\biggr)\Biggr]
			\displaybreak[0]\\[1ex]
		&= z\lim_{x\searrow a}\Biggl[\,\sum_{\substack{l:\,0\le l\le\Delta \\[0.5ex] l\text{ even}}}
			z^l w_l(x)\psi(x)
			+ \hspace*{-1ex}\sum_{\substack{l:\,0\le l\le\Delta-1 \\[0.5ex] l\text{ even}}}
			z^l w_{l+1}(x)p(x)\psi'(x)
			\\[1ex]
		&\quad + \Bigl(\lim_{t \searrow a}p(t)\psi'(t)\Bigr)\biggl(\,
			\sum_{\substack{l:\,0\le l\le\Delta \\[0.5ex] l\text{ even}}}
			\hspace*{1ex}\sum_{\substack{k:\,\Delta+1\le k\le 2\Delta-l \\[0.5ex] k\text{ odd}}}\hspace*{-1.5ex}
			z^{l+k-1}w_l(x)w_k(x)
			\\[1ex]
		&\quad - \sum_{\substack{l:\,1\le l\le\Delta \\[0.5ex] l\text{ odd}}}
			\hspace*{1ex}\sum_{\substack{k:\,\Delta+1\le k\le 2\Delta-l \\[0.5ex] k\text{ even}}}\hspace*{-1.5ex}
			z^{l+k-1}w_l(x)w_k(x)\biggr)\Biggr]
			\displaybreak[0]\\[1ex]
		&= z\lim_{x\searrow a}\Biggl[\,\sum_{\substack{l:\,0\le l\le\Delta \\[0.5ex] l\text{ even}}}
			z^l w_l(x)\psi(x)
			+ \hspace*{-1ex}\sum_{\substack{l:\,0\le l\le\Delta-1 \\[0.5ex] l\text{ even}}}
			z^l w_{l+1}(x)p(x)\psi'(x)
			\\[1ex]
		&\quad + \Bigl(\lim_{t \searrow a}p(t)\psi'(t)\Bigr)\sum_{l=0}^\Delta
			\hspace*{1ex}\sum_{\substack{k:\,\Delta+1\le k\le 2\Delta-l \\[0.5ex] l+k\text{ odd}}}\hspace*{-1.5ex}
			(-1)^l z^{l+k-1}w_l(x)w_k(x)\Biggr]
			\displaybreak[0]\\[1ex]
		&= z\lim_{x\searrow a}\Biggl[
			\,\sum_{\substack{l:\,0\le l\le\Delta-1 \\[0.5ex] l\text{ even}}}\hspace*{-1ex}
			z^l\Bigl(w_l(x)\psi(x)+w_{l+1}(x)p(x)\psi'(x)\Bigr)
			\\[1ex]
		&\quad + \left\{\begin{array}{ll}
				z^\Delta w_\Delta(x)\psi(x) & \text{if $\Delta$ is even}  \\[1ex]
				0 & \text{if $\Delta$ is odd}
			\end{array}\right\}
			\\[1ex]
		&\quad + \Bigl(\lim_{t \searrow a}p(t)\psi'(t)\Bigr)
			\sum_{\substack{m:\,\Delta\le m\le2\Delta-2 \\[0.5ex] m\text{ even}}}
			\sum_{l=0}^{m-\Delta}
			(-1)^l z^m w_l(x)w_{m-l+1}(x)\Biggr],
	\end{align*}
	which proves the statement for $\rbvs\uppsi$.
	Inside the limit only even powers of $z$ appear,
	and hence, as $\Rbv_{-z,2}\widehat\uppsi$ is obtained from $\rbvs\uppsi$ by
	replacing $z$ by $-z$, the equality $\rbvs\uppsi = -\Rbv_{-z,2}\widehat\uppsi$ follows.
\end{proof}

\begin{proof}[Proof of Theorem~\ref{A220}]
	First we settle the case $\lambda=0$.  The solutions of \eqref{A208} with $\lambda=0$
	are of the form
	\[
		\psi(x) = c_1\int_{x_0}^x \frac{\rd t}{p(t)} + c_2
	\]
	with $c_1,c_2\in\bb C$.
	For such a solution the limits in \eqref{A221} and \eqref{A222} exist
	and $\rbvSLr{0}\psi=c_1$ and
	\begin{align*}
		\rbvSLs{0}\psi
		&= \lim_{x\searrow a}\bigl(w_0(x)\psi(x)+w_1(x)p(x)\psi'(x)\bigr)
		\\[1ex]
		&= \lim_{x\searrow a}\biggl(c_1\int_{x_0}^x\frac{\rd t}{p(t)}+c_2
		+\int_x^{x_0}\frac{\rd t}{p(t)}\cdot p(x)\cdot\frac{c_1}{p(x)}\biggr)
		= c_2.
	\end{align*}
	This shows that $\rbvSL{0}:\NSL{0}\to\bb C^2$ is a bijective mapping.
	Moreover, the conditions in (iii) are all equivalent to $c_1=0$ since $\frac1p$
	is not integrable at $a$.

	For the rest of the proof assume that $\lambda\ne0$.

	(i)\,
	Let $\psi$ be a solution of \eqref{A208}, let $z\in\bb C$ with $z^2=\lambda$
	and define $\uppsi$ as in \eqref{A236}.  The existence of the limit in \eqref{A221}
	and the equality
	\begin{equation}\label{A226}
		\rbvSLr{\lambda}\psi = \rbvr\uppsi
	\end{equation}
	are immediate.  The existence of the limit in \eqref{A222} and the relation
	\begin{equation}\label{A237}
		\rbvSLs{\lambda}\psi = \frac{1}{z}\rbvs\uppsi
	\end{equation}
	follows from Lemma~\ref{A235} by observing that $z^2=\lambda$.

	(ii)\,
	Theorem~\ref{A32}\,(ii) and the relations
	in \eqref{A226} and \eqref{A237} show that the
	mapping $\rbvSL{\lambda}:\NSL{\lambda}\to\bb C^2$ is bijective.

	(iii)\,
	The first and the last assertions follow immediately from Theorem~\ref{A32}\,(iii).
	For the second statement note that there is (up to a scalar multiple) a unique solution
	$\uppsi_{\rm reg}$ such that $\uppsi_{\rm reg}|_{(a,x_0)}\in L^2(H|_{(a,x_0)})$.
	Hence $p\psi_{\rm reg}'|_{(a,x_0)}=\uppsi_{\rm reg,1}|_{(a,x_0)}\in L^2\bigl(\frac{1}{p}\big|_{(a,x_0)}\bigr)$.
	Any other solution $\psi$ is such that $\lim_{x\searrow}p(x)\psi'(x)\ne0$ according
	to the already proved third statement of (iii).
	Since $\frac{1}{p}$ is not integrable at $a$ by assumption, such a $\psi$
	satisfies $p\psi'|_{(a,x_0)}\notin L^2\bigl(\frac{1}{p}\big|_{(a,x_0)}\bigr)$.
	Now the claim follows because $p\psi'|_{(a,x_0)}\in L^2\bigl(\frac{1}{p}\big|_{(a,x_0)}\bigr)$
	if and only if $\int_a^{x_0}p|\psi'|^2<\infty$.

	(iv)\,
	It follows from Theorem~\ref{A32}\,(iv) that there exists a polynomial $\hat p$ of degree
	at most $2\Delta$ with real coefficients and no constant term such that
	\begin{equation}\label{A238}
		\rbvspr{z}\uppsi = \Rbv_{z,2}\uppsi + \hat p(z)\Rbv_{z,1}\uppsi
	\end{equation}
	for all $\uppsi\in\mf N_z$.  If we choose $\uppsi$ as in \eqref{A236} for $\psi\in\NSL{z^2}$,
	then, by \eqref{A226} and \eqref{A237}, we have
	\begin{equation}\label{A239}
		\begin{aligned}
			\rbvSLspr{z^2}\psi &= \frac{1}{z}\rbvspr{z}\uppsi
			= \frac{1}{z}\Bigl(\Rbv_{z,2}\uppsi + \hat p(z)\Rbv_{z,1}\uppsi\Bigr)
				\\[1ex]
			&= \rbvSLs{z^2}\psi + \frac{\hat p(z)}{z}\rbvSLr{z^2}\psi.
		\end{aligned}
	\end{equation}
	Since this relation must be true for all $z\in\bb C\setminus\{0\}$ and all $\psi\in\NSL{z^2}$,
	it follows by replacing $z$ by $-z$ that $\hat p$ is an odd polynomial.
	Hence one can define a polynomial $\pxxh$ by the relation $\pxxh(z^2)=\frac{\hat p(z)}{z}$,
	which is a real polynomial of degree at most $\Delta-1$.
	Now the assertion follows from \eqref{A239}.
\end{proof}

\medskip

In the next theorem we establish the existence of a singular
Titchmarsh--Weyl coefficient, which is used in Theorem~\ref{A227} below to
obtain a spectral measure.
Recall from Definition~\ref{A115} that $\mc N_\kappa$,
$\kappa\in\bb N_0$, is the set of all generalized Nevanlinna functions
with $\kappa$ negative squares and that
$\mc N_{<\infty}=\bigcup_{\kappa\in\bb N_0}\mc N_\kappa$.
Further, denote by $\mc N_\kappa^{(\infty)}$, $\kappa\in\bb N_0$, the set of
functions from $\mc N_\kappa$ whose only generalized pole of non-positive type
is infinity and set
$\mc N_{<\infty}^{(\infty)}=\bigcup_{\kappa\in\bb N_0}\mc N_\kappa^{(\infty)}$;
see Definition~\ref{A116}.
Note that a function $q$, defined on $\bb C\setminus\bb R$, belongs
to $\mc N_\kappa^{(\infty)}$ if and only if
\[
	q(z) = p_{2\kappa}(z)q_0(z)
\]
where $p_{2\kappa}$ is a monic real polynomial of degree $2\kappa$
and $q_0\in\mc N_0$, i.e.\ $q_0(\qu z)=\qu{q_0(z)}$ and $\Im q_0(z)\ge0$
for $z\in\bb C^+=\{z\in\bb C:\Im z>0\}$.

\begin{theorem}[\textbf{Singular Titchmarsh--Weyl coefficients}]\thlab{A223}
	Let $(p;w)\in\KSL$ with $\dom(p;w)=(a,b)$ be given.
	Then, for each fixed $x_0\in(a,b)$, the following statements hold.
	\begin{enumerate}[{\rm(i)}]
	\item
		For $\lambda\in\bb C$ let $\theta(\cdot\,;\lambda)$ and $\varphi(\cdot\,;\lambda)$
		be the unique solutions of \eqref{A208} such that
		\begin{equation}\label{A224}
			\rbvSL{\lambda}\theta(\cdot\,;\lambda)=\binom{1}{0}, \qquad
			\rbvSL{\lambda}\varphi(\cdot\,;\lambda)=\binom{0}{1}.
		\end{equation}
		Then, for each $x\in(a,b)$, the functions $\theta(x;\cdot)$ and $\varphi(x;\cdot)$
		are entire of order $\frac12$ and finite \textup{(}positive\textup{)} type
		\[
			\intop_a^x \sqrt{\frac{w(t)}{p(t)}}\,\rd t.
		\]
		Moreover, for each $\lambda\in\bb C$, one has
		$W_p\bigl(\varphi(\cdot\,;\lambda),\theta(\cdot\,;\lambda)\bigr)\equiv1$
		where the weighted Wronskian $W_p$ is as in \eqref{A317}, and the
		following relations hold:
		\begin{equation}\label{A292}
		\begin{alignedat}{2}
			\lim_{x\searrow a}\varphi(x;\lambda) &= 1,
			\qquad &
			\lim_{x\searrow a}\frac{p(x)\varphi'(x;\lambda)}{\int_a^x w(t)\rd t} &= -\lambda,
			\\[1ex]
			\lim_{x\searrow a}\frac{\theta(x;\lambda)}{w_1(x)} &= -1,
			\hspace*{7ex} &
			\lim_{x\searrow a}p(x)\theta'(x;\lambda) &= 1.
		\end{alignedat}
		\end{equation}
	\item
		The limit
		\begin{equation}\label{A225}
			m_{p,w}(\lambda) \defeq \lim_{x\nearrow b}
			\frac{\theta(x;\lambda)}{\varphi(x;\lambda)}\,, \qquad \lambda\in\bb C\setminus[0,\infty),
		\end{equation}
		exists locally uniformly on $\bb C\setminus[0,\infty)$ and defines an
		analytic function in~$\lambda$.
		The function $m_{p,w}$ belongs to the class $\mc N_\kappa^{(\infty)}$ with
		$\kappa=\bigl\lfloor\frac{\DeltaSL(p,w)}{2}\bigr\rfloor$.
	\item
		We have
		\[
			\theta(\cdot\,;\lambda)-m_{p,w}(\lambda)\varphi(\cdot\,;\lambda)\in L^2\bigl(w|_{(x_0,b)}\bigr),
			\qquad \lambda\in\bb C\setminus[0,\infty),	
		\]
		and this property characterizes the value $m_{p,w}(\lambda)$ for each $\lambda\in\bb C\setminus[0,\infty)$.
	\item
		For $\lambda\in\bb C\setminus[0,\infty)$ let $\psi$ be any non-trivial
		solution of \eqref{A208} such that $\psi|_{(x_0,b)}\in L^2(w|_{(x_0,b)})$.
		Then
		\[
			m_{p,w}(\lambda) = -\frac{\rbvSLs{\lambda}\psi}{\rbvSLr{\lambda}\psi}\,.
		\]
	\end{enumerate}
	The function $m_{p,w}$ depends on the choice of $x_0$.
	This dependence is controlled as follows.
	\begin{enumerate}[{\rm(i)}]
	\setcounter{enumi}{4}
	\item
		Let $\hat x_0\in(a,b)$, and let $\wh m_{p,w}$ be
		the correspondingly defined singular Titchmarsh--Weyl coefficient.
		Then there exists a polynomial $\pxxh$ with real coefficients whose degree
		does not exceed $\DeltaSL(p,w)-1$ such that
		\begin{equation}\label{A308}
			\wh m_{p,w}(\lambda) = m_{p,w}(\lambda)-\pxxh(\lambda).
		\end{equation}
	\end{enumerate}
\end{theorem}

\medskip

\noindent
Before we prove Theorem~\ref{A223}, let us introduce some notation.

\begin{definition}\thlab{A264}
	We refer to each function $m_{p,w}$ constructed
	as in Theorem~\ref{A223} as \emph{a singular Titchmarsh--Weyl coefficient}
	associated with the Sturm--Liouville equation \eqref{A208} when $(p;w)\in\KSL$.
	We denote by $[m]_{p,w}$ the equivalence class of $\mc N_{<\infty}^{(\infty)}$-functions
	modulo the relation
	\begin{equation}\label{A287}
		m_1 \;\hat{\sim}\; m_2 \quad\defequiv\quad m_1-m_2\in\bb R[z]
	\end{equation}
	which contains some (and hence any) function $m_{p,w}$ in Theorem~\ref{A223};
	we call $[m]_{p,w}$ \emph{the} singular Titchmarsh--Weyl coefficient.
\end{definition}

\medskip

\pagebreak[3]

\begin{remark}\thlab{A271}\rule{0ex}{1ex}
	\begin{enumerate}[{\rm(i)}]
	\item
		If $\DeltaSL(p,w)=1$, then $m_{p,w}\in\mc N_0$ by Theorem~\ref{A223}\,{\rm(ii)},
		which is in accordance with the classical theory since \eqref{A208}
		is in the limit circle case at $a$;
		see, e.g.\ \cite[Corollary~8.1]{bennewitz.everitt:2005} and
		\cite[Corollary~A.9]{kostenko.sakhnovich.teschl:2012}.
	\item
		According to Theorem~\ref{A220}\,{\rm(iii)}, $\varphi(\cdot\,;\lambda)$ is the
		--- up to a multiplicative scalar --- unique solution of \eqref{A208} that
		satisfies $\int_a^{x_0}p(x)|\varphi'(x;\lambda)|^2\rd x<\infty$ for some $x_0$.
		If \eqref{A208} is in the limit circle case at $a$ (i.e.\ if $\DeltaSL(p,w)=1$),
		then $\varphi(\cdot\,;\lambda)$ is the only solution of \eqref{A208}
		that satisfies the boundary condition \eqref{A290};
		if \eqref{A208} is in the limit point case, then $\varphi(\cdot\,;\lambda)$ is the
		only solution in $L^2(w|_{(a,x_0)})$.
	\end{enumerate}
\end{remark}

\noindent
For the proof of Theorem~\ref{A223} let $\theta(\cdot\,;\lambda)$ and $\varphi(\cdot\,;\lambda)$
be as in the statement of the theorem, i.e.\ the unique solutions of \eqref{A208} that satisfy \eqref{A224}.
Let $z\in\bb C\setminus\{0\}$ be such that $\lambda=z^2$ and let $H$ be as in \eqref{A234}.
Then the functions $\uptheta(\cdot\,;z), \upvarphi(\cdot\,;z)\in\mf N_z$ that satisfy \eqref{A110}
have the following form:
\begin{equation}\label{A254}
	\uptheta(x;z) = \begin{pmatrix} p(x)\theta'(x;z^2) \\[2ex] z\theta(x;z^2) \end{pmatrix},
	\qquad
	\upvarphi(x;z) = \begin{pmatrix} \,\dfrac{1}{z}p(x)\varphi'(x;z^2) \\[2ex] \varphi(x;z^2) \end{pmatrix},
\end{equation}
where $'$ denotes the derivative with respect to $x$; cf.\ \eqref{A226} and \eqref{A237}.
It follows from Lemma~\ref{A235} that
\begin{equation}\label{A255}
	\uptheta(x;-z) = \begin{pmatrix} p(x)\theta'(x;z^2) \\[2ex] -z\theta(x;z^2) \end{pmatrix},
	\qquad
	\upvarphi(x;-z) = \begin{pmatrix} \,-\dfrac{1}{z}p(x)\varphi'(x;z^2) \\[2ex] \varphi(x;z^2) \end{pmatrix}.
\end{equation}
Moreover, let $q_H$ be a singular Weyl coefficient of $H$
as in \eqref{A35}.

\begin{proof}[Proof of Theorem~\ref{A223}] 
	Let $\uptheta$ and $\upvarphi$ be as in \eqref{A254}.

	Item (i) follows directly from Theorem~\ref{A94}\,(i).
    Note that $a_+=a$ in Theorem~\ref{A94} since $w(x)>0$, $x\in(a,b)$ a.e.

	(ii)\,
	It follows from Theorem~\ref{A94}\,(ii) that, for $z\in\bb C\setminus\bb R$,
	\[
		\frac{\theta(x;z^2)}{\varphi(x;z^2)}
		= \frac{\uptheta_2(x;z)}{z\upvarphi_2(x;z)} \to \frac{1}{z}q_H(z)
		\qquad \text{as}\;\;x\nearrow b.
	\]
	The function $q_H$ is an odd function as the following calculation shows
	(where we use \eqref{A254} and \eqref{A255}):
	\begin{equation}\label{A265}
		q_H(-z) = \lim_{x\nearrow b}\frac{\uptheta_2(x;-z)}{\upvarphi_2(x;-z)}
		= \lim_{x\nearrow b}\frac{-\uptheta_2(x;z)}{\upvarphi_2(x;z)}
		= -q_H(z).
	\end{equation}
	Hence
	\begin{equation}\label{A253}
		\lim_{x\nearrow b}\frac{\theta(x;\lambda)}{\varphi(x;\lambda)}
		= \lim_{x\nearrow b}\frac{\uptheta(x;z)}{z\upvarphi(x;z)}
		= \frac{1}{z}q_H(z)
	\end{equation}
	is independent of the choice of $z$ such that $z^2=\lambda$, and therefore $m_{p,w}$
	is well defined as a function of $\lambda$ and analytic on $\bb C\setminus[0,\infty)$.
	It follows from \cite[Theorem~4.4]{langer.woracek:ninfrep}
	that $m_{p,w}\in\mc N_{<\infty}^{(\infty)}$.  Set $m_+(\lambda)\defeq \lambda m_{p,w}(\lambda)$;
	then $m_+\in\mc N_{<\infty}^{(\infty)}$ again by \cite[Theorem~4.4]{langer.woracek:ninfrep}.
	Now \cite[Proposition~4.8]{kaltenbaeck.winkler.woracek:nksym} implies
	that $\ind_- m_{p,w}+\ind_- m_+=\ind_- q_H$.  Moreover, $m_{p,w}$ and $m_+$ have
	infinity as only generalized pole of non-positive type and therefore,
	by the definition of $m_+$ and by \eqref{A7}, we have $\ind_-m_+-1\le\ind_-m_{p,w}\le \ind_-m_+$.
	This, together with the fact that $\ind_-q_H=\Delta(H)=\DeltaSL(p,w)$ by Theorem~\ref{A94}\,(ii)
	and \eqref{A247}, yields $\ind_-m_{p,w}=\bigl\lfloor\frac{\DeltaSL(p,w)}{2}\bigr\rfloor$.

	(iii)\,
	We can write
	\[
		\theta(\cdot\,;z^2)-m_{p,w}(z^2)\varphi(\cdot\,;z^2)
		= \frac{1}{z}\Bigl(\uptheta_2(\cdot\,;z)-q_H(z)\upvarphi_2(\cdot\,;z)\Bigr).
	\]
	By Theorem~\ref{A94}\,(iii) the right-hand side of this equality is in $L^2(w|_{(x_0,b)})$.
	Since (up to a scalar multiple) only one solution is in $L^2(w)$ at $b$
	(because of the limit point assumption at $b$), the value
	of $m_{p,w}(z^2)$ is uniquely determined by the $L^2$ property.

	(iv)\,
	The formula follows from item (iii) since any such $\psi$ is a multiple of
	$\psi_0\defeq\theta(\cdot\,;\lambda)-m_{p,w}(\lambda)\varphi(\cdot\,;\lambda)$
	and $\rbvSLr{\lambda}\psi_0=1$, $\rbvSLs{\lambda}\psi_0=-m_{p,w}(\lambda)$.

	(v)\,
	Let $\psi\in\NSL{\lambda}$ be such that $\psi|_{(x_0,b)}\in L^2(w|_{(x_0,b)})$.
	It follows from item (iv) and Theorem~\ref{A220}\,(iv) with the notation used there that
	\[
		\wh m_{p,w}(\lambda) = -\frac{\rbvSLspr{\lambda}\psi}{\rbvSLrpr{\lambda}\psi}
		= -\frac{\rbvSLs{\lambda}\psi+\pxxh(\lambda)\rbvSLr{\lambda}\psi}{\rbvSLr{\lambda}\psi}
		= m_{p,w}(\lambda)-\pxxh(\lambda).
	\]
	Since $\pxxh$ is the polynomial from Theorem~\ref{A220}\,(iv), it has the
	properties stated there.
\end{proof}

\medskip

\noindent
From \eqref{A253} we obtain the following relation between $m_{p,w}$ and $q_H$
if the same base point $x_0$ is chosen:
\begin{equation}\label{A240}
	m_{p,w}(z^2) = \frac{1}{z}q_H(z), \qquad z\in\bb C\setminus\bb R.
\end{equation}


Next we construct a measure with the help of the Stieltjes inversion formula
and the singular Titchmarsh--Weyl coefficient.  Before we formulate the theorem,
we introduce the following class of measures.

\begin{definition}\thlab{A217}
Let $\nu$ be a Borel measure on $\bb R$.  We say that $\nu\in\bb M^-$ if there exists
an $n\in\bb N_0$ such that
\begin{equation}\label{A218}
	\nu\bigl((-\infty,0)\bigr)=0 \qquad\text{and}\qquad
	\intop_{[0,\infty)}\frac{\rd\nu(t)}{(1+t)^{n+1}}<\infty.
\end{equation}
If $\nu\in\bb M^-$, we denote by $\Delta^-(\nu)$ the minimal $n\in\bb N_0$ such that
\eqref{A218} holds.
\end{definition}

The reason for the use of the minus sign in the notation will become
clearer in the next section; see, in particular, item 3 in \ref{A320}.

In the next theorem a measure is constructed, which will turn out to be
a spectral measure for the Sturm--Liouville equation \eqref{A208}.

\begin{theorem}[\textbf{The spectral measure}]\thlab{A227}
	\rule{0ex}{1ex} \\
	Let $(p;w)\in\KSL$ with $\dom(p;w)=(a,b)$ be given.
	Then there exists a unique Borel measure $\mu_{p,w}$ that satisfies
	\begin{equation}\label{A228}
		\mu_{p,w}\bigl([s_1,s_2]\bigr)=\frac{1}{\pi}\lim_{\eps\searrow 0}
		\lim_{\delta\searrow 0}\intop_{s_1-\eps}^{s_2+\eps}
		\Im m_{p,w}(t+i\delta)\,\rd t,\quad -\infty<s_1<s_2<\infty,
	\end{equation}
	where $m_{p,w}\in[m]_{p,w}$ is any singular Titchmarsh--Weyl coefficient associated
	with \eqref{A208}.
	We have $\mu_{p,w}\in\bb M^-$ and $\Delta^-(\mu_{p,w})=\DeltaSL(p,w)$.

	Moreover, $\mu_{p,w}(\{0\})>0$ if and only if
	\begin{equation}\label{A269}
		\int_a^b w(x)\rd x < \infty.
	\end{equation}
	If \eqref{A269} is satisfied, then
	\begin{equation}\label{A270}
		\mu_{p,w}\bigl(\{0\}\bigr) = \biggl[\,\int_a^b w(x)\rd x\biggr]^{-1}.
	\end{equation}
\end{theorem}


\bigskip

\noindent
We refer to the measure $\mu_{p,w}$ given by \eqref{A228} as the
\emph{spectral measure} associated with the Sturm--Liouville equation \eqref{A208}.
This is justified by Theorem~\ref{A229} below.

\begin{remark}\thlab{A312} \rule{0ex}{1ex}
	\begin{enumerate}[{\rm(i)}]
	\item
		The number $\Delta^-(\mu_{p,w})$, which describes the behaviour of the spectral
		measure $\mu_{p,w}$ at infinity, gives a finer measure of the growth
		of the coefficients $p$ and $w$ at the endpoint $a$ than $\ind_- m_{p,w}$,
		the negative index of the singular Titchmarsh--Weyl coefficient.
	\item
		If $\DeltaSL(p,w)=1$, then $\int_0^\infty \frac{\rd\mu_{p,w}(t)}{(1+t)^2}<\infty$,
		which is the classical growth condition for a spectral measure corresponding to
		a Nevanlinna function.  This is in accordance with the fact that in this case
		\eqref{A208} is in the limit circle case and therefore classical Hilbert space
		theory is sufficient to obtain a scalar spectral measure;
		see, e.g.\ \cite[Theorem~11.1]{bennewitz.everitt:2005}.
	\item
		If one uses the more general functions $\check{w}_l$ instead of $w_l$ as
		in Remark~\ref{A288}\,(iii), then --- similarly to \eqref{A308} --- the singular Titchmarsh--Weyl
		coefficient changes only by the real polynomial $p$ from \eqref{A307},
		which is of degree at most $\Delta-1$, and hence the spectral measure is unchanged.
	\end{enumerate}
\end{remark}

\medskip

\begin{proof}[Proof of Theorem~\ref{A227}]
	Let us first note that \eqref{A142}, \eqref{A95} and \eqref{A228}
	imply that
	\[
		\mu_H = \mu_{q_H} \qquad\text{and}\qquad
		\mu_{p,w} = \mu_{m_{p,w}}
	\]
	with the notation from \S\ref{A50}.
	Define $\tau(s)\defeq s^2$, $s\in\bb R$, and denote by $\mu^\tau$ the corresponding
	push-forward measure for a Borel measure $\mu$, i.e.\ $\mu^\tau(B)=\mu(\tau^{-1}(B))$
	for Borel sets $B\subseteq\bb R$.
	It follows from \cite[Theorem~4.4]{langer.woracek:ninfrep} and \eqref{A240} that
	\begin{equation}\label{A248}
		\mu_{p,w} \ll \mu_H^\tau \qquad\text{and}\qquad
		\frac{\rd\mu_{p,w}}{\rd\mu_H^\tau}(t) = \mathds{1}_{[0,\infty)}(t),
	\end{equation}
	where $\mathds{1}_{[0,\infty)}$ denotes the characteristic function of the
	interval $[0,\infty)$.

	For $m\in\bb N_0$, we obtain from \eqref{A248} that
	\begin{align*}
		\int_{\bb R} \frac{\rd\mu_H(t)}{(1+t^2)^{m+1}}
		&= \mu_H\bigl(\{0\}\bigr) + \int_{(0,\infty)}\frac{\rd\mu_H^\tau(s)}{(1+s)^{m+1}}
			\\[1ex]
		&= \mu_H\bigl(\{0\}\bigr) + \int_{(0,\infty)}\frac{\rd\mu_{p,w}(s)}{(1+s)^{m+1}}\,,
	\end{align*}
	which shows that $\mu_{p,w}\in\bb M^-$ and $\Delta^-(\mu_{p,w})=\Delta(\mu_H)=\Delta(H)=\DeltaSL(p,w)$,
	where the last equality follows from \eqref{A247}.
	The fact that a real polynomial makes no contribution in the Stieltjes inversion formula
	implies that the measure $\mu_{p,w}$ does not depend on the choice of the
	representative from $[m]_{p,w}$.

	Since $\mu_{p,w}(\{0\})=\mu_H(\{0\})$ by \eqref{A248}, the equivalence
	of $\mu_{p,w}(\{0\})>0$ and \eqref{A269} and the relation \eqref{A270}
	follow from Proposition~\ref{A77}.
\end{proof}

In the following theorem a Fourier transform is constructed, which yields the
unitary equivalence of the operator $A_{p,w}$ and the multiplication operator
in the space $L^2(\mu_{p,w})$.  In particular, this implies that the spectrum
of $A_{p,w}$ is simple.
Note that the function $\varphi$ that appears in the Fourier transform is the
(up to a scalar multiple) unique solution of \eqref{A208} that is `regular' at $a$;
cf.\ Remark~\ref{A271}\,(ii).

\begin{theorem}[\textbf{The Fourier transform}]\thlab{A229} \rule{0ex}{1ex}\\
	Let $(p;w)\in\KSL$ with $\dom(p;w)=(a,b)$ be given,
	and let $\mu_{p,w}$ be the spectral measure associated with \eqref{A208} via \eqref{A228}.
	Then the following statements hold.
	\begin{enumerate}[{\rm(i)}]
	\item
		The map defined by
		\begin{equation}\label{A230}
		\begin{aligned}
			(\Theta_{p,w}f)(t) \defeq \int_a^b \varphi(x;t)f(x)w(x)\,\rd x,
			\qquad t\in\bb R, \hspace*{15ex} &
			\\[-0.5ex]
			f\in L^2(w),\;\; \sup(\supp f)<b, &
		\end{aligned}
		\end{equation}
		extends to an isometric isomorphism from $L^2(w)$ onto $L^2(\mu_{p,w})$.
	\item
		The operator $\Theta_{p,w}$ establishes a unitary equivalence between $A_{p,w}$ and
		the operator $M_{\mu_{p,w}}$ of multiplication by the independent variable
		in $L^2(\mu_{p,w})$, i.e.\ we have
		\[
			\Theta_{p,w}A_{p,w} = M_{\mu_{p,w}}\Theta_{p,w}.
		\]
	\item
		For compactly supported functions, the inverse of $\Theta_{p,w}$ acts as an
		integral transformation, namely,
		\begin{equation}\label{A231}
		\begin{aligned}
			(\Theta_{p,w}^{-1}g)(x)=\int_0^\infty \varphi(x;t)g(t)\,\rd\mu_{p,w}(t),
			\qquad x\in(a,b), \hspace*{13ex} &
			\\[-0.5ex]
			g\in L^2(\mu_{p,w}),\;\; \supp g\text{ compact}. &
		\end{aligned}
		\end{equation}
	\end{enumerate}
\end{theorem}

\pagebreak[3]

\noindent
First we need a lemma.

\begin{lemma}\thlab{A212}
	Set
	\begin{align*}
		L^2_{\rm odd}(\mu_H) \defequ \bigl\{g\in L^2(\mu_H): g\text{ is odd}\bigr\},
			\\[1ex]
		L^2_{\rm even}(\mu_H) \defequ \bigl\{g\in L^2(\mu_H): g\text{ is even}\bigr\}.
	\end{align*}
	The Fourier transform $\Theta_H$ from Theorem~\ref{A4} maps
	\[
		\left\{\begin{pmatrix} f_1 \\ f_2 \end{pmatrix}\in L^2(H): f_2\equiv 0 \right\}
		\quad\text{bijectively onto }L^2_{\rm odd}(\mu_H)
	\]
	and
	\[
		\left\{\begin{pmatrix} f_1 \\ f_2 \end{pmatrix}\in L^2(H): f_1\equiv 0 \right\}
		\quad\text{bijectively onto }L^2_{\rm even}(\mu_H).
	\]
\end{lemma}

\begin{proof}
	Since $\upvarphi_1(x;-z) = -\upvarphi_1(x;z)$ and $\upvarphi_2(x;-z) = \upvarphi_2(x;z)$
	by \eqref{A255}, we have, for $f_1\in L^2\bigl(\frac{1}{p}\bigr)$ with $\sup(\supp f_1)<b$,
	that
	\begin{align*}
		\biggl[\Theta_H\binom{f_1}{0}\biggr](-t)
			&= \int_a^b \upvarphi(x;-t)^*H(x)\binom{f_1(x)}{0}\rd x
			= \int_a^b \upvarphi_1(x;-t)\frac{1}{p(x)}f_1(x)\rd x
			\\[1ex]
		&= -\int_a^b \upvarphi_1(x;t)\frac{1}{p(x)}f_1(x)\rd x
			= -\biggl[\Theta_H\binom{f_1}{0}\biggr](t),
	\end{align*}
	and similarly,
	\[
		\biggl[\Theta_H\binom{0}{f_2}\biggr](-t) = \biggl[\Theta_H\binom{0}{f_2}\biggr](t)
	\]
	for $f_2\in L^2(w)$.  Now the result follows from the bijectivity
	of $\Theta_H$.
\end{proof}


\begin{proof}[Proof of Theorem~\ref{A229}]\rule{0ex}{1ex}

	(i)\,
	The operator
	\[
		U:\begin{cases}
			L^2_{\rm even}(\mu_H) \to L^2(\mu_{p,w}) \\[2ex]
			f \mapsto g \quad\text{with } g(s) = f\bigl(\sqrt{s}\bigr), \;
			s\in[0,\infty),
		\end{cases}
	\]
	is well defined since $\mu_{p,w}((-\infty,0))=0$.  It is isometric because \eqref{A248}
	implies that, for $f\in L^2_{\rm even}(\mu_H)$,
	\begin{align*}
		\|Uf\|_{L^2(\mu_{p,w})}^2
		&= \int_{[0,\infty)} \bigl|f\bigl(\sqrt{s}\bigr)\bigr|^2 \rd\mu_{p,w}(s)
		= \int_{[0,\infty)} \bigl|f\bigl(\sqrt{s}\bigr)\bigr|^2 \rd\mu_H^\tau(s)
			\\[1ex]
		&= \int_{\bb R} |f(t)|^2 \rd\mu_H(t)
		= \|f\|_{L^2(\mu_H)}^2,
	\end{align*}
	where $\mu_H^\tau$ is as at the beginning of the proof of \thref{A227}.
	Moreover, $U$ is surjective and the inverse is given by $(U^{-1}g)(t)=g(t^2)$.
	
	Let $\iota_2$ and $P_2$ be as \eqref{A251}.
	By Lemma~\ref{A212} the operator $U\Theta_H\iota_2$ is well defined and isometric
	from $L^2(w)$ onto $L^2(\mu_{p,w})$.
	Moreover,
	\begin{align*}
		\bigl(U\Theta_H\iota_2f\bigr)(t)
		&= U\biggl[\int_a^b \upvarphi_2(x;\cdot)w(x)f(x)\rd x\biggr](t)
			\\[1ex]
		&= \int_a^b \upvarphi_2\bigl(x;\sqrt{t}\,\bigr)w(x)f(x)\rd x
			= \bigl(\Theta_{p,w}f\bigr)(t)
	\end{align*}
	by \eqref{A254}, which shows that
	\begin{equation}\label{A215}
		\Theta_{p,w} = U\Theta_H\iota_2.
	\end{equation}

	(ii)\,
	Let $M_{\mu_H}$ be the multiplication operator by the independent variable
	in $L^2(\mu_H)$ as in Theorem~\ref{A4}\,(ii).
	With Theorem~\ref{A4}\,(ii), the definition of $A_{p,w}$ and \eqref{A215} we obtain
	\begin{equation}\label{A214}
	\begin{aligned}
		\Theta_{p,w}A_{p,w}
		&= U\Theta_H\iota_2P_2 \bigl(T(H)\bigr)^2\iota_2
		= U\Theta_H \bigl(T(H)\bigr)^2\iota_2
			\\[1ex]
		&= UM_{\mu_H}^2\Theta_H\iota_2
		= UM_{\mu_H}^2U^{-1}\Theta_{p,w}
		= M_{\mu_{p,w}}\Theta_{p,w}.
	\end{aligned}
	\end{equation}

	(iii)\,
	It follows from \eqref{A215} that $\Theta_{p,w}^{-1}=P_2\Theta_H^{-1}U^{-1}$.
	If $g\in L^2(\mu_{p,w})$ with compact support, then, by \eqref{A254} and \eqref{A248},
	we have
	\begin{align*}
		(\Theta_{p,w}^{-1}g)(x)
		&= P_2\int_{\bb R} g(t^2)\upvarphi(x;t)\rd\mu_H(t)
		= \int_{\bb R} g(t^2)\varphi(x;t^2)\rd\mu_H(t)
			\\[1ex]
		&= \int_{[0,\infty)} g(s)\varphi(x;s)\rd\mu_H^\tau(s)
		= \int_{[0,\infty)} g(s)\varphi(x;s)\rd\mu_{p,w}(s),
	\end{align*}
	which shows the desired representation for $\Theta_{p,w}^{-1}$.
\end{proof}


\noindent
From Theorem~\ref{A6} we obtain the following uniqueness result, which says that equality of
singular Titchmarsh--Weyl coefficients or spectral measures implies
equality of the coefficients up to a reparameterization of the independent variable.
We cannot prove an existence result since we cannot characterize those spectral
measures that lead to diagonal Hamiltonians with non-vanishing determinant.
If we considered strings and the corresponding Krein--Feller operators, we would
also obtain an existence result: namely every measure from the class $\bb M^-$ is
the spectral measure of a certain string with two singular endpoints.

\begin{theorem}[\textbf{Global Uniqueness Theorem}]\thlab{A232} \rule{0ex}{1ex} \\
	Let $(p_1;w_1),(p_2;w_2)\in\KSL$ be given with $\dom(p_i;w_i)=(a_i,b_i)$, $i=1,2$.
	Then the following statements are equivalent:
	\begin{enumerate}[{\rm(i)}]
	\item
		there exists an increasing bijection $\gamma:(a_2,b_2)\to(a_1,b_1)$
		such that $\gamma$ and $\gamma^{-1}$ are locally absolutely continuous and
		\begin{equation}\label{A233}
			p_2(x) = \frac{1}{\gamma'(x)}p_1\bigl(\gamma(x)\bigr), \qquad
			w_2(x) = \gamma'(x)w_1\bigl(\gamma(x)\bigr)
		\end{equation}
		for $x\in(a_2,b_2)$ a.e.;
	\item
		$[m]_{p_1,w_1}=[m]_{p_2,w_2}$;
	\item
		$\mu_{p_1,w_1}=\mu_{p_2,w_2}$.
	\end{enumerate}
\end{theorem}

\begin{proof}
	The theorem follows from Theorem~\ref{A6} and Proposition~\ref{A97}
	if we recall \eqref{A72} and observe that
	\[
		q_{H_1}(z)-q_{H_2}(z)=z\bigl(m_{p_1,w_1}(z^2)-m_{p_2,w_2}(z^2)\bigr),
	\]
	which vanishes at $0$ and that hence $\alpha=0$ in Theorem~\ref{A6}.
\end{proof}

\begin{remark}\thlab{A298}
	If one considers Sturm--Liouville equations in impedance form, i.e.\ when $p=w$,
	and assumes that the left endpoints coincide,
	then equality of spectral measures is equivalent to the equality of
	coefficients a.e.\ because in this case one has $\gamma'=1$ if \eqref{A233}
	is satisfied.
	We refer to \cite{albeverio.hryniv.mykytyuk:2005} and the references therein
	for other types of inverse spectral theorems
	for Sturm--Liouville equations in impedance form.
\end{remark}

We also obtain a local version of the uniqueness result, which is an extension
of \cite[Theorem~1.5]{langer.woracek:lokinv} to the case of two singular endpoints.
See also \cite[\S4.4]{langer:2016} for a result on strings with regular left endpoint.
The next theorem follows immediately from Theorem~\ref{A31}.

\begin{theorem}[\textbf{Local Inverse Spectral Theorem}]\thlab{A249} \rule{0ex}{1ex} \\
	Let $(p_1;w_1),(p_2;w_2)\in\KSL$ be given with $\dom(p_i;w_i)=(a_i,b_i)$, $i=1,2$,
	and let $\tau>0$.  Moreover, for $i=1,2$, let $s_i(\tau)$ be the unique value $s_i$
	such that
	\[
		\intop_{a_i}^{s_i}\sqrt{\frac{w_i(\xi)}{p_i(\xi)}}\,\rd\xi = \tau
	\]
	if $\int_{a_i}^{b_i}\!\sqrt{\frac{w_i(\xi)}{p_i(\xi)}}\,\rd\xi > \tau$
	and set $s_i(\tau)\defeq b_i$ otherwise.
	Then the following statements are equivalent.
	\begin{enumerate}[{\rm(i)}]
	\item
		There exists an increasing bijection $\gamma:(a_2,s_2(\tau))\to(a_1,s_1(\tau))$
		such that $\gamma$ and $\gamma^{-1}$ are locally absolutely continuous and \eqref{A233}
		holds for  $x\in(a_2,s_2(\tau))$ a.e.
	\item
		There exist singular Titchmarsh--Weyl coefficients $m_{p_1,w_1}$ and $m_{p_2,w_2}$
		and there exists a $\beta\in(0,2\pi)$ such that, for each $\eps>0$,
		\[
			m_{p_1,w_1}\bigl(re^{i\beta}\bigr)-m_{p_2,w_2}\bigl(re^{i\beta}\bigr)
			= \rmO\Bigl(e^{(-2\tau+\eps)\sqrt{r}\sin\frac{\beta}{2}}\Bigr),
			\qquad r\to\infty.
		\]
	\item
		There exist singular Titchmarsh--Weyl coefficients $m_{p_1,w_1}$ and $m_{p_2,w_2}$
		and there exists a $k\ge0$ such that, for each $\delta\in(0,\pi)$,
		\begin{multline*}
			m_{p_1,w_1}(\lambda)-m_{p_2,w_2}(\lambda)
			= \rmO\Bigl(|\lambda|^k e^{-2\tau\Im\sqrt{\lambda}}\Bigr),
			\\[1ex]
			|\lambda|\to\infty,\;\;\lambda\in\bigl\{z\in\bb C:\delta\le\arg z\le2\pi-\delta\bigr\},
		\end{multline*}
		where $\sqrt{\lambda}$ is chosen so that $\Im\sqrt{\lambda}>0$.
	\end{enumerate}
\end{theorem}

%
%

\medskip

\noindent
In the next proposition we provide a sufficient condition for $(p;w)\in\KSL$.
This result is also used in Section~\ref{sec-Schroedinger} below.
Let us recall the following notation: we write $f(x)\asymp g(x)$ as $x\searrow0$
if there exist $c,C>0$ and $x_0>0$ such that $cg(x)\le f(x)\le Cg(x)$ for all $x\in(0,x_0)$.

\begin{proposition}\thlab{A242}
	Let $\alpha\ge1$ and assume that $p$ and $w$, defined on $(0,b)$ with $b>0$ or $b=\infty$
	satisfy \eqref{A210} with $a=0$ and
	\begin{equation}\label{A259}
		p(x)\asymp x^\alpha,\;\; w(x)\asymp x^\alpha \qquad\text{as}\;\; x\searrow0.
	\end{equation}
	If \eqref{A208} is in the limit point case at $b$, then $(p;w)\in\KSL$
	and $\DeltaSL(p,w)=\bigl\lfloor\frac{\alpha+1}{2}\bigr\rfloor$.
\end{proposition}

\noindent
Before we prove Proposition~\ref{A242}, we need a lemma.

\begin{lemma}\thlab{A243}
	Let $\alpha>1$ and assume that $p$ and $w$, defined on $(0,b)$ with $b>0$ or $b=\infty$
	satisfy \eqref{A259}. Moreover, choose $x_0\in(0,b)$.  Then
	\begin{equation}\label{A244}
		w_l(x) \asymp
		\begin{cases}
			x^l & \text{if\, $l$ is even and $l<\alpha+1$}, \\[1ex]
			x^{-\alpha+l} & \text{if\, $l$ is odd and $l<\alpha$},
		\end{cases}
	\end{equation}
	as $x\searrow0$.
\end{lemma}

\begin{proof}
	We prove the lemma by induction.  For $l=0$ the statement is clear from the definition of $w_0$.
	Now assume that \eqref{A244} is true for $l\in\bb N_0$.  If $l$ is even and $l+1<\alpha$, then
	\[
		w_{l+1}(x) \asymp \intop_x^{x_0} t^{-\alpha}t^l \rd t
		= \frac{1}{\alpha-l-1}\bigl(x^{-\alpha+l+1}-x_0^{-\alpha+l+1}\bigr)
		\asymp x^{-\alpha+l+1}.
	\]
	If $l$ is odd and $l+1<\alpha+1$, then
	\[
		w_{l+1}(x) \asymp \intop_0^x t^\alpha t^{-\alpha+l}\rd t
		= \frac{1}{l+1}x^{l+1} \asymp x^{l+1}.
	\]
	In both cases it follows that \eqref{A244} is true for $l+1$ instead of $l$.
	Hence the statement follows by induction.
\end{proof}

\begin{proof}[Proof of Proposition~\ref{A242}]
	The conditions (i) and (ii) in Definition~\ref{A213} are easy to check.
	For (iii) let us first consider the case $\alpha=1$.  Then $w_1(x) \asymp -\ln x$
	and hence $w_1\in L^2(w|_{(0,x_0)})$, which shows
	that $\DeltaSL(p,w)=1=\bigl\lfloor\frac{1+1}{2}\bigr\rfloor$.

	Now let $\alpha>1$.  If $l$ is even and $l<\alpha+1$, then
	\begin{equation}\label{A245}
		w_l \in L^2\Bigl(\frac{1}{p}\big|_{(0,x_0)}\Bigr)
		\quad\iff\quad \intop_0^{x_0} x^{-\alpha}x^{2l}\rd x < \infty
		\quad\iff\quad
		l>\frac{\alpha-1}{2}\,.
	\end{equation}
	If $l$ is odd and $l<\alpha$, then
	\begin{equation}\label{A246}
		w_l \in L^2(w|_{(0,x_0)})
		\quad\iff\quad \intop_0^{x_0} x^\alpha x^{2(-\alpha+l)}\rd x < \infty
		\quad\iff\quad l>\frac{\alpha-1}{2}\,.
	\end{equation}
	The minimal integer $l$ that satisfies $l>\frac{\alpha-1}{2}$
	is $\bigl\lfloor\frac{\alpha+1}{2}\bigr\rfloor$.  Since $\bigl\lfloor\frac{\alpha+1}{2}\bigr\rfloor<\alpha$
	for $\alpha>1$, the asymptotic relations \eqref{A244} are
	valid for $w_l$, $l\le\lfloor\frac{\alpha+1}{2}\rfloor$.
\end{proof}

\bigskip

\begin{example}\thlab{A299}
	Equations of the form
	\[
		-a_2y''-a_1y' = \lambda y,
	\]
	where $a_1$ and $a_2$ are continuous functions on $(a,b)$ and $a_2(x)>0$ for $x\in(a,b)$,
	can be written in the form \eqref{A208} with
	\[
		p(x) = \exp\biggl(\int_{x_0}^x\frac{a_1(t)}{a_2(t)}\rd t\biggr), \qquad
		w(x) = \frac{1}{a_2(x)}\exp\biggl(\int_{x_0}^x\frac{a_1(t)}{a_2(t)}\rd t\biggr)
	\]
	with some $x_0\in[a,b]$.
	As an example we consider the associated Laguerre equation
	\[
		-xy''(x)-(1+\alpha-x)y'(x) = \lambda y(x), \qquad x\in(0,\infty),
	\]
	with $\alpha\ge0$.  For $p$ and $w$ one obtains
	\[
		p(x) = x^{\alpha+1}e^{-x}, \qquad w(x) = x^\alpha e^{-x}.
	\]
	It can be shown in a similar way as in Proposition~\ref{A242} that $(p;w)\in\KSL$
	with $\Delta=\lfloor\alpha+1\rfloor$.  Hence the singular Titchmarsh--Weyl
	coefficient belongs to $\mc N_\kappa^{(\infty)}$ with
	$\kappa=\bigl\lfloor\frac{\alpha+1}{2}\bigr\rfloor$.
	This is in agreement with \cite{dijksma.shondin:2002} where a model for this
	singular Titchmarsh--Weyl coefficient was constructed.
	For $\alpha<-1$ the associated Laguerre equation was studied with the help of
	Pontryagin spaces in \cite{krall:1979, krall:1982, derkach:1998};
	in this case the results of the next subsection can be applied.
\end{example}

%
%
\section{Sturm--Liouville equations without \\ potential: singular {\boldmath$w$}}
\label{sec-SLw}
%
%

In this section we consider the case when $w$ is not integrable at $a$
but $\frac{1}{p}$ is.  In Definition~\ref{A213} and most theorems one just has to
swap the roles of $\frac{1}{p}$ and $w$.
Let us state the definition of the class of coefficients explicitly.

\begin{definition}\thlab{A318}
	We say that $(p;w)\in\KSLpl$ if $p$ and $w$ are defined on some interval $(a,b)$
	and they satisfy \eqref{A210} and the following conditions.
	\begin{enumerate}[{\rm(i)}]
	\item
		For one (and hence for all) $x_0\in(a,b)$,
		\[
			\intop_a^{x_0} w(x)\rd x = \infty
			\qquad\text{and}\qquad
			\intop_a^{x_0} \frac{1}{p(x)}\rd x < \infty.
		\]
	\item
		For one (and hence for all) $x_0\in(a,b)$,
		\[
			\intop_a^{x_0} \intop_a^x \frac{1}{p(t)}\rd t\,w(x)\rd x < \infty.
		\]
	\item
		Let $x_0\in(a,b)$ and define functions $w_l$, $l=0,1,\dots$, recursively by
		\begin{align*}
			w_0(x) &= 1,
				\\[1ex]
			w_{l}(x) &= \begin{cases}
				\displaystyle \,\intop_x^{x_0} w(t)w_{l-1}(t)\rd t \quad
				& \text{if $l$ is odd},
				\\[4ex]
				\displaystyle \,\intop_a^x \frac{1}{p(t)}w_{l-1}(t)\rd t & \text{if $l$ is even}.
			\end{cases}
		\end{align*}
		There exists an $n\in\bb N_0$ such that
		\begin{equation}\label{A319}
			w_n\big|_{(a,x_0)} \in \begin{cases}
				L^2\bigl(w\big|_{(a,x_0)}\bigr) & \text{if $n$ is even},
				\\[2ex]
				L^2\bigl(\frac{1}{p}\big|_{(a,x_0)}\bigr) & \text{if $n$ is odd}.
			\end{cases}
		\end{equation}
	\item
		Equation \eqref{A208} is in the limit point case at $b$, i.e.\
		for $\lambda\in\bb C\setminus\bb R$, equation \eqref{A208} has
		(up to a scalar multiple) only one solution in $L^2(w|_{(x_0,b)})$
		for $x_0\in(a,b)$.
	\end{enumerate}
	If $(p;w)\in\KSLpl$, we denote by $\DeltaSLpl(p,w)$ the minimal $n\in\bb N_0$ such that
	\eqref{A319} holds.
\end{definition}

\medskip

\begin{nremark}{Differences between the classes $\KSL$ and $\KSLpl$}\label{A320}
In the following list we mention the major differences that occur in theorems
and other important statements corresponding to coefficients in $\KSL$ and $\KSLpl$,
respectively.
\begin{enumerate}[1.]
	\item
		The analogue of Remark~\ref{A207}\,(ii) is not true;
		the equation \eqref{A208} is always in the limit point case at $a$
		if $\DeltaSLpl(p,w)\ge1$.
		This follows from the fact that $1\notin L^2(w|_{(a,x_0)})$.
	\item
		The regularized boundary values have the following form (with $\Delta=\DeltaSLpl(p,w)$):
		\begin{align*}
			\rbvSLplr{\lambda}\psi &= \lim_{x\searrow a}\psi(x),
				\\[1ex]
			\rbvSLpls{\lambda}\psi &=  \lim_{x\searrow a}\Biggl[p(x)\psi'(x)
				+\,\sum_{k=1}^{\lfloor\frac{\Delta}{2}\rfloor}
				\lambda^k\Bigl(w_{2k}(x)p(x)\psi'(x)-w_{2k-1}(x)\psi(x)\Bigr)
				\\[1ex]
			&\quad -
				\left\{\begin{array}{ll}
					\lambda^{\frac{\Delta+1}{2}} w_\Delta(x)\psi(x) & \text{if $\Delta$ is odd}  \\[1ex]
					0 & \text{if $\Delta$ is even}
				\end{array}\right\}
				\\[1ex]
			&\quad + \bigl(\lim_{t\searrow a}\psi(t)\bigr)\biggl(\,
				\sum_{k=\lfloor\frac{\Delta+3}{2}\rfloor}^\Delta
				\sum_{l=0}^{2k-\Delta-2}
				(-1)^{l+1} \lambda^k w_l(x)w_{2k-l-1}(x)\biggr)\Biggr].
		\end{align*}
	\item
		The singular Titchmarsh--Weyl coefficient $m_{p,w}^+$, which is defined as in \eqref{A225},
		is connected with the singular Weyl coefficient of the corresponding canonical system via
		\begin{equation}\label{A283}
			m_{p,w}^+(z^2) = zq_H(z).
		\end{equation}
		This relation explains the use of the notation with $+$ as this was used,
		e.g.\ in \cite{kaltenbaeck.winkler.woracek:nksym} and \cite{langer.woracek:ninfrep}.
		The singular Titchmarsh--Weyl coefficient $m_{p,w}^+$ belongs
		to $\mc N_\kappa^{(\infty)}$
		where $\kappa=\bigl\lfloor\frac{\DeltaSLpl(p,w)+1}{2}\bigr\rfloor$.
		The equivalence classes $[m]_{p,w}^+$ are defined not with respect to
		the equivalence relation $\hat{\sim}$ defined in \eqref{A287} but with the
		equivalence relation $\sim$ defined in \eqref{A98}, which is
		\[
			m_1 \sim m_2\quad \defequiv\quad m_1-m_2\in\bb R[z],\; (m_1-m_2)(0)=0.
		\]
	\item
		The spectral measure $\mu_{p,w}^+$ belongs to the class $\bb M^+$, which is the set of
		Borel measures on $\bb R$ such that
		\begin{equation}\label{A281}
			\nu\bigl((-\infty,0]\bigr)=0 \qquad\text{and}\qquad
			\int_{(0,\infty)}\frac{\rd\nu(t)}{t(1+t)^{n+1}}<\infty.
		\end{equation}
		If $\nu\in\bb M^+$, we denote by $\Delta^+(\nu)$ the minimal $n\in\bb N_0$
		such that \eqref{A281} holds.
		Then $\Delta^+(\mu_{p,w}^+)=\DeltaSLpl(p,w)$.
	\item
		Instead of \eqref{A270} one has
		\begin{equation}\label{A282}
			-\lim_{\lambda\nearrow0} m_{p,w}^+(\lambda)
			= \biggl[\,\int_a^b \frac{1}{p(x)}\,\rd x\biggr]^{-1},
		\end{equation}
		and that the left-hand side is zero if and only
		if the integral on the right-hand side is infinite.
		Relation \eqref{A282} follows from \eqref{A79}, \eqref{A130}
		and \eqref{A283}.
\end{enumerate}
\end{nremark}

\medskip

\noindent
The global uniqueness result is different from the one in the previous subsection
since adding a constant to the singular Titchmarsh--Weyl coefficient corresponds
to a more complicated transformation;
cf.\ also \cite[Corollary~3.6]{eckhardt.gesztesy.nichols.teschl:2013}
for the case when the equations are in impedance form.

\pagebreak[3]
\begin{theorem}[\textbf{Global Uniqueness Theorem}]\thlab{A284} \rule{0ex}{1ex}
	\begin{enumerate}[{\rm(i)}]
	\item
		Let $(p_1;w_1),(p_2;w_2)\in\KSLpl$ be given with $\dom(p_i;w_i)=(a_i,b_i)$, $i=1,2$.
		Assume that there exist singular Titchmarsh--Weyl coefficients $m_{p_i,w_i}^+$ corresponding
		to $(p_i;w_i)$ for $i=1,2$ such that
		\begin{equation}\label{A257}
			m_{p_1,w_1}^+(\lambda)-m_{p_2,w_2}^+(\lambda)
			= c_n\lambda^n+\ldots+c_1\lambda+c_0
		\end{equation}
		with $c_0,\dots,c_n\in\bb R$.  Then there exists an increasing bijection
		$\gamma:(a_2,b_2)\to(a_1,b_1)$ such that $\gamma$ and $\gamma^{-1}$ are locally absolutely
		continuous and
		\begin{equation}\label{A285}
		\begin{aligned}
			p_2(x) &= \frac{1}{\gamma'(x)}\biggl(1+c_0\int_{a_1}^{\gamma(x)} \frac{1}{p_1(t)}\,\rd t\biggr)^2
			p_1\bigl(\gamma(x)\bigr),
				\\[2ex]
			w_2(x) &= \gamma'(x)\biggl(1+c_0\int_{a_1}^{\gamma(x)} \frac{1}{p_1(t)}\,\rd t\biggr)^2
			w_1\bigl(\gamma(x)\bigr)
		\end{aligned}
		\end{equation}
		for $x\in(a_2,b_2)$;
		for all $x\in(a_2,b_2)$ one has
		\begin{equation}\label{A297}
			1+c_0\int_{a_1}^{\gamma(x)}\frac{1}{p_1(t)}\rd t>0.
		\end{equation}
		Moreover, $\DeltaSLpl(p_1,w_1)=\DeltaSLpl(p_2,w_2)$.
	\item
		Let $(p_1;w_1)\in\KSLpl$ be given with $\dom(p_1;w_1)=(a_1,b_1)$.
		Let $(a_2,b_2)\subseteq\bb R$ be an open interval, $\gamma:(a_2,b_2)\to(a_1,b_1)$
		an increasing bijection such that $\gamma$ and $\gamma^{-1}$ are locally absolutely
		continuous, and let $c_0\in\bb R$ such that
		\begin{equation}\label{A252}
			1+c_0\int_{a_1}^{b_1} \frac{1}{p_1(t)}\,\rd t \ge 0.
		\end{equation}
		Define functions $p_2,w_2$ by \eqref{A285}.  Then $(p_2;w_2)\in\KSLpl$
		with $\DeltaSLpl(p_1,w_1)=\DeltaSLpl(p_2,w_2)$, and there exist singular Titchmarsh--Weyl
		coefficients $m_{p_i,w_i}^+$, $i=1,2$, such that
		\[
			m_{p_1,w_1}^+(\lambda)-m_{p_2,w_2}^+(\lambda) = c_0.
		\]
	\item
		Let $(p_1;w_1),(p_2;w_2)\in\KSLpl$ be given with $\dom(p_i;w_i)=(a_i,b_i)$, $i=1,2$.
		Then $\mu_{p_1,w_1}^+=\mu_{p_2,w_2}^+$ if and only if there exists $\gamma$ as above and $c_0\in\bb R$
		such that \eqref{A285} and \eqref{A252} hold.
	\end{enumerate}
\end{theorem}

\medskip

\noindent
Before we can prove the theorem we need a lemma about a transformation of diagonal
Hamiltonians.

\begin{lemma}\thlab{A260}
	Let $H\in\bb H$ be a diagonal Hamiltonian with $\dom H=(a,b)$ of the form
	\[
		H(x) = \begin{pmatrix} h_{11}(x) & 0 \\[1ex] 0 & h_{22}(x) \end{pmatrix}.
	\]
	Assume that $h_{22}(x)>0$ for almost all $x\in(a,b)$ and let $q_H$ be a singular
	Weyl coefficient and $\mu_H$ the corresponding spectral measure.
	Moreover, let $c\in\bb R$ and define the functions
	\begin{align}
		\alpha(x) \defequ 1+c\int_a^x h_{22}(t)\rd t, \qquad x\in(a,b),
		\label{A262}
		\\[1ex]
		\widetilde q(z) \defequ q_H(z)-\frac{c}{z}\,.
		\label{A263}
	\end{align}
	Then the following statements are equivalent:
	\begin{enumerate}[{\rm(i)}]
	\item
		$\alpha(x)>0$ for all $x\in(a,b)$;
	\item
		$c+\mu_H\bigl(\{0\}\bigr) \ge 0$;
	\item
		$\widetilde q \in \mc N_{<\infty}^{(\infty)}$.
	\end{enumerate}
	If these conditions are satisfied, then
	\begin{equation}\label{A261}
		\widetilde H(x) \defeq
		\begin{pmatrix}
			\bigl(\alpha(x)\bigr)^2h_{11}(x) & 0 \\[2ex]
			0 & \dfrac{h_{22}(x)}{\bigl(\alpha(x)\bigr)^2}
		\end{pmatrix},
		\qquad x\in(a,b),
	\end{equation}
	belongs to $\bb H$ and $\widetilde q$ is a singular Weyl coefficient for $\widetilde H$.
\end{lemma}

\begin{proof}
	First note that the integral in \eqref{A262} is positive for all $x\in(a,b)$ and
	strictly increasing in $x$.  Hence (i) is equivalent to
	\[
		c+\biggl[\int_a^b h_{22}(t)\rd t\biggr]^{-1} \ge 0,
	\]
	where we use $1/\infty=0$ in the case when the integral is infinite.
	Since $\mu_H(\{0\})=\bigl[\int_a^b h_{22}(t)\rd t\bigr]^{-1}$ by Proposition~\ref{A77},
	the equivalence of (i) and (ii) follows.
	The function $\widetilde q$ is in $\mc N_{<\infty}$.  The only possible
	finite generalized pole of non-positive type of $\wt q$ is $0$.
	Hence $\widetilde q\in\mc N_{<\infty}^{(\infty)}$ if and only if
	\[
		\lim_{\varepsilon\searrow0} i\varepsilon \widetilde q(i\varepsilon) \le 0;
	\]
	see, e.g.\ \cite[Theorem~3.1]{langer:1986}.
	It follows from \cite[Theorem~3.9\,(ii)]{langer.woracek:ninfrep} that
	\[
		c+\mu_H\bigl(\{0\}\bigr)
		= c-\lim_{\varepsilon\searrow0} i\varepsilon q_H(i\varepsilon)
		= -\lim_{\varepsilon\searrow0} i\varepsilon \widetilde q(i\varepsilon),
	\]
	which implies the equivalence of (ii) and (iii).

	For the rest of the proof assume that $\alpha(x)>0$ for all $x\in(a,b)$.
	One can easily show that the assertions of the lemma are unaffected by reparameterizations.
	So we can assume that $H$ is defined on $(0,\infty)$.
	Let $\mf h$ be the indefinite Hamiltonian associated with $H$ as
	in \S\ref{A108} with Weyl coefficient $q_{\mf h}$
	such that $q_{\mf h}=q_H$, and let $\omega_{\mf h}$ be the corresponding maximal chain
	of matrices as in \S\ref{A52}.
	It follows from Lemma~\ref{A80} that
	\[
		\widetilde \alpha(x) \defeq 1-c\frac{\partial}{\partial z}\omega_{\mf h,21}(x;z)\Big|_{z=0}
		= \begin{cases}
			1, & x\in[-1,0), \\[2ex]
			\alpha(x), & x\in(0,\infty).
		\end{cases}
	\]
	The transformation $\mf T_c$ from \cite[Definition~4.1]{kaltenbaeck.woracek:p3db}
	applied to $\omega_{\mf h}$ yields a maximal chain of matrices $\widetilde W$, where
	\[
		\widetilde W(x;z) \defeq \bigl(\mf T_c\omega_{\mf h}\bigr)(x;z)
		= \begin{pmatrix} 1 & -\dfrac{c}{z} \\[2ex] 0 & 1 \end{pmatrix}\omega_{\mf h}(x;z)
		\begin{pmatrix} \dfrac{1}{\widetilde\alpha(x)} & \dfrac{c}{z} \\[3ex] 0 & \widetilde\alpha(x) \end{pmatrix},
	\]
	so that $\widetilde W$ corresponds to an indefinite Hamiltonian $\widetilde{\mf h}$
	whose Weyl coefficient is
	\[
		q_{\widetilde{\mf h}}(z) = q_{\mf h}(z)-\frac{c}{z}\,;
	\]
	see \cite[Theorem~4.4]{kaltenbaeck.woracek:p3db}.
	Differentiating $\widetilde W$ with respect to $x$ and using the
	differential equation \eqref{A70}, i.e.\
	\[
		\frac{\partial}{\partial x}\widetilde W(x;z)J
		= z\widetilde W(x;z)\widetilde H(x),
	\]
	one can easily show that the Hamiltonian function $\widetilde H$ that
	corresponds to $\widetilde{\mf h}$ is given by \eqref{A261};
	cf.\ \cite[Rule~4]{winkler:1995a}.
\end{proof}

\begin{proof}[Proof of Theorem~\ref{A284}]
	Throughout the proof let $H_1$ and $H_2$ be the corresponding Hamiltonians
	\[
		H_i(x) = \begin{pmatrix} w_i(x) & 0 \\[1ex] 0 & \dfrac{1}{p_i(x)} \end{pmatrix},
		\qquad i=1,2,
	\]
	and let $q_{H_1}$, $q_{H_2}$ be the singular
	Weyl coefficients such that $m_{p_i,w_i}^+(z^2)=zq_{H_i}(z)$, $i=1,2$ as in \eqref{A283}.

	(i)\,
	Suppose that \eqref{A257} holds.  Define $\widetilde q$ as in \eqref{A263}
	with $H=H_1$ and $c=c_0$, i.e.\
	\begin{equation}\label{A258}
		\widetilde q(z) = q_{H_1}(z)-\frac{c_0}{z}
		= q_{H_2}(z)+c_nz^{2n-1}+c_{n-1}z^{2n-3}+\ldots+c_1z.
	\end{equation}
	The equality of the first and the last expression in \eqref{A258} implies that
	$\widetilde q \in \mc N_{<\infty}^{(\infty)}$, i.e.\ condition (iii)
	in Lemma~\ref{A260} is satisfied.  Hence we can apply Lemma~\ref{A260}, which
	yields a Hamiltonian $\widetilde H$ with corresponding singular Weyl coefficient $\widetilde q$.
	Since $\widetilde q$ and $q_{H_2}$ differ only by a real polynomial without constant term,
	Theorem~\ref{A5}\,(iii) implies that $H_2$ is a reparameterization of $\widetilde H$.
	This, together with \eqref{A261} and \eqref{A262}, shows \eqref{A285}.
	Relation \eqref{A297} follows from Lemma~\ref{A260}\,(i).
	Since $\mu_{p_1,w_1}^+=\mu_{p_2,w_2}^+$, we obtain
	\[
		\DeltaSLpl(p_1,w_1) = \Delta^+(\mu_{p_1,w_1}^+)
		= \Delta^+(\mu_{p_2,w_2}^+) = \DeltaSLpl(p_2,w_2).
	\]

	(ii)\,
	Condition \eqref{A252} implies that (i) in Lemma~\ref{A260} is satisfied with $H=H_1$
	and $c=c_0$.  Hence we can apply Lemma~\ref{A260}, which yields $\widetilde H$.
	The assertion follows since $H_2$ is a reparameterization of $\widetilde H$.

	(iii)\,
	It follows from \cite[Theorem~3.9\,(iv)]{langer.woracek:ninfrep} that
	$\mu_{p_1,w_1}^+=\mu_{p_2,w_2}^+$ if and only if $m_{p_1,w_1}^+$ and $m_{p_2,w_2}^+$ differ
	by a real polynomial.  Now the claim follows from (i) and (ii).
\end{proof}

Let us finally point out that Proposition~\ref{A242} remains valid for the situation
in this section if \eqref{A259} is replaced by
\[
	p(x)\asymp x^{-\alpha},\;\; w(x)\asymp x^{-\alpha} \qquad\text{as}\;\; x\searrow0
\]
and $\KSL$ is replaced by $\KSLpl$.

%
%
\section{Schr\"odinger equations}
\label{sec-Schroedinger}
%
%

Let $V\in L^1_{\rm loc}(0,b)$ with $b>0$ or $b=\infty$ and consider
the one-dimensional Schr\"odinger equation
\begin{equation}\label{A200}
	-u''(x)+V(x)u(x) = \lambda u(x).
\end{equation}
In this section the left endpoint needs to be finite, which without loss
of generality we assume to be $0$.
In the following we write $\dom(V)\defeq(0,b)$.

Assume that, for $\lambda=0$, equation \eqref{A200} has a solution $\phi$
(i.e.\ $V=\frac{\phi''}{\phi}$) in $W^{2,1}_{\rm loc}(0,b)$ that satisfies
\begin{equation}\label{A201}
	\begin{aligned}
		& \phi(x)>0 \;\;\text{for all}\;\;x\in(0,b),
		\\[1ex]
		& \phi\big|_{(0,x_0)}\in L^2(0,x_0),\;\frac{1}{\phi}\Big|_{(0,x_0)}\notin L^2(0,x_0)
		\quad\text{for some}\;\;x_0\in(0,b).
	\end{aligned}
\end{equation}

\noindent
A similar approach, namely to assume the existence of a particular solution
instead of explicit conditions on the coefficients, was used in \cite{davies:2013}.

Note that $\phi$ with the above properties is determined only up to a multiplicative
positive constant; see Remark~\ref{A204} below for a further discussion
of this non-uniqueness.

Let $x_0\in(0,b)$ and define functions $\wt w_l$, $l=0,1,\ldots,$ recursively by
\begin{equation}\label{A311}
\begin{aligned}
	\wt w_0(x) &= \frac{1}{\phi(x)}\,,
	\\[1ex]
	\wt w_k(x) &= \begin{cases}
		\displaystyle \phi(x)\intop_x^{x_0} \frac{1}{\phi(t)}\wt w_{k-1}(t)\rd t \quad
		& \text{if $k$ is odd}, \\[4ex]
		\displaystyle \frac{1}{\phi(x)}\intop_0^x \phi(t)\wt w_{k-1}(t)\rd t
		& \text{if $k$ is even}.
	\end{cases}
\end{aligned}
\end{equation}

\begin{remark}\thlab{A302}
	It follows from the last condition in \eqref{A201}
	and \cite[Theorem~2.2]{niessen.zettl:1992} that $\phi$ is a principal solution
	of \eqref{A200} with $\lambda=0$.  The function $\wt w_1$ is a non-principal
	solution of \eqref{A200} with $\lambda=0$, and one has
	\begin{equation}\label{A294}
		-\wt w_{k+2}''+V\wt w_{k+2} = \wt w_k \qquad\text{when $k\in\bb N$ is odd}.
	\end{equation}
	In \cite{kurasov.luger:2011} a sequence of functions $g_k$ was used which satisfy
	the relations
	\[
		-g_{k+1}''+Vg_{k+1} -\mu_kg_{k+1} = g_k
	\]
	with pairwise distinct numbers $\mu_k$.
\end{remark}

\begin{remark}\thlab{A310}
	Instead of $\wt w_k$ one can use more general functions $\check{w}_k$ that
	are defined as in \eqref{A311} but with the relation
	\[
		\check{w}_k(x) = \phi(x)\left[\,\intop_x^{x_0}\frac{1}{\phi(t)}
		\check{w}_{k-1}(t)\rd t+c_k\right], \qquad
		k\;\text{odd},
	\]
	with arbitrary constants $c_k\in\bb R$;
	cf.\ Remarks~\ref{A288}\,(iii) and \ref{A312}\,(iii).
	The spectral measure that is constructed below remains the same.
\end{remark}

\medskip

\noindent
In this section we consider the following class of potentials.

\begin{definition}\thlab{A275}
	We say that $V\in\KSchr$ if $V\in L^1_{\rm loc}(0,b)$, there exists a $\phi$
	so that $\phi$ is a solution of \eqref{A200} with $\lambda=0$, that \eqref{A201}
	holds and that the following conditions are satisfied.
	\begin{enumerate}[{\rm(i)}]
	\item
		For one (and hence for all) $x_0\in(0,b)$,
		\[
			\intop_0^{x_0} \phi(x)\wt w_1(x)\rd x < \infty.
		\]
	\item
		There exists an $n\in\bb N$ such that
		\begin{equation}\label{A277}
			\wt w_n\big|_{(0,x_0)} \in L^2(0,x_0).
		\end{equation}
	\item
		Equation \eqref{A200} is in the limit point case at $b$.
	\end{enumerate}
	If $V\in\KSchr$, we denote by $\DeltaSchr(V)$ the minimal $n\in\bb N$
	such that \eqref{A277} holds.
\end{definition}


\begin{remark}\thlab{A293}\rule{0ex}{1ex}
	\begin{enumerate}[{\rm(i)}]
	\item
		Since the function $\wt w_1$ is a non-principal solution of \eqref{A200}
		with $\lambda=0$, we have $\Delta(V)=1$ if and only if
		equation \eqref{A200} is regular or in the limit circle case at the left endpoint.
	\item
		One can also consider the case when \eqref{A200} is regular or in the limit circle
		case at the right endpoint $b$.  In this case one has to impose a fixed self-adjoint
		boundary condition at $b$; cf.\ Remark~\ref{A274}.
		This situation was considered, e.g.\ in \cite{silva.toloza:2014} in connection
		with $n$-entire operators.
	\end{enumerate}
\end{remark}

\medskip

\noindent
In order to apply the results from Section~\ref{sec-SL}, we set
\begin{equation}\label{A276}
	p(x)=w(x)\defeq \bigl(\phi(x)\bigr)^2, \qquad x\in(0,b).
\end{equation}
With $w_k$ defined as in \eqref{A206} we have
\begin{equation}\label{A329}
	\wt w_k(x) =
	\begin{cases}
		\dfrac{w_k(x)}{\phi(x)}  & \text{if $k$ is even}, \\[3ex]
		\phi(x)w_k(x) & \text{if $k$ is odd}.
	\end{cases}
\end{equation}
It is easy to see that $V\in\KSchr$ if and only if $(p;w)\in\KSL$
with $p,w$ from \eqref{A276}; for the equivalence of the limit point property
at $b$ see \eqref{A315} below.
Moreover, if $V\in\KSchr$, then $\DeltaSchr(V)=\DeltaSL(p,w)$.

\medskip

\begin{example}\thlab{A266}\rule{0ex}{1ex}
\begin{enumerate}[{\rm(i)}]
	\item
		A large subclass of $\KSchr$ is the following.
		Let $b>0$ or $b=\infty$, let $V_0\in L^1_{\rm loc}(0,b)$
		and let $l\in\bigl[-\frac12,\infty\bigr)$.  Moreover, set
		\begin{equation}\label{A205}
			V(x) = \frac{l(l+1)}{x^2}+V_0(x)
		\end{equation}
		and assume that
		\begin{equation}\label{A303}
		\begin{alignedat}{2}
			& xV_0(x)\big|_{(0,x_0)} \in L^1(0,x_0) \qquad & & \text{if $l>-\frac12$}\,,
			\\[1ex]
			& (\ln x)xV_0(x)\big|_{(0,x_0)} \in L^1(0,x_0) \qquad & & \text{if $l=-\frac12$}
		\end{alignedat}
		\end{equation}
		with some $x_0\in(0,b)$.
		Moreover, suppose that the minimal operator associated with \eqref{A201} is
		bounded from below and that \eqref{A200} is in the limit point case at $b$.
		Under the assumption that \eqref{A303} is valid, it follows
		from \cite[Lemma~3.2]{kostenko.teschl:2011} that there exists a solution $\phi$
		of \eqref{A200} with $\lambda=0$ such that
		\begin{equation}\label{A291}
			\phi(x) = x^{l+1}\bigl(1+\rmo(x)\bigr), \qquad x\searrow 0.
		\end{equation}
		Assume that $\phi(x)>0$ for $x\in(0,b)$, which is satisfied, e.g.\ if
		the minimal operator is uniformly positive, which can be achieved by
		a shift of the spectral parameter.
		Now it follows from Proposition~\ref{A242} that $(p;w)\in\KSL$, and hence $V\in\KSchr$ and
		$\DeltaSchr(V)=\DeltaSL(p,w)=\bigl\lfloor l+\frac{3}{2}\bigr\rfloor$.

		Since $l=0$ is allowed in \eqref{A205}, the class $\KSchr$ contains potentials where
		$0$ is a regular endpoint.
		If $l\in\bigl[-\frac12,\frac12\bigr)\setminus\{0\}$, then \eqref{A200} is
		in the limit circle case at $0$ and $\DeltaSchr(V)=1$.

		Potentials of the form \eqref{A205} have been studied in many papers; see, e.g.\
		\cite{fulton:2008, fulton.langer:2010,
		hryniv.sacks:2010, kostenko.sakhnovich.teschl:2010,
		kostenko.teschl:2011, kurasov.luger:2011,
		fulton.langer.luger:2012, albeverio.hryniv.mykytyuk:2012,
		kostenko.sakhnovich.teschl:2012a, kostenko.teschl:2013,
		silva.toloza:2014, silva.teschl.toloza:2015, eckhardt:2014, luger.neuner:2015}.
	\item
		The class $\KSchr$ contains also potentials that have a stronger singularity
		at the left endpoint than those considered in (i).
		If
		\begin{equation}\label{A328}
			V(x)=\frac{\phi''(x)}{\phi(x)} \qquad\text{where}\qquad
			\phi(x)\asymp x^\beta, \quad x\searrow0
		\end{equation}
		with $\beta\ge\frac{1}{2}$, $\phi(x)>0$ for $x\in(0,b)$ and \eqref{A200}
		is in the limit point case at $b$,
		then $V\in\KSchr$ with $\DeltaSchr(V)=\bigl\lfloor\beta+\frac{1}{2}\bigr\rfloor$;
		cf.\ Proposition~\ref{A242}.
		For instance, functions of the form
		\[
			\phi(x)=x^\beta\Bigl(2+\sin\frac{1}{x}\Bigr), \qquad x\in(0,x_0),
		\]
		with $\beta>0$ lead to oscillatory potentials that do not satisfy \eqref{A303},
		namely,
		\[
			V(x) = -\frac{1}{x^4}\cdot\frac{\sin\frac{1}{x}}{2+\sin\frac{1}{x}}
			+\rmO\Bigl(\frac{1}{x^3}\Bigr), \qquad x\searrow0.
		\]
		It follows from Lemma~\ref{A243} and \eqref{A329} that if $V$ is of
		the form in \eqref{A328}, then
		\begin{equation}\label{A330}
			\wt w_k(x) \asymp x^{-\beta+k}, \qquad x\searrow0,\;
			k\in\bb N_0,\; k<2\beta.
		\end{equation}
		In particular, the relation in \eqref{A330} is valid for
		$k=0,1,\ldots,2\Delta-1$ if $\beta$ is not an odd integer multiple
		of $\frac{1}{2}$, and it is valid for $k=0,1,\ldots,2\Delta-2$, otherwise.
	\item
		The function $V(x)=\frac{1}{x^4}$ does not belong to $\KSchr$.
		It can easily be checked that the only possible choice for $\phi$
		(up to scalar multiples) is $\phi(x)=xe^{-1/x}$.  Moreover, one can show that
		\[
			\wt w_n(x) \sim C_n x^{\alpha_n}e^{\frac1x}, \qquad x\searrow0,
		\]
		with some $C_n>0$, $n\in\bb N$, and $\alpha_n=\frac{3n-2}{2}$ when $n$ is even and
		$\alpha_n=\frac{3n-1}{2}$ when $n$ is odd.
		Hence condition (i) in Definition~\ref{A275} is satisfied, but there is
		no $n\in\bb N$ such that \eqref{A277} holds.
		This potential was also studied in \cite{luger.neuner:2016},
		where it was shown that the approach with super-singular perturbations,
		as developed in \cite{kurasov.luger:2011} and \cite{luger.neuner:2015},
		cannot be applied to this potential.
	\item
		Potentials from the class $H^{-1}_{\rm loc}(0,b)$ could also be treated
		by our method if we relaxed the assumption $V\in L^1_{\rm loc}(0,b)$.
		In this case, one would only have $\phi\in H^1_{\rm loc}(0,b)$.
		Operators with such potentials were considered, e.g.\
		in \cite{hryniv.mykytyuk:2012, savchuk.shkalikov:1999,
		eckhardt.gesztesy.nichols.teschl:2013}.
		Note that the class $H^{-1}_{\rm loc}(0,b)$ includes measure coefficients.
\end{enumerate}
\end{example}

\noindent
Let us introduce the unitary operator
\begin{equation}\label{A323}
	U:\left\{\begin{array}{rcl}
		L^2(0,b) & \to & L^2(w),
		\\[1ex]
		u & \mapsto & \dfrac{u}{\phi}
	\end{array}\right.
\end{equation}
and define the self-adjoint operator
\begin{equation}\label{A267}
	\AV \defeq U^{-1}A_{p,w}U
\end{equation}
with $p$ and $w$ as in \eqref{A276} and $A_{p,w}$ from \eqref{A286}.
For $u\in W^{2,1}(0,b)$ with compact support we have
\begin{equation}\label{A315}
	\begin{aligned}
		\AV u &= -\phi\frac{1}{w}\Bigl(p\Bigl(\frac{u}{\phi}\Bigr)'\Bigr)'
		= -\frac{1}{\phi}\Bigl(\phi^2\frac{\phi u'-\phi'u}{\phi^2}\Bigr)'
		\\[1ex]
		&= -\frac{\phi u''-\phi''u}{\phi} = -u''+Vu.
	\end{aligned}
\end{equation}
Therefore $\AV$ is the Friedrichs extension of the minimal operator
connected with the equation \eqref{A200}; cf.\ the discussion below \eqref{A289}.
In particular, if $\DeltaSchr(V)=1$, then a possible boundary condition at $0$
to characterize $A_V$ is
\[
	\lim_{x\searrow0}\frac{u(x)}{\widetilde w_1(x)} = 0;
\]
see \cite[Theorem~4.3]{niessen.zettl:1992}.
As mentioned above, if $\DeltaSchr(V)\ge2$, then \eqref{A200} is in the limit point
case at $0$ and hence no boundary condition is needed there.

We can apply all theorems from Section~\ref{sec-SL}.
In order to rewrite these results in a more intrinsic form,
we define regularized boundary values by $\rbvSchr{\lambda}u \defeq \rbvSL{\lambda}Uu$
for $\lambda\in\bb C$ and $u$ a solution of \eqref{A200}.
Then Theorem~\ref{A220}, together with a straightforward calculation,
yields the following theorem.

\begin{theorem}[\textbf{Regularized boundary values}]\thlab{A321}
	Let $V\in\KSchr$ with $\dom(V)=(0,b)$, set $\Delta\defeq\DeltaSchr(V)$
	and let $\NSchr{\lambda}$ be the set of all solutions of \eqref{A200}.
	Then, for $x_0\in(0,b)$, the following statements hold.
	\begin{enumerate}[{\rm(i)}]
	\item
		For each $\lambda\in\bb C$ and each solution $u\in\NSchr{\lambda}$
		the boundary value
		\[
			\rbvSchrr{\lambda}u
			= \lim_{x\searrow 0}\bigl(\phi(x)u'(x)-\phi'(x)u(x)\bigr),
		\]
		and the regularized boundary value
		\begin{align*}
			\rbvSchrs{\lambda}u
			&= \lim_{x\searrow 0}\Biggl[
				\,\sum_{k=0}^{\lfloor\frac{\Delta-1}{2}\rfloor}
				\lambda^k\Bigl(\wt w_{2k+1}(x)u'(x)-\wt w_{2k+1}'(x)u(x)\Bigr)
				\\[1ex]
			&\quad +
				\left\{\begin{array}{ll}
					\lambda^{\frac{\Delta}{2}} \wt w_\Delta(x)u(x) & \text{if $\Delta$ is even}  \\[1ex]
					0 & \text{if $\Delta$ is odd}
				\end{array}\right\}
				\\[1ex]
			&\quad + \bigl(\rbvSchrr{\lambda}u\bigr)\biggl(\,
				\sum_{k=\lfloor\frac{\Delta+1}{2}\rfloor}^{\Delta-1}
				\sum_{l=0}^{2k-\Delta}
				(-1)^l \lambda^k \wt w_l(x)\wt w_{2k-l+1}(x)\biggr)\Biggr].
		\end{align*}
		exist.
	\item
		For each $\lambda\in\bb C$ we define
		\[
			\rbvSchr{\lambda}:
			\left\{
			\begin{array}{ccl}
			\NSchr{\lambda} & \to & \bb C^2 \\[1ex]
			u & \mapsto & \bigr(\rbvSchrr{\lambda}u,\;\rbvSchrs{\lambda}u\bigr)^T.
			\end{array}
			\right.
		\]
		Then $\rbvSchr{\lambda}$ is a bijection from $\NSchr{\lambda}$ onto $\bb C^2$.
	\item
		For each $\lambda\in\bb C$ there exists an {\rm(}up to scalar multiples{\rm)} unique
		solution $u\in\NSchr{\lambda}\setminus\{0\}$ such that
		$\lim_{x\searrow 0}\frac{u(x)}{\phi(x)}$ exists.

		This solution is characterized by the property that
		$\int_0^{x_0}\phi^2\big|\bigl(\frac{u}{\phi}\bigr)'\big|^2<\infty$
		and also by the property that
		$\rbvSchrr{\lambda}u=0$ {\rm(}and $u\not\equiv0${\rm)}.

		If $u$ is a solution such that $\lim\limits_{x\searrow 0}\frac{u(x)}{\phi(x)}$
		exists, then
		$\rbvSchrs{\lambda}u = \lim\limits_{x\searrow 0}\frac{u(x)}{\phi(x)}$.
	\end{enumerate}
	The regularized boundary value $\rbvSchrs{\lambda}$ depends on the choice of\, $x_0$
	in the following way.
	\begin{enumerate}[{\rm(i)}]
	\setcounter{enumi}{3}
	\item Let $x_0,\hat x_0\in(0,b)$, and let $\rbvSchr{\lambda}$ and $\rbvSchrpr{\lambda}$
		be the correspondingly defined regularized boundary value mappings.
		Then there exists a polynomial $\pxxh(z)$ with real coefficients whose degree
		does not exceed $\Delta-1$ such that
		\[
			\rbvSchrspr{\lambda}u = \rbvSchrs{\lambda}u
			+ \pxxh(\lambda)\rbvSchrr{\lambda}u,\qquad u\in\NSchr{\lambda},\;\; \lambda\in\bb C.
		\]
		Moreover, clearly, $\rbvSchrrpr{\lambda} = \rbvSchrr{\lambda}$.
	\end{enumerate}
\end{theorem}

\medskip

\noindent
The next theorem about a fundamental system of solutions of \eqref{A200} and
the existence of a singular Titchmarsh--Weyl coefficient follows from
Theorem~\ref{A223}\,(i) with the help of the unitary operator $U$ from \eqref{A323}.

\begin{theorem}[\textbf{Singular Titchmarsh--Weyl coefficients}]\thlab{A322}
	Let $V\in\KSchr$ with $\dom(V)=(0,b)$ be given.
	Then, for each fixed $x_0\in(0,b)$, the following statements hold.
	\begin{enumerate}[{\rm(i)}]
	\item
		For $\lambda\in\bb C$ let $\wt\theta(\cdot\,;\lambda)$ and $\wt\varphi(\cdot\,;\lambda)$
		be the unique solutions of \eqref{A200} such that
		\[
			\rbvSchr{\lambda}\wt\theta(\cdot\,;\lambda)=\binom{1}{0}, \qquad
			\rbvSchr{\lambda}\wt\varphi(\cdot\,;\lambda)=\binom{0}{1}.
		\]
		Then, for each $x\in(0,b)$, the functions $\wt\theta(x;\cdot)$ and $\wt\varphi(x;\cdot)$
		are entire of order $\frac12$ and finite type $x$.
		Moreover, for each $\lambda\in\bb C$, one has
		$W\bigl(\wt\varphi(\cdot\,;\lambda),\wt\theta(\cdot\,;\lambda)\bigr)\equiv1$
		where $W\defeq W_1$ denotes the Wronskian as in \eqref{A317} with $p\equiv1$, and the
		following relations hold:
		\begin{alignat*}{2}
			&\lim_{x\searrow 0}\frac{\wt\varphi(x;\lambda)}{\phi(x)} = 1,
			\qquad &
			&\lim_{x\searrow 0}\frac{\phi(x)\wt\varphi'(x;\lambda)-\phi'(x)\wt\varphi(x;\lambda)}{\int_a^x \bigl(\phi(t)\bigr)^2\rd t} = -\lambda,
			\\[1ex]
			&\lim_{x\searrow 0}\frac{\wt\theta(x;\lambda)}{\wt w_1(x)} = -1,
			\hspace*{7ex} &
			&\lim_{x\searrow 0}\bigl(\phi(x)\wt\theta'(x;\lambda)-\phi'(x)\wt\theta(x;\lambda)\bigr) = 1.
		\end{alignat*}
		Further, one has $\wt\varphi(x;0)=\phi(x)$ and $\wt\theta(x;0)=-\wt w_1(x)$, $x\in(0,b)$.
	\item
		The limit
		\[
			\wt m_V(\lambda) \defeq \lim_{x\nearrow b}\frac{\wt\theta(x;\lambda)}{\wt\varphi(x;\lambda)},
			\qquad \lambda\in\bb C\setminus[0,\infty),
		\]
		exists locally uniformly on $\bb C\setminus[0,\infty)$ and defines an
		analytic function in~$\lambda$.
		The function $\wt m_V$ belongs to the class $\mc N_\kappa^{(\infty)}$ with
		$\kappa=\bigl\lfloor\frac{\DeltaSchr(V)}{2}\bigr\rfloor$.
	\item
		We have
		\[
			\wt\theta(\cdot\,;\lambda)-\wt m_V(\lambda)\wt\varphi(\cdot\,;\lambda) \in L^2(x_0,b),
			\qquad \lambda\in\bb C\setminus[0,\infty),
		\]
		and this property characterizes the value $\wt m_V(\lambda)$
		for each $\lambda\in\bb C\setminus[0,\infty)$.
	\item
		For $\lambda\in\bb C\setminus[0,\infty)$ let $u$ be any non-trivial
		solution of \eqref{A200} such that $u|_{(x_0,b)}\in L^2(x_0,b)$.
		Then
		\[
			\wt m_V(\lambda) = -\frac{\rbvSchrs{\lambda}u}{\rbvSchrr{\lambda}u}\,.
		\]
	\end{enumerate}
	The function $\wt m_V$ depends on the choice of $x_0$.
	This dependence is controlled as follows.
	\begin{enumerate}[{\rm(i)}]
	\setcounter{enumi}{4}
	\item
		Let $\hat x_0\in(a,b)$, and let $\wh{\wt m}_V$ be
		the correspondingly defined singular Titchmarsh--Weyl coefficient.
		Then there exists a polynomial $\pxxh$ with real coefficients whose degree
		does not exceed $\DeltaSchr(V)-1$ such that
		\[
			\wh{\wt m}_V(\lambda) = \wt m_V(\lambda)-\pxxh(\lambda).
		\]
	\end{enumerate}
\end{theorem}

\medskip

\noindent
The functions $\wt\theta$ and $\wt\varphi$ are related to
the functions $\theta$ and $\varphi$ corresponding to \eqref{A208} with $p$ and $w$
as in \eqref{A276} as follows:
\begin{equation}\label{A304}
	\wt\theta(x;\lambda) = \phi(x)\theta(x;\lambda), \qquad
	\wt\varphi(x;\lambda) = \phi(x)\varphi(x;\lambda), \qquad
	x\in(0,b).
\end{equation}
The function $\wt m_V$ is called \emph{singular Titchmarsh--Weyl coefficient}.
It follows from \eqref{A225} and \eqref{A304} that $\wt m_V=m_{p,w}$.
As in Section~\ref{sec-SL} one defines equivalence classes $[\wt m]_V$
with respect to the equivalence relation $\hat\sim$ defined in \eqref{A287}.

The next theorem about the existence of a spectral measure follows from Theorem~\ref{A227}.
For the definition of the class $\bb M^-$ see Definition~\ref{A217}.

\begin{theorem}[\textbf{The spectral measure}]\thlab{A324}
	\rule{0ex}{1ex} \\
	Let $V\in\KSchr$ with $\dom(V)=(0,b)$ be given.
	Then there exists a unique Borel measure $\wt\mu_V$ that satisfies
	\[
		\wt\mu_V\bigl([s_1,s_2]\bigr)=\frac{1}{\pi}\lim_{\eps\searrow 0}
		\lim_{\delta\searrow 0}\intop_{s_1-\eps}^{s_2+\eps}
		\Im \wt m_V(t+i\delta)\,\rd t,\quad -\infty<s_1<s_2<\infty,
	\]
	where $\wt m_V\in[\wt m]_V$ is any singular Titchmarsh--Weyl coefficient associated
	with \eqref{A200}.
	We have $\wt\mu_V\in\bb M^-$ and $\Delta^-(\wt\mu_V)=\DeltaSchr(V)$.

	Moreover, $\wt\mu_V(\{0\})>0$ if and only if
	\begin{equation}\label{A325}
		\int_0^b \bigl(\phi(x)\bigr)^2\rd x < \infty.
	\end{equation}
	If \eqref{A325} is satisfied, then
	\[
		\wt\mu_V\bigl(\{0\}\bigr) = \biggl[\,\int_0^b \bigl(\phi(x)\bigr)^2\rd x\biggr]^{-1}.
	\]
\end{theorem}

\medskip

\noindent
Clearly, we have $\wt\mu_V=\mu_{p,w}$ where $p$ and $w$ are as in \eqref{A276}.

\begin{example}\thlab{A327}
Consider $V$ as in Example~\ref{A266}\,(i).
Since $\DeltaSchr(V)=\bigl\lfloor l+\frac{3}{2}\bigr\rfloor$,
we obtain from Theorem~\ref{A322} that $\wt m_V\in\mc N_\kappa^{(\infty)}$
with $\kappa=\bigl\lfloor\frac{l}{2}+\frac{3}{4}\bigr\rfloor$.
Moreover, Theorem~\ref{A324} yields that $\wt\mu_V\in\bb M^-$
with $\Delta^-(\wt\mu_V)=\DeltaSchr(V)=\bigl\lfloor l+\frac{3}{2}\bigr\rfloor$.
\end{example}

In the next theorem we consider the corresponding Fourier transform
and its inverse.  This theorem follows from Theorem~\ref{A229}.

\begin{theorem}[\textbf{The Fourier transform}]\thlab{A326} \rule{0ex}{1ex}\\
	Let $V\in\KSchr$ with $\dom(V)=(0,b)$ be given,
	and let $\wt\mu_V$ be the spectral measure associated with \eqref{A200} as
	in Theorem~\ref{A324}.
	Then the following statements hold.
	\begin{enumerate}[{\rm(i)}]
	\item
		The map defined by
		\begin{equation}\label{A278}
		\begin{aligned}
			\bigl(\wt\Theta_{V}f\bigr)(t) \defeq \int_0^b \wt\varphi(x;t)f(x)\,\rd x,
			\qquad t\in\bb R, \hspace*{20ex} &
			\\[-1ex]
			f\in L^2(0,b),\;\; \sup(\supp f)<b, &
		\end{aligned}
		\end{equation}
		extends to an isometric isomorphism from $L^2(0,b)$ onto $L^2(\wt\mu_V)$
	\item
		The operator $\wt\Theta_V$ establishes a unitary equivalence between $A_V$ and
		the operator $M_{\wt\mu_V}$ of multiplication by the independent variable
		in $L^2(\wt\mu_V)$, i.e.\ we have
		\[
			\wt\Theta_V A_V = M_{\wt\mu_V}\wt\Theta_V.
		\]
	\item
		For compactly supported functions, the inverse of $\wt\Theta_V$ acts as an
		integral transformation, namely,
		\begin{align*}
			\bigl(\wt\Theta_V^{-1}g\bigr)(x) = \int_0^\infty \wt\varphi(x;t)g(t)\,\rd\wt\mu_V(t),
			\qquad x\in(a,b), \hspace*{12ex} &
			\\
			g\in L^2(\wt\mu_V),\;\; \supp g\text{ compact}. &
		\end{align*}
	\end{enumerate}
\end{theorem}

\medskip

\noindent
The existence of a Fourier transform into a scalar $L^2$-space shows in particular
that the spectrum of $A_V$ is simple.

\begin{remark}\thlab{A204}
	Recall that the solution $\phi$ is not unique.  If $\phi$ is multiplied
	by a positive constant $r$, then $\wt w_l$ and $\wt\theta$ are divided by $r$,
	$\wt\varphi$ is multiplied by $r$, and $\wt m_V$ and $\wt\mu_V$ are divided by $r^2$.
	However, in the situation of Example~\ref{A266}\,{\rm(i)} one can normalize $\phi$
	such that \eqref{A291} holds.
\end{remark}

\medskip

\noindent
Finally, let us state global and local uniqueness theorems.
For the case of Bessel-type potentials as in Example~\ref{A266}\,(i)
see, e.g.\ \cite[Theorem~5.1]{eckhardt:2014}.

\begin{theorem}[\textbf{Global Uniqueness Theorem}]\thlab{A279} \rule{0ex}{1ex} \\
	Let $V_1,V_2\in\KSchr$ be given with $\dom(V_i)=(0,b_i)$, $i=1,2$.
	Then the following statements are equivalent:
	\begin{enumerate}[{\rm(i)}]
	\item
		$b_1=b_2$ and $V_1(x)=V_2(x)$, $x\in(0,b_1)$ a.e.;
	\item
		there exists a $c>0$ such that $[\wt m]_{V_1}=c[\wt m]_{V_2}$;
	\item
		there exists a $c>0$ such that $\wt\mu_{V_1}=c\wt\mu_{V_2}$.
	\end{enumerate}
\end{theorem}

\begin{proof}
	For the implication (i)\,$\Rightarrow$\,(ii) see Remark~\ref{A204}.
	The equivalence of (ii) and (iii) is clear from the definition of $\wt\mu_V$.
	Now suppose that (ii) holds and let $p_i=w_i=\phi_i^2$ be as in \eqref{A276}.
	By rescaling $\phi_2$ we may assume that $c=1$.  Then we have $[m]_{p_1,w_1}=[m]_{p_2,w_2}$.
	It follows from Theorem~\ref{A232} that there
	exists $\gamma:(0,b_2)\to(0,b_1)$ such that \eqref{A233} holds.
	However, this implies that $\gamma'(x)=1$ a.e., and hence
	$b_1=b_2$ and $\phi_1=\phi_2$.  This shows that $V_1=V_2$, i.e.\ (i) is satisfied.
\end{proof}

Local uniqueness theorems for Schr\"odinger equations have attracted a
lot of attention recently.
For the case of a regular left endpoint B.~Simon proved the first version of
such a theorem in \cite[Theorem~1.2]{simon:1999}; alternative proofs were given
in \cite{gesztesy.simon:2000_CMP}, \cite{bennewitz:2001} and \cite{langer:2016}.
For Bessel-type operators with potentials as in Example~\ref{A266}\,(i)
a local uniqueness theorem was proved in \cite[Theorem~4.1]{kostenko.teschl:2013}.

\begin{theorem}[\textbf{Local Uniqueness Theorem}]\thlab{A280} \rule{0ex}{1ex} \\
	Let $V_1,V_2\in\KSchr$ be given with $\dom(V_i)=(0,b_i)$, $i=1,2$.
	Then, for $\tau>0$, the following statements are equivalent:
	\begin{enumerate}[{\rm(i)}]
	\item
		one has $V_1(x)=V_2(x)$,\; $x\in\bigl(0,\min\{\tau,b_1,b_2\}\bigr)$ a.e.;
	\item
		there exist singular Titchmarsh--Weyl coefficients $\wt m_{V_1}$ and $\wt m_{V_2}$
		and there exist $c>0$ and $\beta\in(0,2\pi)$ such that, for each $\eps>0$,
		\[
			\wt m_{V_1}\bigl(re^{i\beta}\bigr)-c\wt m_{V_2}\bigl(re^{i\beta}\bigr)
			= \rmO\Bigl(e^{(-2\tau+\eps)\sqrt{r}\sin\frac{\beta}{2}}\Bigr),
			\qquad r\to\infty;
		\]
	\item
		there exist singular Titchmarsh--Weyl coefficients $\wt m_{V_1}$ and $\wt m_{V_2}$
		and there exist $c>0$ and $k\ge0$ such that, for each $\delta\in(0,\pi)$,
		\begin{multline*}
			\wt m_{V_1}(\lambda)-c\wt m_{V_2}(\lambda)
			= \rmO\Bigl(|\lambda|^k e^{-2\tau\Im\sqrt{\lambda}}\Bigr),
			\\[1ex]
			|\lambda|\to\infty,\;\;\lambda\in\bigl\{z\in\bb C:\delta\le\arg z\le2\pi-\delta\bigr\},
		\end{multline*}
		where $\sqrt{\lambda}$ is chosen so that $\Im\sqrt{\lambda}>0$.
	\end{enumerate}
\end{theorem}

\begin{proof}
	This theorem follows from Theorem~\ref{A249}; we only have to observe
	that $s_i(\tau)=\min\{\tau,b_i\}$
	and that the validity of \eqref{A233} with $p_i=w_i$ implies that $\gamma(x)=x$
	for $x\in(0,\tau)$.
\end{proof}

\medskip

\noindent
Let us conclude this section with two remarks about possible extensions.

\begin{remark}\thlab{A306}
	With a similar method one can also treat general Sturm--Liouville equations
	of the form
	\begin{equation}\label{A314}
		-\bigl(Py'\bigr)'+Qy = \lambda Wy
	\end{equation}
	where $1/P$, $Q$ and $W$ are locally integrable.
	If a positive solution $\phi$ of \eqref{A314} with $\lambda=0$ exists such that
	$\phi\in L^2(W|_{(a,x_0)})$, then one can use the mapping $u\mapsto \frac{u}{\phi}$
	to transform \eqref{A314} to an equation of the form \eqref{A208} with
	\[
		p \defeq P\phi^2, \qquad w \defeq W\phi^2;
	\]
	cf.\ \cite[Lemma~3.2]{niessen.zettl:1992}.
	Using Theorem~\ref{A232} one can show that if the spectral measures corresponding
	two equations of the form \eqref{A314} coincide, then the coefficients
	are related via a Liouville transform;
	see \cite[Theorem~3.4]{eckhardt.gesztesy.nichols.teschl:2013}
	for a related result and \cite[Theorem~4.2]{bennewitz:2003} for the case
	when $W\equiv1$ and the left endpoint is regular.
\end{remark}

\begin{remark}\thlab{A313}
One can also apply the results of the first part of the paper to Dirac systems
of the form
\begin{equation}\label{A300}
	-Ju'+Vu = zu
\end{equation}
on an interval $(a,b)$, where $V$ is a real-valued, symmetric
and locally integrable $2\times2$-matrix function, $z\in\bb C$
is the spectral parameter and $u$ is a $2$-vector function.
Assume that there exists a solution $\phi$ of \eqref{A300} with $z=0$
(i.e.\ $J\phi'=V\phi$) which is in $(L^2(a,x_0))^2$ for some $x_0\in(a,b)$.
Under this assumption we can transform \eqref{A300} into a canonical system
\eqref{A30} as it was done in \cite[Section~4.1, pp.~336, 337]{krein.langer:1985}.
Let $\Phi$ be a $2\times2$-matrix solution of $J\Phi'=V\Phi$
(i.e.\ columns of $\Phi$ are solutions of \eqref{A300} with $z=0$) such that
\[
	\binom{\Phi_{12}}{\Phi_{22}} = \phi \qquad\text{and}\qquad
	\det\Phi(x_0) = 1.
\]
From the second relation it follows that $\Phi(x_0)^TJ\Phi(x_0) = J$.
Since $\frac{\rd}{\rd x}(\Phi^TJ\Phi)=0$, we have
\begin{equation}\label{A301}
	\Phi^TJ\Phi = J \qquad\text{on}\;\;(a,b),
\end{equation}
and hence $\det\Phi(x)=1$, $x\in(a,b)$.
Set
\[
	H \defeq \Phi^T\Phi,
\]
which is clearly symmetric and non-negative and does not vanish on
any set of positive measure.
It is easy to see that $y$ is a solution of \eqref{A30} if and only
if $u\defeq\Phi y$ is a solution of \eqref{A300}.
Since $h_{22} = \Phi_{12}^2+\Phi_{22}^2 = \phi_1^2+\phi_2^2$,
condition~(I) in Definition~\ref{A1} is satisfied.
If $H\in\bb H$, i.e.\ also (HS) and ($\Delta$) are fulfilled,
then one can apply the results from Sections~\ref{sec-measure}--\ref{sec-inverse}.

In order to write the results in a more intrinsic form, one can
use the unitary transformation
\[
	U: \begin{cases}
		L^2(H) \to (L^2(a,b))^2, \\[0.5ex]
		y \mapsto \Phi u,
	\end{cases}
\]
whose inverse acts like $U^{-1}u = \Phi^{-1}u = -J\Phi^TJu$.
For instance, one can define regularized boundary values by
$\rbvDirac{z}u \defeq \rbv U^{-1}u$ as in Section~\ref{sec-Schroedinger}.
Details are left to the reader.
See also, e.g.\ \cite{brunnhuber.eckhardt.kostenko.teschl:2014,
fritzsche.kirstein.sakhnovich:2012} for different approaches to Dirac operators.
\end{remark}



\begin{thebibliography}{999}

\bibitem{albeverio.gesztesy.hoeghkrohn.holden:2005}
S.~Albeverio, F.~Gesztesy, R.~H{\o}egh-Krohn and H.~Holden,
\textit{Solvable Models in Quantum Mechanics}. Second edition.
With an appendix by Pavel Exner.
AMS Chelsea Publishing, Providence, RI, 2005.

\bibitem{albeverio.hryniv.mykytyuk:2005}
S.~Albeverio, R.~Hryniv and Ya.~Mykytyuk,
Inverse spectral problems for Sturm--Liouville operators in impedance form.
\textit{J.\ Funct.\ Anal.} \textbf{222} (2005), 143--177.

\bibitem{albeverio.hryniv.mykytyuk:2012}
S.~Albeverio, R.~Hryniv and Ya.~Mykytyuk,
Scattering theory for Schr\"odinger operators with Bessel-type potentials.
\textit{J.\ Reine Angew.\ Math.} \textbf{666} (2012), 83--113.

\bibitem{arnold:1989}
V.I.~Arnol{\cprime}d,
\textit{\cyr Matematicheskie metody klassicheskoi mekhaniki} [Russian].
Third edition. `Nauka', Moscow, third edition edition, 1989.
English translation: \textit{Mathematical Methods of Classical Mechanics},
Graduate Texts in Mathematics, vol.~60, Springer-Verlag, New York, 1989.

\bibitem{arov.dym:2008}
D.Z.~Arov and H.~Dym,
\textit{{$J$}-Contractive Matrix Valued Functions and Related Topics}.
\textit{Encyclopedia of Mathematics and its Applications}, vol.~116.
Cambridge University Press, Cambridge, 2008.

\bibitem{atkinson:1964}
F.V.~Atkinson,
\textit{Discrete and Continuous Boundary Problems}.
Mathematics in Science and Engineering, vol.~8.
Academic Press, New York, 1964.

\bibitem{bennewitz:2001}
C.~Bennewitz,
A proof of the local Borg--Marchenko theorem.
\textit{Comm.\ Math.\ Phys.} \textbf{218} (2001), 131--132.

\bibitem{bennewitz:2003}
C.~Bennewitz,
A Paley--Wiener theorem with applications to inverse spectral theory.
In: \textit{Advances in Differential Equations and Mathematical Physics (Birmingham, AL, 2002)}.
Contemp.\ Math., vol.~327, Amer.\ Math.\ Soc., Providence, RI, 2003, pp.~21--31.

\bibitem{bennewitz.everitt:2005}
C.~Bennewitz and W.N.~Everitt,
The Titchmarsh--Weyl eigenfunction expansion theorem for
Sturm--Liouville differential equations.
In: \textit{Sturm--Liouville Theory}.
Birkh\"auser, Basel, 2005, pp.~137--171.

\bibitem{bognar:1974}
J.~Bogn{\'a}r,
\textit{Indefinite Inner Product Spaces}.
Ergebnisse der Mathematik und ihrer Grenzgebiete, vol.~78.
Springer-Verlag, New York, 1974.

\bibitem{debranges:1960}
L.~de~Branges,
Some Hilbert spaces of entire functions.
\textit{Trans.\ Amer.\ Math.\ Soc.} \textbf{96} (1960), 259--295.

\bibitem{debranges:1961}
L.~de~Branges,
Some Hilbert spaces of entire functions. II.
\textit{Trans.\ Amer.\ Math.\ Soc.} \textbf{99} (1961), 118--152.

\bibitem{debranges:1961a}
L.~de~Branges,
Some Hilbert spaces of entire functions. III.
\textit{Trans.\ Amer.\ Math.\ Soc.} \textbf{100} (1961), 73--115.

\bibitem{debranges:1962a}
L.~de~Branges,
Some Hilbert spaces of entire functions. IV.
\textit{Trans.\ Amer.\ Math.\ Soc.} \textbf{105} (1962), 43--83.

\bibitem{debranges:1968}
L.~de~Branges,
\textit{Hilbert Spaces of Entire Functions}.
Prentice-Hall Inc., Englewood Cliffs, N.J., 1968.

\bibitem{brunnhuber.eckhardt.kostenko.teschl:2014}
R.~Brunnhuber, J.~Eckhardt, A.~Kostenko and G.~Teschl,
Singular Weyl--Titchmarsh--Kodaira theory for one-dimensional Dirac operators.
\textit{Monatsh.\ Math.} \textbf{174} (2014), 515--547.

\bibitem{bube.burridge:1983}
K.P.~Bube and R.~Burridge,
The one-dimensional inverse problem of reflection seismology.
\textit{SIAM Rev.} \textbf{25} (1983), 497--559.

\bibitem{daho.langer:1985}
K.~Daho and H.~Langer,
Matrix functions of the class $N_\kappa$.
\textit{Math.\ Nachr.} \textbf{120} (1985), 275--294.

\bibitem{davies:2013}
E.B.~Davies,
Singular Schr\"odinger operators in one dimension.
\textit{Mathematika} \textbf{59} (2013), 141--159.

\bibitem{derkach:1998}
V.A.~Derkach,
On extensions of the Laguerre operator in spaces with an indefinite metric.
\textit{Mat.\ Zametki} \textbf{63} (1998), 509--521.

\bibitem{dijksma.snoo:1987b}
A.~Dijksma and H.~de~Snoo,
Symmetric and selfadjoint relations in Kre\u\i n spaces. II.
\textit{Ann.\ Acad.\ Sci.\ Fenn.\ Ser.\ A I Math.} \textbf{12} (1987), 199--216.

\bibitem{dijksma.langer.luger.shondin:2000}
A.~Dijksma, H.~Langer, A.~Luger and Y.~Shondin,
A factorization result for generalized Nevanlinna functions of the
class $\mathcal N_\kappa$.
\textit{Integral Equations Operator Theory} \textbf{36} (2000), 121--125.

\bibitem{dijksma.langer.shondin:2004}
A.~Dijksma, H.~Langer and Y.~Shondin,
Rank one perturbations at infinite coupling in Pontryagin spaces.
\textit{J.\ Funct.\ Anal.} \textbf{209} (2004), 206--246.

\bibitem{dijksma.langer.shondin.zeinstra:2000}
A.~Dijksma, H.~Langer, Y.~Shondin and C.~Zeinstra,
Self-adjoint operators with inner singularities and Pontryagin spaces.
In: \textit{Operator Theory and Related Topics, vol.~II (Odessa, 1997)}.
Oper.\ Theory Adv.\ Appl., vol.~118.
Birkh\"auser, Basel, 2000, pp.~105--175.

\bibitem{dijksma.luger.shondin:2006}
A.~Dijksma, A.~Luger and Y.~Shondin,
Minimal models for $\mathcal N_\kappa^\infty$-functions.
In: \textit{Operator Theory and Indefinite Inner Product Spaces}.
Oper.\ Theory Adv.\ Appl., vol.~163.
Birkh\"auser, Basel, 2006, pp.~97--134.

\bibitem{dijksma.luger.shondin:2009}
A.~Dijksma, A.~Luger and Y.~Shondin,
Approximation of $\mathcal N_\kappa^\infty$-functions. I:
Models and regularization.
In: \textit{Spectral Theory in Inner Product Spaces and Applications}.
Oper.\ Theory Adv.\ Appl., vol.~188.
Birkh\"auser Verlag, Basel, 2009, pp.~87--112.

\bibitem{dijksma.luger.shondin:2010}
A.~Dijksma, A.~Luger and Y.~Shondin,
Approximation of $\mathcal N_\kappa^\infty$-functions II:
Convergence of models.
In: \textit{Recent Advances in Operator Theory in Hilbert and Krein Spaces}.
Oper.\ Theory Adv.\ Appl., vol.~198.
Birkh\"auser Verlag, Basel, 2010, pp.~125--169.

\bibitem{dijksma.shondin:2000}
A.~Dijksma and Y.~Shondin,
Singular point-like perturbations of the Bessel operator in a Pontryagin space.
\textit{J.\ Differential Equations} \textbf{164} (2000), 49--91.

\bibitem{dijksma.shondin:2002}
A.~Dijksma and Y.~Shondin,
Singular point-like perturbations of the Laguerre operator in a Pontryagin space.
In: \textit{Operator Methods in Ordinary and
Partial Differential Equations (Stockholm, 2000)}.
Oper.\ Theory Adv.\ Appl., vol.~132.
Birkh\"auser, Basel, 2002, pp.~141--181.

\bibitem{eckhardt:2014}
J.~Eckhardt,
Inverse uniqueness results for Schr\"odinger operators using de Branges theory.
\textit{Complex Anal.\ Oper.\ Theory} \textbf{8} (2014), 37--50.

\bibitem{eckhardt.gesztesy.nichols.teschl:2013}
J.~Eckhardt, F.~Gesztesy, R.~Nichols and G.~Teschl,
Inverse spectral theory for Sturm--Liouville operators with distributional potentials.
\textit{J.\ Lond.\ Math.\ Soc.\ (2)} \textbf{88} (2013), 801--828.

\bibitem{eckhardt.teschl:2013}
J.~Eckhardt and G.~Teschl,
Uniqueness results for one-dimensional Schr\"odinger operators with
purely discrete spectra.
\textit{Trans.\ Amer.\ Math.\ Soc.} \textbf{365} (2013), 3923--3942.

\bibitem{flanders:1989}
H.~Flanders,
\textit{Differential Forms with Applications to the Physical Sciences}.
Second edition.
Dover Books on Advanced Mathematics. Dover Publications Inc.,
New York, 1989.

\bibitem{fritzsche.kirstein.sakhnovich:2012}
B.~Fritzsche, B.~Kirstein and A.L.~Sakhnovich,
Weyl functions of generalized Dirac systems: integral representation,
the inverse problem and discrete interpolation.
\textit{J.\ Anal.\ Math.} \textbf{116} (2012), 17--51.

\bibitem{fulton:2008}
C.~Fulton,
Titchmarsh--Weyl $m$-functions for second-order Sturm--Liouville problems
with two singular endpoints.
\textit{Math.\ Nachr.} \textbf{281} (2008), 1418--1475.

\bibitem{fulton.langer:2010}
C.~Fulton and H.~Langer,
Sturm--Liouville operators with singularities and generalized Nevanlinna functions.
\textit{Complex Anal.\ Oper.\ Theory} \textbf{4} (2010), 179--243.

\bibitem{fulton.langer.luger:2012}
C.~Fulton, H.~Langer and A.~Luger,
Mark Krein's method of directing functionals and singular potentials.
\textit{Math.\ Nachr.} \textbf{285} (2012), 1791--1798.

\bibitem{gesztesy.simon:2000_CMP}
F.~Gesztesy and B.~Simon,
On local Borg--Marchenko uniqueness results.
\textit{Comm.\ Math.\ Phys.} \textbf{211} (2000), 273--287.

\bibitem{gesztesy.zinchenko:2006}
F.~Gesztesy and M.~Zinchenko,
On spectral theory for Schr\"odinger operators with strongly singular potentials.
\textit{Math.\ Nachr.} \textbf{279} (2006), 1041--1082.

\bibitem{gilbert:1998}
D.J.~Gilbert,
On subordinacy and spectral multiplicity for a class of singular differential operators.
\textit{Proc.\ Roy.\ Soc.\ Edinburgh Sect.\ A} \textbf{128} (1998), 549--584.

\bibitem{gohberg.krein:1967}
I.C.~Gohberg and M.G.~Krein,
\textit{{\cyr Teoriya vol{\cprime}terrovykh operatorov v gil{\cprime}bertovom
prostranstve i ee prilozheniya}} [Russian].
Izdat.\ `Nauka', Moscow, 1967.
English translation:
\textit{Theory and Applications of Volterra Operators in Hilbert Space},
Translations of Mathematical Monographs, vol.~24.
American Mathematical Society, Providence, R.I., 1970.

\bibitem{hassi.snoo.winkler:2000}
S.~Hassi, H.~de~Snoo and H.~Winkler,
Boundary-value problems for two-dimensional canonical systems.
\textit{Integral Equations Operator Theory} \textbf{36} (2000), 445--479.

\bibitem{hassi.remling.snoo:2000}
S.~Hassi, C.~Remling and H.~de~Snoo,
Subordinate solutions and spectral measures of canonical systems.
\textit{Integral Equations Operator Theory} \textbf{37} (2000), 48--63.

\bibitem{herglotz:1911}
G.~Herglotz,
\"Uber Potenzreihen mit positivem, reellem Teil im Einheitskreis [German].
\textit{Leipz.\ Ber.} \textbf{63} (1911), 501--511.

\bibitem{hoermander:1990}
L.~H{\"o}rmander,
\textit{The Analysis of Linear Partial Differential Operators. I.
Distribution Theory and Fourier Analysis.}
Second edition.
Grundlehren der Mathematischen Wissenschaften [Fundamental
Principles of Mathematical Sciences], vol.~256.
Springer-Verlag, Berlin, 1990.

\bibitem{hryniv.mykytyuk:2012}
R.O.~Hryniv and Y.V.~Mykytyuk,
Self-adjointness of Schr\"odinger operators with singular potentials.
\textit{Methods Funct.\ Anal.\ Topology} \textbf{18} (2012), 152--159.

\bibitem{hryniv.sacks:2010}
R.~Hryniv and P.~Sacks,
Numerical solution of the inverse spectral problem for Bessel operators.
\textit{J.\ Comput.\ Appl.\ Math.} \textbf{235} (2010), 120--136.

\bibitem{jonas.langer.textorius:1992}
P.~Jonas, H.~Langer and B.~Textorius,
Models and unitary equivalence of cyclic selfadjoint operators in Pontrjagin spaces.
In: \textit{Operator Theory and Complex Analysis (Sapporo, 1991)}.
Oper. Theory Adv. Appl., vol.~59.
Birkh\"auser, Basel, 1992, pp.~252--284.

\bibitem{kac:1950}
I.~Kac,
On Hilbert spaces generated by monotone Hermitian matrix-functions.
\textit{Har\cprime kov Gos.\ Univ.\ U\v c.\ Zap.\ 34 = Zap.\ Mat.\ Otd.\
Fiz.-Mat.\ Fak.\ i Har\cprime kov.\ Mat.\ Ob\v s\v c.\ (4)} \textbf{22} (1951), 95--113.
  1950.

\bibitem{kac:1962a}
I.S.~Kac,
On the spectral multiplicity of a second-order differential operator.
\textit{Dokl.\ Akad.\ Nauk SSSR} \textbf{145} (1962), 510--513.

\bibitem{kac:1963a}
I.S.~Kac,
Spectral multiplicity of a second-order differential operator and
expansion in eigenfunction.
\textit{Izv.\ Akad.\ Nauk SSSR Ser.\ Mat.} \textbf{27} (1963), 1081--1112.

\bibitem{kac:1984}
I.S.~Kac,
Linear relations, generated by a canonical differential equation on
an interval with a regular endpoint, and expansibility in eigenfunctions.
\textit{VINITI Deponirovannye Nauchnye Raboty} \textbf{195} (1985), 50 pp., b.o.~720,
Deposited in Ukr NIINTI, no.~1453, 1984.

\bibitem{kac:1999}
I.S.~Kac,
Inclusion of the Hamburger power moment problem in the spectral
theory of canonical systems.
\textit{Zap.\ Nauchn.\ Sem.\ S.-Peterburg.\ Otdel.\ Mat.\ Inst.\ Steklov.\
(POMI)} \textbf{262} (1999)
(Issled.\ po Linein.\ Oper.\ i Teor.\ Funkts.\ \textbf{27}), 147--171, 234.
English translation: \textit{J.\ Math.\ Sci.\ (New York)} \textbf{110} (2002), 2991--3004.

\bibitem{kac.krein:1968}
I.S.~Kac and M.G.~Krein,
\textit{{On spectral functions of a string}}, pp.~648--737.
Izdat.\ `Mir', Moscow, 1968.
Addition II in F.V.~Atkinson,
{\cyr Diskretnye i nepreryvnye granichnye zadachi} (Russian translation of
`Discrete and Continuous Boundary Problems').
English translation: \textit{Amer.\ Math.\ Soc.\ Transl.\ (2)} \textbf{103} (1974), 19--102.

\bibitem{kaltenbaeck.winkler.woracek:nksym}
M.~Kaltenb\"ack, H.~Winkler and H.~Woracek,
Generalized Nevanlinna functions with essentially positive spectrum.
\textit{J.\ Operator Theory} \textbf{55} (2006), 17--48.

\bibitem{kaltenbaeck.woracek:krall}
M.~Kaltenb\"ack and H.~Woracek,
Generalized resolvent matrices and spaces of analytic functions.
\textit{Integral Equations Operator Theory} \textbf{32} (1998), 282--318.

\bibitem{kaltenbaeck.woracek:p2db}
M.~Kaltenb\"ack and H.~Woracek,
Pontryagin spaces of entire functions. II.
\textit{Integral Equations Operator Theory} \textbf{33} (1999), 305--380.

\bibitem{kaltenbaeck.woracek:p3db}
M.~Kaltenb\"ack and H.~Woracek,
Pontryagin spaces of entire functions. III.
\textit{Acta Sci.\ Math.\ (Szeged)} \textbf{69} (2003), 241--310.

\bibitem{kaltenbaeck.woracek:p4db}
M.~Kaltenb\"ack and H.~Woracek,
Pontryagin spaces of entire functions. IV.
\textit{Acta Sci.\ Math.\ (Szeged)} \textbf{72} (2006), 709--835.

\bibitem{kaltenbaeck.woracek:hskansys}
M.~Kaltenb\"ack and H.~Woracek,
Canonical differential equations of Hilbert--Schmidt type.
In: \textit{Operator Theory in Inner Product Spaces}.
Oper. Theory Adv. Appl., vol.~175.
Birkh\"auser, Basel, 2007, pp.~159--168.

\bibitem{kaltenbaeck.woracek:p6db}
M.~Kaltenb\"ack and H.~Woracek,
Pontryagin spaces of entire functions. VI.
\textit{Acta Sci.\ Math.\ (Szeged)} \textbf{76} (2010), 511--560.

\bibitem{kaltenbaeck.woracek:p5db}
M.~Kaltenb\"ack and H.~Woracek,
Pontryagin spaces of entire functions. V.
\textit{Acta Sci.\ Math.\ (Szeged)} \textbf{77} (2011), 223--336.

\bibitem{kodaira:1949}
K.~Kodaira,
The eigenvalue problem for ordinary differential equations of the
second order and Heisenberg's theory of $S$-matrices.
\textit{Amer.\ J.\ Math.} \textbf{71} (1949), 921--945.

\bibitem{kostenko.sakhnovich.teschl:2010}
A.~Kostenko, A.~Sakhnovich and G.~Teschl,
Inverse eigenvalue problems for perturbed spherical Schr\"odinger operators.
\textit{Inverse Problems} \textbf{26} (2010), 105013, 14~pp.

\bibitem{kostenko.sakhnovich.teschl:2012a}
A.~Kostenko, A.~Sakhnovich and G.~Teschl,
Commutation methods for Schr\"odinger operators with strongly singular potentials.
\textit{Math.\ Nachr.} \textbf{285} (2012), 392--410.

\bibitem{kostenko.sakhnovich.teschl:2012}
A.~Kostenko, A.~Sakhnovich and G.~Teschl,
Weyl--Titchmarsh theory for Schr\"odinger operators with strongly singular potentials.
\textit{Int.\ Math.\ Res.\ Not.\ IMRN} \textbf{2012} (2012), 1699--1747.

\bibitem{kostenko.teschl:2011}
A.~Kostenko and G.~Teschl,
On the singular Weyl--Titchmarsh function of perturbed spherical Schr\"odinger operators.
\textit{J.\ Differential Equations} \textbf{250} (2011), 3701--3739.

\bibitem{kostenko.teschl:2013}
A.~Kostenko and G.~Teschl,
Spectral asymptotics for perturbed spherical Schr\"odinger
operators and applications to quantum scattering.
\textit{Comm.\ Math.\ Phys.} \textbf{322} (2013), 255--275.

\bibitem{krall:1979}
A.M.~Krall,
Laguerre polynomial expansions in indefinite inner product spaces.
\textit{J.\ Math.\ Anal.\ Appl.} \textbf{70} (1979), 267--279.

\bibitem{krall:1982}
A.M.~Krall,
On boundary values for the Laguerre operator in indefinite inner product spaces.
\textit{J.\ Math.\ Anal.\ Appl.} \textbf{85} (1982), 406--408.

\bibitem{krein.langer:1973}
M.G.~Krein and H.~Langer,
\"Uber die $Q$-Funktion eines $\pi$-hermiteschen Operators im Raume $\Pi_\kappa$ [German].
\textit{Acta Sci.\ Math.\ (Szeged)} \textbf{34} (1973), 191--230.

\bibitem{krein.langer:1977}
M.G.~Krein and H.~Langer,
\"Uber einige Fortsetzungsprobleme, die eng mit der Theorie
hermitescher Operatoren im Raume $\Pi_\kappa$ zusammenh\"angen. I.
Einige Funktionenklassen und ihre Darstellungen [German].
\textit{Math.\ Nachr.} \textbf{77} (1977), 187--236.

\bibitem{krein.langer:1985}
M.G.~Krein and H.~Langer,
On some continuation problems which are closely related to the theory
of operators in spaces $\Pi_\kappa$. IV.
Continuous analogues of orthogonal polynomials on the unit circle
with respect to an indefinite weight and related continuation problems
for some classes of functions.
\textit{J.\ Operator Theory} \textbf{13} (1985), 299--417.

\bibitem{krein.langer:2014}
M.G.~Krein and H.~Langer,
Continuation of Hermitian positive definite functions and related questions.
\textit{Integral Equations Operator Theory} \textbf{78} (2014), 1--69.

\bibitem{kurasov.luger:2011}
P.~Kurasov and A.~Luger,
An operator theoretic interpretation of the generalized
Titchmarsh--Weyl coefficient for a singular Sturm--Liouville problem.
\textit{Math.\ Phys.\ Anal.\ Geom.} \textbf{14} (2011), 115--151.

\bibitem{langer:1982}
H.~Langer,
Spectral functions of definitizable operators in Krein spaces.
In: \textit{Functional Analysis (Dubrovnik, 1981)}.
Lecture Notes in Math., vol.~948.
Springer, Berlin, 1982, pp.~1--46.

\bibitem{langer:1986}
H.~Langer,
A characterization of generalized zeros of negative type of functions
of the class $N_\kappa$.
In: \textit{Advances in Invariant Subspaces and other Results of Operator Theory
(Timi\c soara and Herculane, 1984)}.
Oper.\ Theory Adv.\ Appl., vol.~17.
Birkh\"auser, Basel, 1986, pp.~201--212.

\bibitem{langer:2016}
H.~Langer,
Transfer functions and local spectral uniqueness for Sturm--Liouville operators,
canonical systems and strings.
\textit{Integral Equations Operator Theory} \textbf{85} (2016), 1--23.

\bibitem{langer.luger.matsaev:2011}
H.~Langer, A.~Luger and V.~Matsaev,
Convergence of generalized Nevanlinna functions.
\textit{Acta Sci.\ Math.\ (Szeged)} \textbf{77} (2011), 425--437.

\bibitem{langer.woracek:esmod}
M.~Langer and H.~Woracek,
A function space model for canonical systems with an inner singularity.
\textit{Acta Sci.\ Math.\ (Szeged)} \textbf{77} (2011), 101--165.

\bibitem{langer.woracek:lokinv}
M.~Langer and H.~Woracek,
A local inverse spectral theorem for Hamiltonian systems.
\textit{Inverse Problems} \textbf{27} (2011), 055002, 17~pp.

\bibitem{langer.woracek:gpinf}
M.~Langer and H.~Woracek,
Indefinite Hamiltonian systems whose Titchmarsh--Weyl coefficients
have no finite generalized poles of non-positive type.
\textit{Oper.\ Matrices} \textbf{7} (2013), 477--555.

\bibitem{langer.woracek:expty}
M.~Langer and H.~Woracek,
The exponential type of the fundamental solution of an indefinite Hamiltonian system.
\textit{Complex Anal.\ Oper.\ Theory} \textbf{7} (2013), 285--312.

\bibitem{langer.woracek:ninfrep}
M.~Langer and H.~Woracek,
Distributional representations of generalized Nevanlinna functions.
\textit{Math.\ Nachr.} \textbf{288} (2015), 1127--1149.

\bibitem{lifschitz:1989}
A.E.~Lifschitz,
\textit{Magnetohydrodynamics and Spectral Theory}.
Developments in Electromagnetic Theory and Applications, vol.~4.
Kluwer Academic Publishers Group, Dordrecht, 1989.

\bibitem{luger.neuner:2015}
A.~Luger and C.~Neuner,
An operator theoretic interpretation of the generalized
Titchmarsh--Weyl function for perturbed spherical Schr\"odinger operators.
\textit{Complex Anal.\ Oper.\ Theory} \textbf{9} (2015), 1391--1410.

\bibitem{luger.neuner:2016}
A.~Luger and C.~Neuner,
On the Weyl solution of the 1-dim Schr\"odinger operator with inverse fourth power potential.
\textit{Monatsh.\ Math.} \textbf{180} (2016), 295--303.

\bibitem{mclaughlin:1986}
J.R.~McLaughlin,
Analytical methods for recovering coefficients in differential
equations from spectral data.
\textit{SIAM Rev.} \textbf{28} (1986), 53--72.

\bibitem{niessen.zettl:1992}
H.-D.~Niessen and A.~Zettl,
Singular Sturm--Liouville problems: the Friedrichs extension and
comparison of eigenvalues.
\textit{Proc.\ London Math.\ Soc.\ (3)} \textbf{64} (1992), 545--578.

\bibitem{orcutt:1969}
B.C.~Orcutt,
\textit{Canonical Differential Equations}.
ProQuest LLC, Ann Arbor, MI, 1969.
PhD Thesis, University of Virginia.

\bibitem{pick:1915}
G.~Pick,
\"Uber die Beschr\"ankungen analytischer Funktionen, welche durch
vorgegebene Funktionswerte bewirkt werden [German].
\textit{Math.\ Ann.} \textbf{77} (1915), 7--23.

\bibitem{remling:2002}
C.~Remling,
Schr\"odinger operators and de Branges spaces.
\textit{J.\ Funct.\ Anal.} \textbf{196} (2002), 323--394.

\bibitem{remling:2018}
C.~Remling,
Spectral Theory of Canonical Systems.
De Gruyter Studies in Mathematics, vol.~70. De Gruyter, Berlin, 2018.

\bibitem{romanov:1408.6022v1}
R.~Romanov,
Canonical systems and de Branges spaces,
arXiv:1408.6022v1, 2014.

\bibitem{sakhnovich:1999}
L.A.~Sakhnovich,
\textit{Spectral Theory of Canonical Differential Systems.
Method of Operator Identities}.
Translated from the Russian manuscript by E.~Melnichenko.
Oper.\ Theory Adv.\ Appl., vol.~107.
Birkh\"auser Verlag, Basel, 1999.

\bibitem{savchuk.shkalikov:1999}
A.M.~Savchuk and A.A.~Shkalikov,
Sturm--Liouville operators with singular potentials.
\textit{Mat.\ Zametki} \textbf{66} (1999), 897--912.

\bibitem{silva.teschl.toloza:2015}
L.O.~Silva, G.~Teschl and J.H.~Toloza,
Singular Schr\"odinger operators as self-adjoint extensions of $N$-entire operators.
\textit{Proc.\ Amer.\ Math.\ Soc.} \textbf{143} (2015), 2103--2115.

\bibitem{silva.toloza:2014}
L.O.~Silva and J.H.~Toloza,
A class of $n$-entire Schr\"odinger operators.
\textit{Complex Anal.\ Oper.\ Theory} \textbf{8} (2014), 1581--1599.

\bibitem{simon:1999}
B.~Simon,
A new approach to inverse spectral theory. I. Fundamental formalism.
\textit{Ann.\ of Math.\ (2)} \textbf{150} (1999), 1029-1057.

\bibitem{snoo.winkler:2005}
H.~de~Snoo and H.~Winkler,
Canonical systems of differential equations with self-adjoint
interface conditions on graphs.
\textit{Proc.\ Roy.\ Soc.\ Edinburgh Sect.\ A} \textbf{135} (2005), 297--315.

\bibitem{snoo.winkler:2005a}
H.~de~Snoo and H.~Winkler,
Two-dimensional trace-normed canonical systems of differential
equations and selfadjoint interface conditions.
\textit{Integral Equations Operator Theory} \textbf{51} (2005), 73--108.

\bibitem{winkler:1995}
H.~Winkler,
The inverse spectral problem for canonical systems.
\textit{Integral Equations Operator Theory} \textbf{22} (1995), 360--374.

\bibitem{winkler:1995a}
H.~Winkler,
On transformations of canonical systems.
In: \textit{Operator Theory and Boundary Eigenvalue Problems (Vienna, 1993)}.
Oper.\ Theory Adv.\ Appl., vol.~80.
Birkh\"auser, Basel, 1995, pp.~276--288.

\bibitem{winkler.woracek:nnham}
H.~Winkler and H.~Woracek,
Reparametrizations of non trace-normed Hamiltonians.
In: \textit{Spectral Theory, Mathematical System Theory, Evolution Equations,
Differential and Difference equations}.
Oper.\ Theory Adv.\ Appl., vol.~221.
Birkh\"auser/Springer Basel AG, Basel, 2012, pp.~667--690.

\bibitem{winkler.woracek:del}
H.~Winkler and H.~Woracek,
A growth condition for Hamiltonian systems related with Krein strings.
\textit{Acta Sci.\ Math.\ (Szeged)} \textbf{80} (2014), 31--94.

\bibitem{woracek:nass}
H.~Woracek,
Existence of zerofree functions $N$-associated to a de Branges Pontryagin space.
\textit{Monatsh.\ Math.} \textbf{162} (2011), 453--506.

\end{thebibliography}

\def\cprime{$'$}

{\small 
\begin{flushleft}
	M.\,Langer \\
	Department of Mathematics and Statistics \\
	University of Strathclyde \\
	26 Richmond Street \\
	Glasgow G1 1XH \\
	UNITED KINGDOM \\
	email: \texttt{m.langer@strath.ac.uk} \\[5mm]
	H.\,Woracek \\
	Institute for Analysis and Scientific Computing \\
	Vienna University of Technology \\
	Wiedner Hauptstra{\ss}e.\ 8--10/101 \\
	1040 Wien \\
	AUSTRIA \\
	email: \texttt{harald.woracek@tuwien.ac.at}
\end{flushleft}
}

\end{document}